%% file: infcoh.tex
\documentclass[10pt,letterpaper,twoside]{article}

\usepackage[letterpaper,left=1in,right=1in,top=1.5in,bottom=1.5in]{geometry}

\title{Filtrations and cohomology I: crystallization}\author{Benjamin Antieau}\date{\today}

\input{preamble}

\begin{document}

\maketitle
\begin{abstract}
    \noindent
    We compare several different notions of filtered derived commutative ring, discussing
    HKR-filtered Hochschild homology, Hodge-filtered de Rham cohomology, and the lesser-known
    Hodge-filtered infinitesimal cohomology. Our main result is that de Rham cohomology is the crystallization
    of infinitesimal cohomology.
\end{abstract}

\renewcommand{\baselinestretch}{0.75}\small
\tableofcontents
\renewcommand{\baselinestretch}{1.0}\normalsize

\section{Introduction}\label{sec:intro}

This is the first paper in a series on the interaction between filtrations and cohomology. It
introduces the main objects which will be studied in the sequels in the setting of derived
commutative rings, a nonconnective extension of the theory of animated commutative rings, which was discovered by Bhatt and Mathew and explained in Raksit's
thesis~\cite{raksit}.

For $k\rightarrow R$ a map of
derived commutative rings, we have
\begin{itemize}
    \item $\HH_\fil(R/k)$, the HKR-filtered Hochschild homology;
    \item $\F^\star_\H\dR_{R/k}$, the Hodge-filtered (derived) de Rham cohomology;
    \item $\F^\star_\H\Inf_{R/k}$, the Hodge-filtered (derived) infinitesimal cohomology.
\end{itemize}
Each of these admits defining universal properties and de Rham cohomology and infinitesimal
cohomology also admit completions $\F^\star_\H\dRhat_{R/k}$ and
$\F^\star_\H\Infhat_{R/k}$.\footnote{When $k$ and $R$ are connective, $\HH_\fil(R/k)$ is already a
complete (and exhaustive) filtration on $\HH(R/k)$, so there is no need to complete.} HKR-filtered Hochschild homology is the initial
{\em infinitesimal} filtered derived commutative $k$-algebra with a filtered circle action. Hodge-filtered de Rham
cohomology is the initial {\em crystalline} filtered derived commutative $k$-algebra with a map
$R\rightarrow\gr^0$. This generalizes the usual universal property of the de Rham complex
$\Omega^\bullet_{R/k}$ as the initial cdga over $k$ with a map of graded objects
$R\rightarrow\Omega^*_{R/k}$.
Hodge-filtered infinitesimal cohomology is the initial {\em infinitesimal} filtered derived
commutative $k$-algebra with a map $R\rightarrow\gr^0$. The associated graded pieces of these
filtered derived commutative rings are
\begin{itemize}
    \item   $\gr^i\HH(R/k)\we\Lambda^i\L_{R/k}[i]\we\LSym^i_R(\L_{R/k}[1])$;
    \item   $\gr^i\dR_{R/k}\we\Lambda^i\L_{R/k}[-i]\we\LSym^i_R(\L_{R/k}[1])[-2i]$;
    \item   $\gr^i\Inf_{R/k}\we\LSym^i_R(\L_{R/k}[-1])$.
\end{itemize}
Here, $\LSym^i$ is a certain nonabelian derived functor arising in the derived commutative ring
monad. These graded pieces do not change upon completion.

Infinitesimal filtered derived commutative rings are the naive way to `derive' the notion of a
filtered commutative ring. Indeed, the $\infty$-category of infinitesimal filtered derived commutative
rings which are connective with respect to the neutral $t$-structure is the animation of the category of finitely
generated filtered polynomial rings. Crystalline filtered derived commutative rings are trickier to
define, but their graded cousins are obtained by shearing (in the sense of
Section~\ref{sec:grdalg}) the naive notion of graded derived
commutative ring. The difference between infinitesimal and crystalline derived commutative rings
vanishes over $\bQ$. We will see in F\&C.II that many of the properties enjoyed by de
Rham cohomology in characteristic $0$ extend to infinitesimal cohomology over $\bZ$.
One manifestation of that principle is proved below: there is a Poincar\'e lemma for infinitesimal cohomology
which asserts that the natural map $k\rightarrow\Inf_{R/k}$ is an equivalence for any $R$.

The names infinitesimal and crystalline are justified. For a smooth $\bF_p$-algebra $R$,
it is known that $\dR_{R/\bF_p}\we\dRhat_{R/\bF_p}$ and these compute the cohomology of the structure sheaf of the crystalline site $(R/\bF_p)_\crys$.
Crystalline filtered derived commutative rings are also related to animated pd-pairs, which can be
used to define derived crystalline cohomology~\cite{magidson-divided,mao-crystalline}.
Similarly, for Grothendieck's infinitesimal site $(R/\bF_p)_\inf$, one has
$\R\Gamma((R/\bF_p)_\inf,\Oscr)\we R^\flat$, the inverse limit perfection of $R$. We prove in
Section~\ref{sec:charp} that
$\Infhat_{R/\bF_p}\we R^\flat$ for any quasisyntomic $\bF_p$-algebra $R$. This fact has also been observed
by Jiaqi Fu and Akhil Mathew in private communications.

The filtrations which arise in prismatic cohomology are typically infinitesimal in nature. This
fact has been used by Holeman in~\cite{holeman-derived} to state a universal property for derived prismatic
cohomology relative to a prism.

Our main theorem in this paper asserts that there is a
natural map $$\F^\star_\H\Inf_{R/k}\rightarrow\F^\star_\H\dR_{R/k}.$$ In fact, it is the unit map
of the crystallization adjunction. To state it, we must restrict to non-negatively filtered objects
$\F^+\Mod_k$.

\begin{theorem}\label{thm:intromain}
    Let $k$ be a derived commutative ring.
    \begin{enumerate}
        \item[{\em (1)}]
            There is a map of monads $\LSym^\inf\rightarrow\LSym^\crys$ on $\F^+\Mod_k$ inducing an adjunction
            $$\divideontimes\colon\F^+\DAlg_k^\inf\rightleftarrows\F^+\DAlg_k^\crys\colon\text{{\em forget}}$$
            for any derived commutative ring $k$ and there are canonical
            identifications $\divideontimes(\F^\star_\H\Inf_{R/k})\we\F^\star_\H\dR_{R/k}$ for all derived
            commutative rings $k$ and $R$.
        \item[{\em (2)}]
            There is a map of monads $\LSym^\inf\rightarrow\LSym^\crys$ on $\widehat{\F^+\Mod}_k$ inducing an adjunction
            $$\divideontimes\colon\widehat{\F^+\DAlg}_k^\inf\rightleftarrows\widehat{\F^+\DAlg}_k^\crys\colon\text{{\em forget}}$$
            for any derived commutative ring $k$ and there are canonical
            identifications $\divideontimes(\F^\star_\H\Infhat_{R/k})\we\F^\star_\H\dRhat_{R/k}$ for all derived
            commutative rings $k$ and $R$.
    \end{enumerate}
\end{theorem}

Part (1) of the theorem will follow from part (2). To prove part (2) we
study the process of deriving monads and how this process interacts with various $t$-structures.
This suffices for the graded analogue of Theorem~\ref{thm:intromain}. More subtle arguments are
needed to upgrade from the graded to the filtered setting.

This paper arose from our study of the work~\cite{raksit} of Raksit who defined the filtered
circle $\bT_\fil$ and established the universal properties of $\HH_\fil(R/k)$ and the complete
version of de Rham cohomology, $\F^\star_\H\dRhat_{R/k}$. Raksit also explains conceptually the
close relationship between Hochschild homology and de Rham cohomology which goes back to work of
many authors~\cite{antieau-derham,bzn,bms2,feigin-tsygan-additive,loday,toen-vezzosi-simpliciales}.
Raksit's algebraic approach complements the stacky approach of Moulinos--Robalo--To\"en~\cite{mrt}.

We wanted to know if $\dR_{R/k}$ and $\dRhat_{R/k}$ are derived commutative rings and if
$\F^\star_\H\dR_{R/k}$ and $\F^\star_\H\dRhat_{R/k}$
are infinitesimal filtered derived commutative rings. The present paper answers these and related
questions in the affirmative via the crystallization adjunction of Theorem~\ref{thm:intromain}.

Raksit pointed out to us another way to see that $\widehat{\dR}_{R/k}$ admits the
structure of a derived commutative $k$-algebra.
Upon taking filtered $S^1$-fixed points on $\HH_\fil(R/k)$, one obtains an
infinitesimal filtered derived commutative $k$-algebra structure $\F^\star_\mot\HC^-(R/k)$ on
$\HC^-(R/k)$ such that
$\gr^0\HC^-(R/k)$ is equivalent to $\dRhat_{R/k}$, so that the latter is a derived commutative ring.
This argument can be derived from~\cite{mrt,raksit}.

The existence of the filtered circle has consequences for the computation of differentials in the HKR spectral
sequence as discovered by Mundinger~\cite{mundinger}.
The filtered circle has been further explored, in synthetic settings, by
Hedenlund--Moulinos~\cite{hedenlund-moulinos} and in our joint work with
Riggenbach~\cite{antieau-riggenbach}.

The forgetful functor from derived commutative rings to $\bE_\infty$-rings preserves limits, from
which it follows that if $(\Xscr,\Oscr)$ is a topos with a sheaf of commutative rings $\Oscr$, then $\R\Gamma(\Xscr,\Oscr)$
is a derived commutative ring. This can be used in some cases to give another proof that
$\dR_{R/k}$ or $\dRhat_{R/k}$ is a derived commutative ring.

\paragraph{Related work.}
As mentioned above, Adam Holeman used our definition of infinitesimal cohomology
in his work on prismatic cohomology~\cite{holeman-derived}.
Jiaqi Fu has used our notion of infinitesimal cohomology in work on infinitesimal foliations~\cite{fu-inf}. To\"en and
Vezzosi discuss an infinitesimal version of derived foliation in~\cite{tv-inf} and in their
book~\cite{tv-book} from a stacky
perspective.
Theorem~\ref{thm:intromain} answers~\cite[Q.~5.4]{tv-inf} in the affirmative and provides an
alternative approach to the de Rhamification functor of~\cite{tv-book} which is restricted to the
case of foliations.
These derived foliations are related in the work of Brantner--Magidson--Nuiten on deformation theory,
which connects infinitesimal filtered derived commutative rings to the partition Lie algebras
of Brantner--Mathew~\cite{brantner-mathew}.

\paragraph{Outline.}
Sections~\ref{sec:filtered_monad},~\ref{sec:polynomial_monad},~\ref{sec:derived_algebras}, and~\ref{sec:dalg}
develop the algebraic theory necessary to make our definitions and construct various maps of
monads. Sections~\ref{sec:grdalg} and~\ref{sec:fdalg} introduce various flavors of graded and
filtered derived commutative rings. In Section~\ref{sec:raksit}, we give a summary, with proofs, of
the main points of Raksit's work~\cite{raksit}. We introduce infinitesimal cohomology in
Section~\ref{sec:inf} and study it for quasisyntomic $\bF_p$-algebras in Section~\ref{sec:charp}.
Finally, the crystallization adjunction is established in Section~\ref{sec:crystallization}.

\paragraph{Future work.} Later titles in this series are
\begin{enumerate}
    \item[II:] {\em the Gauss--Manin connection};
    \item[III:] {\em cohomology theories for $\bE_\infty$-rings};
    \item[IV:] {\em modules with integrable connection}.
\end{enumerate}
We will refer to them in the text as F\&C.II, etc.
In F\&C.II, we establish a fully algebraic and derived version of the Gauss--Manin connection and use it
to establish completeness conditions for de Rham and infinitesimal cohomology, allowing for example
for simple proofs of the main results of~\cite{bhatt-completions}. In F\&C.III, we give a survey of
$\bE_\infty$-infinitesimal cohomology and its properties. No proofs will be required because they are
identical to those in the derived infinitesimal case. Part IV develops the theory of modules with
integrable connection over general filtered derived commutative rings and relates these to
classical notions of $D$-module.

In more distant future work, we hope to return the project begun by Holeman
in~\cite{holeman-derived} to understand the nature of the filtrations arising in prismatic
cohomology. This is related to forthcoming work of Kirill Magidson on universal properties for the de Rham--Witt
complex.

\paragraph{Notation.}
We work with $\infty$-categories throughout this paper. However, `category' means a
$1$-category and commutative ring means an ordinary (static) commutative ring. Unless otherwise
specified, all algebraic operations, such as $\otimes$ or $\Lambda^i$, are derived. If $k$ is an
$\bE_\infty$-ring, we let $\Mod_k=\Mod_k(\Sp)$ denote the $\infty$-category of $k$-modules.
There is a natural $t$-structure on $\Mod_k$ whose connective objects $\Mod_k^\cn=(\Mod_k)_{\geq
0}$ is the full subcategory of $\Mod_k$ generated under colimits by $k$ itself and closed under
extensions. The coconnective objects $(\Mod_k)_{\leq 0}$ are those $M$ where $\pi_i M=0$ for $i>0$.
When $k$ is connective, $\Mod_k^\heart=\Mod_{\pi_0k}(\Ascr\mathrm{b})$. An object
$M\in\Mod_k^\heart$ is called a static $k$-module. Objects in the $\infty$-category $\Sscr$ of Kan complexes
up to weak equivalence are called anima.

\paragraph{Acknowledgments.} This project was born out of a conversation
with Achim Krause in March 2020 and heavily influenced by discussions with
Bhargav Bhatt, Lukas Brantner, Haoyang Guo, Akhil Mathew, Thomas Nikolaus, Joost Nuiten, Nick Rozenblyum,
and Bertrand To\"en.
I want to especially thank Carlos Cortez, Sanath Devalapurkar, Adam Holeman, Kirill Magidson, Deven Manam,
Tasos Moulinos, and Arpon Raksit for many discussions when the ideas here were just forming and Elden Elmanto for introducing me to
Sanath at a critical point. Arpon in particular patiently answered too many questions about his work.
Finally, several people have provided exceptionally helpful feedback on early drafts, including Jiaqi Fu, Adam
Holeman, Kirill Magidson, Tasos Moulinos, and Arpon Raksit.

Some of the results here were explained at Oberwolfach in 2022 at the workshop {\em Non-commutative
geometry and cyclic homology} and reported on in~\cite{antieau_whatis}. We apologize for the long
gestational period for this paper.

This project was supported by NSF grants DMS-2120005, DMS-2102010, and
DMS-2152235, Simons Fellowships 666565 and 00005925, and the Simons Collaboration on Perfection.

\section{Filtered monads}\label{sec:filtered_monad}

A \defidx{monad} on an $\infty$-category $\Cscr$ is an algebra
object $T$ of the monoidal $\infty$-category
$\End(\Cscr)=\Fun(\Cscr,\Cscr)$, where the monoidal structure arises from composition
of functors. If $\Escr\subseteq\End(\Cscr)$ is a full monoidal subcategory,
then an algebra object of $\Escr$ is called an
\longdefidx{$\Escr$-monad}{monad!$\Escr$-monad}. The action of endomorphisms makes $\Cscr$ into a left
$\End(\Cscr)$-module in the $\infty$-category $\Cat_\infty$ of $\infty$-categories (specifically, $\Cscr$ is left-tensored over $\End(\Cscr)$ in the
sense of~\cite[Def.~4.2.1.19]{ha}), so for any algebra object (i.e.,
monad) $T$ in
$\End(\Cscr)$ there is an $\infty$-category $\LMod_T(\Cscr)$ of left
$T$-modules in $\Cscr$ generalizing the classical notation of a module over a monad in a
$1$-category; see~\cite[Def.~4.2.1.13]{ha}.

There is a basic categorical lemma on the interplay of
endomorphisms and adjunctions which will be used many times below.

\begin{lemma}\label{lem:oplax}
    Let $F\colon\Cscr\rightleftarrows\Dscr\colon G$ be an adjunction between
    $\infty$-categories. Then, the natural intertwining functor
    $\End(\Cscr)\rightarrow\End(\Dscr)$ given by $T\mapsto F\circ
    T\circ G$ is naturally oplax monoidal (using the unit of the adjunction) with right adjoint the lax monoidal
    functor given by $T'\mapsto G\circ
    T'\circ F$ (using the counit of the adjunction).
\end{lemma}

\begin{proof}
    To prove the adjointness of these functors, we follow the first paragraph of the proof
    of~\cite[Lem.~3.10]{brantner-campos-nuiten}. The adjunction $\Cscr\rightleftarrows\Dscr$ is
    classified by a functor $\Mscr\rightarrow\Delta^1$ which is both a Cartesian and coCartesian
    fibration and where the fiber over $0$ is identified with $\Cscr$ and the fiber over $1$ is
    identified with $\Dscr$. Consider the $\infty$-category $\End_{/\Delta^1}(\Mscr)$ of
    endomorphisms of $\Mscr\rightarrow\Delta^1$. There are then restriction functors
    $$\End(\Cscr)\leftarrow\End_{\Delta^1}(\Mscr)\rightarrow\End(\Dscr).$$
    The left arrow admits a left adjoint and the right arrow a right adjoint, by relative left and
    right Kan extensions over $\Delta^1$. The compositions are identified
    in~\cite{brantner-campos-nuiten} as the given functors $T\mapsto F\circ T\circ G$ and $T'\mapsto
    G\circ T'\circ F$ which are thus adjoint as claimed.

    The lax monoidality of the right adjoint is
    in~\cite[Lem.~3.10]{brantner-campos-nuiten}. In general, the left adjoint
    of a lax monoidal functor is oplax monoidal.
\end{proof}

\begin{remark}
    As a consequence of Lemma~\ref{lem:oplax}, there is a natural functor
    $\Alg(\End(\Cscr))\leftarrow\Alg(\End(\Dscr))$, which takes monads on $\Dscr$ to monads on
    $\Cscr$. Given a monad $T$ on $\Dscr$, the underlying functor of the resulting monad on $\Cscr$
    is $UTF$. The right adjoint $U$ induces a forgetful functor
    $\LMod_{UTF}(\Cscr)\leftarrow\LMod_T(\Dscr)$.
\end{remark}

The following definition is due to Raksit.

\begin{definition}[Filtered monad]\label{def:filteredmonad}
    A \longdefidx{filtered monad}{monad!filtered} on an $\infty$-category $\Cscr$ is a lax monoidal
    functor $\F_\star T\colon\bZ_{\geq 0}^\times\rightarrow\End(\Cscr)$, where $\bZ_{\geq
    0}^\times$ is the set of nonnegative integers equipped with the
    multiplicative monoidal structure. If $\Escr\subseteq\End(\Cscr)$ is
    a full monoidal subcategory, then a \longdefidx{filtered
    $\Escr$-monad}{monad!filtered $\Escr$-monad} is a
    lax monoidal functor $\F_\star T\colon\bZ_{\geq 0}^\times\rightarrow\Escr$.
\end{definition}

Under mild conditions, the colimit $T$ of a filtered $\Escr$-monad $\F_\star T$ exists
and is itself an algebra object in $\Escr$ and hence a monad on $\Cscr$. For
example, Raksit proves the following lemma; see~\cite[Prop.~4.1.4]{raksit}.

\begin{lemma}\label{lem:colimit}
    Suppose that $\Cscr$ is an $\infty$-category admitting colimits, that
    $\Escr\subseteq\End(\Cscr)$ is a monoidal subcategory admitting sequential
    colimits, and that each $F\in\Escr$ preserves sequential colimits as a
    functor $\Cscr\rightarrow\Cscr$. If $\F_\star T$ is a filtered
    $\Escr$-monad, then $\colim\F_\star T$ is naturally an $\Escr$-monad.
\end{lemma}

\begin{definition}[Filtered $\Sigma$-monads]
    Suppose that $\Cscr$ admits sifted colimits. In this case, the full
    subcategory $\End_\Sigma(\Cscr)\subseteq\End(\Cscr)$ of endomorphisms which preserve sifted
    colimits is a monoidal
    subcategory. A
    \longdefidx{$\Sigma$-monad}{monad!$\Sigma$-monad} is an
    $\End_\Sigma(\Cscr)$-monad. A \longdefidx{filtered $\Sigma$-monad}{monad!filtered
    $\Sigma$-monad} is a filtered
    $\End_\Sigma(\Cscr)$-monad. If $\Cscr$ has all small colimits, then the
    colimit lemma applies to filtered $\Sigma$-monads.
\end{definition}

The point of introducing filtered monads is that the techniques (explained below) used to extend each
$\LSym^d$ from $\Mod_\bZ^\cn$ to $\Mod_\bZ$ do not immediately apply to
$\LSym\we\oplus_{d\geq 0}\LSym^d$ because the infinite coproduct is not
polynomial. The use of filtered monads gives a way of making the appropriate
extensions for monads which are exhaustively filtered by polynomial functors.
Another, closely related approach, is taken in forthcoming work with H\"ubner, Kubrak,
and Nuiten.

\section{Polynomial monads}\label{sec:polynomial_monad}

We discuss various approximations to monads on a stable $\infty$-category $\Cscr$, most with respect to a
suitable $t$-structure $(\Cscr_{\geq 0},\Cscr_{\leq 0})$. Some of these games require the monad to be built out of certain polynomial
functors $F$. Polynomiality ensures that the values of the functor $F$ on general objects of $\Cscr$
can be deduced from its values on $\Cscr_{\geq 0}$.

\begin{definition}[Polynomial functors]\label{def:polynomial}
    Let $F\colon\Cscr\rightarrow\Dscr$ be a functor of additive
    $\infty$-categories where $\Dscr$ is idempotent complete.
    \begin{enumerate}
        \item[(a)] Say that $F$ is additively polynomial of degree
            $\leq 0$ if it is constant. Inductively, for $d\geq 1$, $F$ is
            additively polynomial of degree $\leq d$ if, for every $X$ in
            $\Cscr$, the derivative functor
            $$D_XF\we\mathrm{fib}(F(X\oplus(-))\rightarrow F(-))$$ is additively
            polynomial of degree $\leq d-1$. Let
            $\Fun^{\ap}(\Cscr,\Dscr)\subseteq\Fun(\Cscr,\Dscr)$ be the full
            subcategory of additively polynomial
            functors, i.e., those which
            are additively polynomial of degree $\leq d$ for some $d\geq 0$.
        \item[(b)] If $\Cscr$ and $\Dscr$ are prestable (in particular, if they
            are stable), $F$ is
            excisively polynomial of degree $\leq d$ if it is additively polynomial
            of degree $\leq d$ and if it preserves finite geometric
            realizations. Let
            $\Fun^{\ep}(\Cscr,\Dscr)\subseteq\Fun(\Cscr,\Dscr)$ denote the
            full subcategory of excisively polynomial
            functors, i.e., those
            which are excisively polynomial of degree $d$ for some $d\geq 0$.
    \end{enumerate}
\end{definition}

\begin{notation}
    If $\Cscr$ and $\Dscr$ are additive $\infty$-categories which admit sifted colimits, let
    $\Fun^{\ap}_\Sigma(\Cscr,\Dscr)\subseteq\Fun^\ap(\Cscr,\Dscr)$ be the full subcategory of
    additively polynomial functors which preserve sifted colimits. If $\Cscr$ and $\Dscr$ are
    additionally prestable, let $\Fun^\ep_\Sigma(\Cscr,\Dscr)\subseteq\Fun^\ep(\Cscr,\Dscr)$ be the full
    subcategory of excisively polynomial functors which preserve sifted colimits.
\end{notation}

\begin{example}
    If $F$ is additively (resp. excisively) polynomial of degree $\leq d$ for some integer $d\geq 0$,
    then it is additively (resp. excisively) polynomial of degree $\leq d'$ for all integers $d'\geq d$.
\end{example}

\begin{example}\label{ex:composition}
    Let $F\colon\Cscr\rightarrow\Dscr$ and $G\colon\Dscr\rightarrow\Escr$ be additively polynomial
    functors of degrees $\leq d$ and $\leq e$ between additive $\infty$-categories for some
    $d,e\geq 0$.
    Then, $G\circ F$ is additively polynomial of degree $\leq de$. If $\Cscr$, $\Dscr$, and $\Escr$
    are prestable and $F$ and $G$ are additionally excisively polynomial, then $G\circ F$ is
    excisively polynomial of degree $\leq de$.
\end{example}

The proof that when $\Dscr$ is stable the definition of excisive polynomial agrees with the original
definition of Goodwillie~\cite{goodwillie-calc2} can be found
in~\cite[Prop.~2.15]{bgmn}, but see also~\cite{johnson-mccarthy}. The definition of
additively polynomial functors of additive (1-)categories goes back to Eilenberg and Mac Lane
in~\cite{eml2}. Moreover, a functor $F\colon\Cscr\rightarrow\Dscr$ of additive $\infty$-categories
is additively polynomial of degree $\leq d$ if and only if the corresponding functor
$\Ho(F)\colon\Ho(\Cscr)\rightarrow\Ho(\Dscr)$ on (additive) homotopy categories is additively polynomial of
degree $\leq d$.

\begin{warning}\label{warning:prestable}
    When $\Dscr$ is prestable but not stable, the definition given here of excisive polynomial
    functor is not in general the same as the classical definition in terms of
    coCartesian cubes. However, if $\Dscr$ is prestable, then
    $\Dscr$ admits a fully faithful embedding
    $\Dscr\hookrightarrow\mathrm{SW}(\Dscr)$ into its Spanier--Whitehead, a
    stable $\infty$-category, and $\Dscr$ is closed under finite colimits in
    $\mathrm{SW}(\Dscr)$; see~\cite[Prop.~C.1.2.2]{sag}. Thus, a functor
    $F\colon\Cscr\rightarrow\Dscr$ is excisively polynomial in the sense of
    Definition~\ref{def:polynomial} if and only the composition
    $\Cscr\xrightarrow{F}\Dscr\hookrightarrow\mathrm{SW}(\Dscr)$ is excisively
    polynomial, which as remarked above agrees with Goodwillie's definition.
\end{warning}

The following lemma will be used throughout the following discussion.

\begin{lemma}\label{lem:easy}
    Let $\Cscr$ and $\Dscr$ be prestable $\infty$-categories admitting small colimits.
    An additively polynomial functor $F\colon\Cscr\rightarrow\Dscr$
    preserves sifted colimits if and only if it is excisively polynomial and
    preserves filtered colimits.
\end{lemma}

\begin{proof}
    If $F$ preserves sifted colimits, then it preserves finite geometric realizations and filtered
    colimits; if it is additionally additively polynomial it follows it is excisively polynomial.
    Conversely, if $F$ preserves filtered colimits and finite geometric realizations, then it
    preserves all geometric realizations as these can be written as filtered colimits of geometric
    realizations (using the skeletal filtration). Functors which preserve filtered colimits and
    geometric realizations preserve all sifted colimits by~\cite[Cor.~5.5.8.17]{htt}.
\end{proof}

\begin{definition}
    Fix a $t$-structure $(\Cscr_{\geq 0},\Cscr_{\leq 0})$ on a stable
    presentable $\infty$-category $\Cscr$. We assume that the $t$-structure is accessible, which
    means that $\Cscr_{\geq 0}$ is presentable too~\cite[Def.~1.4.4.12]{ha}.
    \begin{enumerate}
        \item[(i)]
            Let $\End_\Sigma^\ep(\Cscr)\subseteq\End(\Cscr)$ be the full
            subcategory of sifted colimit preserving excisively polynomial
            endofunctors. Similarly, let $\End^\ep(\Cscr)$ be the full subcategory
            of excisively polynomial endofunctors and let $\End_\Sigma(\Cscr)$
            be the full subcategory of sifted colimit preserving endofunctors.
            All three subcategories are closed under composition and hence are
            monoidal by Example~\ref{ex:composition}.
        \item[(ii)] Let $\End^\ep_\Sigma(\Cscr_{\geq 0})$ be the full monoidal subcategory of sifted
            colimit preserving excisively polynomial endofunctors.
        \item[(iii)] Fix an integer $n\geq 1$ and let $\Cscr_{[0,n-1]}$ be the
            Grothendieck abelian $n$-category associated to $\Cscr_{\geq 0}$ in the sense
            of~\cite[App.~C]{sag}. Let $\End^\ap_\Sigma(\Cscr_{[0,n-1]})$ be
            the full monoidal subcategory of sifted colimit preserving
            additively polynomial endofunctors.
    \end{enumerate}
    A \longdefidx{filtered excisively polynomial $\Sigma$-monad}{monad!filtered
    excisively polynomial $\Sigma$-monad} on $\Cscr$ or $\Cscr_{\geq 0}$ is a filtered
    $\End_\Sigma^\ep(\Cscr)$-monad or $\End_\Sigma^\ep(\Cscr_{\geq 0})$-monad in the sense of
    Definition~\ref{def:filteredmonad}. An \longdefidx{additively polynomial
    $\Sigma$-monad}{monad!additively polynomial $\Sigma$-monad} on
    $\Cscr_{[0,n-1]}$ is an $\End^\ap_\Sigma(\Cscr_{[0,n-1]})$-monad.
\end{definition}

\begin{remark}
    If $\Cscr$ is a cocomplete stable $\infty$-category,
    then Lemma~\ref{lem:colimit} applies to
    $\End_\Sigma(\Cscr)$, but not to $\End_\Sigma^\ep(\Cscr)$ or
    $\End^\ep(\Cscr)$. Indeed, sequential colimits of polynomial functors are
    not typically polynomial. This could possibly be fixed by formally adding sequential colimits
    of polynomial functors
    to obtain some kind of notion of analytic monad (in a different sense from
    how analytic is sometimes used in monad theory, for example
    in~\cite{gepner-haugseng-kock}), but it will not be necessary to make this
    precise thanks to the use of filtered monads. See our forthcoming work with H\"ubner, Kubrak,
    and Nuiten for an alternative approach,
    which uses the notion of right extendability from~\cite{brantner-mathew}.
\end{remark}

The following observation is the key to
deriving monads: it gives a functorial extension of certain functors on
$\Cscr_{\geq 0}$ to $\Cscr$. Recall that a $t$-structure $(\Cscr_{\geq 0},\Cscr_{\leq 0})$ on a
stable presentable $\infty$-category $\Cscr$ is compatible with filtered colimits if the inclusion
$\Cscr_{\leq 0}\subseteq\Cscr$ preserves filtered colimits. It is called right complete if the
natural map $\Cscr\rightarrow\lim_{m\mapsto-\infty}\Cscr_{\geq m}$ is an
equivalence~\cite[Rem.~1.2.1.20]{ha}.

\begin{proposition}\label{prop:calc}
    If $(\Cscr_{\geq 0},\Cscr_{\leq 0})$ is an accessible right complete
    $t$-structure which is compatible with filtered colimits on a
    stable presentable $\infty$-category $\Cscr$ and if $\Dscr$ is a stable
    $\infty$-category with colimits, then the restriction functor
    $$\Fun^\ep_\Sigma(\Cscr,\Dscr)\rightarrow\Fun^\ep_\Sigma(\Cscr_{\geq
    0},\Dscr)$$ is an equivalence.
\end{proposition}

\begin{proof}
    See the proofs of~\cite[Prop.~4.2.15]{raksit} and~\cite[Thm.~3.36]{brantner-mathew}.
    If $F\colon\Cscr_{\geq 0}\rightarrow\Dscr$ is an excisively polynomial functor
    of degree $n$ which preserves sifted colimits, consider
    the $d$th Goodwillie approximation $\P_d(F\circ\tau_{\geq 0})$ of the
    composition $F\circ\tau_{\geq 0}$. This functor is excisively polynomial of degree $d$, by
    construction, and preserves sifted colimits
    by~\cite[Prop.~3.37]{brantner-mathew} since $\tau_{\geq 0}$ preserves filtered
    colimits by the assumption that the $t$-structure is compatible with filtered
    colimits. In~\cite{brantner-mathew, raksit}, this construction is shown to be an inverse to the
    restriction functor.
\end{proof}

\begin{proposition}\label{prop:siftedapproximation}
    If $\Cscr$ is a compactly generated stable presentable
    $\infty$-category,
    then the inclusion of the full monoidal subcategory
    $\End_\Sigma^\ep(\Cscr)\subseteq\End^\ep(\Cscr)$ admits a lax
    monoidal right adjoint.
\end{proposition}

\begin{proof}
    Let $\Cscr^\omega\subseteq\Cscr$ be the full subcategory of compact
    objects. Given an excisively polynomial functor $F\colon\Cscr\rightarrow\Cscr$, 
    restriction along $\Cscr^\omega\rightarrow\Cscr$ results in a functor
    $F'\colon\Cscr^\omega\rightarrow\Cscr$ which is excisively polynomial.
    The left Kan extension $\L F'\colon\Cscr\rightarrow\Cscr$ is
    excisively polynomial by~\cite[Prop.~6.1.5.4]{ha} and now preserves
    filtered colimits and hence all geometric realizations
    by~\cite[Prop.~3.37]{brantner-mathew}.
    One checks that the construction $F\mapsto\L F'$ defines a colocalization
    functor using (3) implies (1) of~\cite[Prop.~5.2.7.4]{htt}. Thus, $F\mapsto\L F'$ is a right
    adjoint to the inclusion; this right adjoint canonically admits the structure of a lax monoidal
    functor since the inclusion is monoidal.
\end{proof}

\begin{definition}[$\Sigma$-approximation]
    If $T$ is an excisively
    polynomial monad, then the counit map $T_\Sigma\rightarrow T$ from the adjunction of
    Proposition~\ref{prop:siftedapproximation} is called the
    \longdefidx{$\Sigma$-approximation}{monad!$\Sigma$-approximation} of $T$;
    it is a map of $\End^\ep(\Cscr)$-monads on $\Cscr$.
    Given a filtered excisively polynomial monad $\F_\star T$ the counit map
    $(\F_\star T)_\Sigma\rightarrow \F_\star T$ is called the
    \longdefidx{$\Sigma$-approximation}{monad!$\Sigma$-approximation} of $\F_\star T$. 
    It is a map of filtered $\End^\ep(\Cscr)$-monads on $\Cscr$.
\end{definition}

\begin{proposition}\label{prop:connective}
    Let $\Cscr_{\geq 0}$ be a Grothendieck prestable $\infty$-category and let $\Dscr$ be a
    stable presentable $\infty$-category with an accessible $t$-structure
    $(\Dscr_{\geq 0},\Dscr_{\leq 0})$.
    If $\Cscr_{\geq 0}$ is compact
    projectively generated, then composition with $\Dscr_{\geq 0}\rightarrow\Dscr$ induces a fully
    faithful left adjoint functor
    $$ \Fun_\Sigma^\ep(\Cscr_{\geq 0},\Dscr_{\geq 0})\rightarrow\Fun_\Sigma^\ep(\Cscr_{\geq 0},\Dscr)$$
    whose essential image consists of those functors which take values in $\Dscr_{\geq 0}$.
\end{proposition}

\begin{proof}
    Let $F\colon\Cscr_{\geq 0}\rightarrow\Dscr$ be an excisively polynomial
    functor which preserves sifted colimits. Then, the composition $\tau_{\geq 0}F=F\circ\tau_{\geq
    0}\colon\Cscr_{\geq 0}\rightarrow\Dscr_{\geq 0}$ is additively polynomial.
    Let $\Cscr^0\subseteq\Cscr_{\geq 0}$ be the full subcategory of compact
    projective objects. The restriction of $\tau_{\geq 0}F$ to $\Cscr^0$ is
    additively polynomial and there is a natural map of functors $\L\tau_{\geq
    0}F\rightarrow\tau_{\geq 0}F\colon\Cscr_{\geq 0}\rightarrow\Dscr_{\geq 0}$ obtained by left Kan extension. By
    construction, $\L\tau_{\geq 0}F$ preserves sifted colimits and is
    additively polynomial (see~\cite[Prop.~5.10]{johnson-mccarthy}
    or~\cite[Prop.~3.35]{brantner-mathew}); hence it is excisively polynomial by
    Lemma~\ref{lem:easy}. The assignment $F\mapsto\L\tau_{\geq 0}F$ is a
    colocalization as one checks again using~\cite[Prop.~5.2.7.4]{htt}.
\end{proof}

\begin{definition}
    In the setting of Proposition~\ref{prop:connective}, we let $F\mapsto\tau_{\geq 0}F$ denote 
    the right adjoint of $\Fun_\Sigma^\ep(\Cscr_{\geq 0},\Dscr_{\geq
    0})\rightarrow\Fun_\Sigma^\ep(\Cscr_{\geq 0},\Dscr)$.
\end{definition}

\begin{warning}
    Note that despite the name,
    $\tau_{\geq 0}F$ does not take connective values in general, although it will take
    connective values on connective objects.
\end{warning}

Suppose that $\Cscr$ is a stable $\infty$-category with an accessible
$t$-structure $(\Cscr_{\geq 0},\Cscr_{\leq 0})$ which is right complete and compatible with
filtered colimits. The equivalence of Proposition~\ref{prop:calc} makes
$\Fun^\ep_\Sigma(\Cscr_{\geq 0},\Cscr)\we\Fun^\ep_\Sigma(\Cscr,\Cscr)$ into
a monoidal $\infty$-category where the monoidal structure is induced from
the composition monoidal structure on the right-hand side. The unit object,
corresponding to the identity functor on the right-hand side, is the inclusion
$\Cscr_{\geq 0}\rightarrow\Cscr$ on the left-hand side.

\begin{corollary}\label{cor:connective}
    Let $\Cscr$ be a stable $\infty$-category with an accessible $t$-structure
    $(\Cscr_{\geq 0},\Cscr_{\leq 0})$ which is right complete and compatible
    with filtered colimits. If $\Cscr_{\geq 0}$ is compact projectively generated, then
    the fully faithful inclusion
    $$\Fun_\Sigma^\ep(\Cscr_{\geq 0},\Cscr_{\geq
    0})\rightarrow\Fun_\Sigma^\ep(\Cscr_{\geq
    0},\Cscr)\we\Fun_\Sigma^\ep(\Cscr,\Cscr)$$ is monoidal and admits a lax
    monoidal right adjoint $F\mapsto \tau_{\geq 0}F$.
\end{corollary}

\begin{proof}
    The existence of the right adjoint is the content of
    Proposition~\ref{prop:connective}. Its lax monoidality will follow from the monoidality of the
    the inclusion.
    If $F,G\colon\Cscr_{\geq 0}\rightarrow\Cscr$
    are excisively polynomial functors of degrees $\leq d$ and $\leq e$ which preserve
    sifted colimits, then the recipe for computing $F\circ G$ is to take
    $(\P_d(F\circ\tau_{\geq 0}))\circ(\P_e(G\circ\tau_{\geq 0}))$ and to restrict
    it to $\Cscr_{\geq 0}$. Note that if $\iota\colon\Cscr_{\geq 0}\hookrightarrow\Cscr$ denotes the
    canonical inclusion, then $\P_e(G\circ\tau_{\geq 0})\circ\iota\we\P_e(G\circ\tau_{\geq
    0}\circ\iota)\we\P_e(G)\we G$ by~\cite[Rem.~6.1.1.30]{ha}.
    Thus, the restriction of $(\P_d(F\circ\tau_{\geq 0}))\circ(\P_e(G\circ\tau_{\geq 0}))$ to
    $\Cscr_{\geq 0}$ is naturally equivalent to
    $\P_d(F\circ\tau_{\geq 0})\circ G$. If $G$ takes connective values, then
    this is naturally equivalent to $F\circ G$ for the same reason. By the remark before
    the proposition, the left adjoint functor preserves the units of the
    monoidal structures. It follows that the full subcategory $\Escr$ of $\Fun_\Sigma^\ep(\Cscr_{\geq
    0},\Cscr)\we\Fun_\Sigma^\ep(\Cscr,\Cscr)$
    consisting of functors which take connective values on connective objects is in fact closed
    under the monoidal structure on $\Fun_\Sigma^\ep(\Cscr,\Cscr)$. Let $\Escr^\otimes$ denote the
    induced $\infty$-operad structure on $\Escr$. By Proposition~\ref{prop:connective}, the
    underlying $\infty$-category of this $\infty$-operad is $\Fun_\Sigma^\ep(\Cscr_{\geq
    0},\Cscr_{\geq 0})$. It suffices to see that the $\infty$-operad structure $\Escr^\otimes$
    agrees with the composition monoidal structure. For this, note that the action of $\Escr$ on
    $\Cscr$ restricts to an action on $\Cscr_{\geq 0}$. This makes $\Cscr_{\geq 0}$ left-tensored
    over $\Ascr$. To see this rigorously, if $\Fscr\rightarrow\mathbf{LMod}^\otimes$ is the
    fibration of $\infty$-operads classifying the left-tensoring of $\Cscr$ over
    $\End_\Sigma^\ep(\Cscr)$, then one can look at the sub-$\infty$-operad of $\Fscr$ generated by objects in
    $\Cscr_{\geq 0}$ or $\Ascr$ in the sense of~\cite[Sec.~2.2.1]{ha}. The resulting object is
    still $\mathbf{LMod}$-monoidal by~\cite[Prop.~2.2.1.1]{ha} and it classifies an action of $\Ascr$ on $\Cscr_{\geq 0}$.
    By universality of the endomorphism $\infty$-category construction
    in~\cite[Sec.~4.7.1]{ha}, applied to $\Cat_\infty$ itself, it follows that there is a functor of
    monoidal $\infty$-categories $\Escr^\otimes\rightarrow\Fun_\Sigma^\ep(\Cscr_{\geq
    0},\Cscr_{\geq 0})$. At the level of underlying $\infty$-categories, this functor is an inverse
    to the equivalence arising from the identification of the essential image in
    Proposition~\ref{prop:connective}. However, a monoidal functor which is an equivalence on underlying categories
    is an equivalence.
\end{proof}

\begin{definition}[Connective cover]
    If $T\in\End_\Sigma^\ep(\Cscr)$ is
    a monad, then the associated map $\tau_{\geq 0}T\rightarrow T$ of $\End_\Sigma^{\ep}(\Cscr)$-monads is called
    the \longdefidx{connective cover}{monad!connective cover} of $T$ and
    similarly for filtered monads taking values in $\End_\Sigma^\ep(\Cscr)$. 
    A (filtered) monad $T\in\End_\Sigma^\ep(\Cscr)$ is connective if the counit $\tau_{\geq
    0}T\rightarrow T$ is an equivalence. Equivalently, $T$ is connective if it preserves connective
    objects.
\end{definition}

Given a stable $\infty$-category $\Cscr$ with a $t$-structure $(\Cscr_{\geq 0},\Cscr_{\leq 0})$ and
$1\leq n\leq\infty$, we
let $\Cscr_{[0,n-1]}=\Cscr_{\geq 0}\cap\Cscr_{\leq n-1}$. This is an additive $\infty$-category which is
typically neither stable nor prestable. We have $\Cscr_{[0,0]}\we\Cscr^\heart$.

\begin{proposition}\label{prop:derivedapproximation}
    Fix a stable presentable $\infty$-category $\Cscr$ with an accessible
    $t$-structure $(\Cscr_{\geq 0},\Cscr_{\leq 0})$ which is compatible with
    filtered colimits and such that $\Cscr_{\geq 0}$ is compact projectively generated.
    Fix an integer $n\geq 1$.
    Let $\iota\colon\Cscr_{[0,n-1]}\rightarrow\Cscr_{\geq 0}$ be the right adjoint
    to $\tau_{\leq n-1}\colon\Cscr_{\geq 0}\rightarrow\Cscr_{[0,n-1]}$.
    The oplax monoidal functor $\Fun(\Cscr_{\geq 0},\Cscr_{\geq
    0})\rightarrow\Fun(\Cscr_{[0,n-1]},\Cscr_{[0,n-1]})$ taking $F$ to
    $\tau_{\leq n-1}\circ F\circ\iota$  induces an oplax monoidal functor
    $$\Fun_\Sigma^\ep(\Cscr_{\geq 0},\Cscr_{\geq
    0})\rightarrow\Fun_\Sigma^\ap(\Cscr_{[0,n-1]},\Cscr_{[0,n-1]}).$$ This
    functor is in fact monoidal and admits a lax monoidal right adjoint $\L_n$.
\end{proposition}

\begin{proof}
    By Lemma~\ref{lem:connectivity} below,
    the assignment $F\mapsto\tau_{\leq n-1}\circ F\circ\iota$ preserves the
    properties of being additively polynomial and preserving sifted colimits,
    so there is an induced functor $\Fun_\Sigma^\ep(\Cscr_{\geq 0},\Cscr_{\geq
    0})\rightarrow\Fun_\Sigma^\ap(\Cscr_{[0,n-1]},\Cscr_{[0,n-1]})$ as claimed,
    and it inherits an oplax monoidal structure from Lemma~\ref{lem:oplax}.
    For monoidality, the
    fact that $\Cscr_{\geq 0}\rightarrow\Cscr_{[0,n-1]}$ is a localization
    implies that the identity functor on $\Cscr_{\geq 0}$ is mapped to the identity functor
    on $\Cscr_{[0,n-1]}$.
    Lemma~\ref{lem:connectivity} below implies that the maps in the
    oplax structure are in fact equivalences. Specifically, if $F$ and $G$ are
    in $\Fun_\Sigma^\ep(\Cscr_{\geq 0})$, then the natural map $$\tau_{\leq
    n-1}\circ(F\circ G)\circ\iota\rightarrow(\tau_{\leq n-1}\circ
    F\circ\iota)\circ(\tau_{\leq n-1}\circ G\circ\iota)$$ is an equivalence by
    Lemma~\ref{lem:connectivity} applied to $G(\iota X)\rightarrow\tau_{\leq
    n-1}G(\iota X)$ for $X\in\Cscr_{[0,n-1]}$.
\end{proof}

The next lemma appears in the case $n=1$ in~\cite[Lem.~4.2.17]{raksit}.

\begin{lemma}\label{lem:connectivity}
    Let $\Cscr_{\geq 0}$ be a compact projectively generated prestable
    $\infty$-category.
    Let $F\in\Fun_\Sigma(\Cscr_{\geq 0},\Cscr_{\geq
    0})$. For any $X\in\Cscr_{\geq 0}$ and any $n\geq 0$, the natural map
    $\tau_{\leq n-1}F(X)\rightarrow\tau_{\leq n-1}F(\tau_{\leq n-1}X)$ is an
    equivalence in the Grothendieck abelian $n$-category $\Cscr_{[0,n-1]}$.
\end{lemma}

\begin{proof}
    Let $\Cscr^0\subseteq\Cscr_{\geq 0}$ be the full subcategory of compact
    projective objects.
    Let $\Ind(\Cscr^0)\subseteq\Cscr_{\geq 0}$ be the full subcategory of
    `flat' objects. Write any $X\in\Cscr_{\geq 0}$ as $|X_\bullet|$
    where $X_\bullet$ is a simplicial object of $\Ind(\Cscr^0)$. The map
    $X\rightarrow\tau_{\leq n-1}X$ is the geometric realization of a map
    $X_\bullet\rightarrow X_\bullet'$, where $X_s\simeq X_s'$ for $s\leq n$. 
    Indeed, by attaching `$s$-cells' for $s>n$, one can kill the homotopy objects of
    $X$ starting in degree $n$. By assumption, $F(X)\we|F(X_\bullet)|$ and
    $F(\tau_{\leq n-1}X)\we|F(X_\bullet')|$. But, the cofiber of
    $F(X)\rightarrow F(\tau_{\leq n-1}X)$ is then computed by
    $|\cofib(F(X_\bullet)\rightarrow F(X_\bullet'))|$. The cofiber
    $\cofib(F(X_s)\rightarrow F(X_s'))$ is nullhomotopic for $s\leq n$, so the
    cofiber is $(n+1)$-connective, the fiber is $n$-connective, and $\tau_{\leq
    n-1}F(X)\rightarrow\tau_{\leq n-1}F(\tau_{\leq n-1}X)$ is an equivalence,
    as desired.
\end{proof}

\begin{definition}[The $n$-derived
    approximation]\label{def:derivedapproximation}
    If $T$ is a monad in
    $\End_\Sigma^\ep(\Cscr_{\geq 0})$, the unit map $T\rightarrow\L_n T$
    of the adjunction from Proposition~\ref{prop:derivedapproximation} is
    called the \longdefidx{$n$-derived approximation}{monad!derived approximation}
    to $T$ and similarly for filtered monads valued in
    $\End_\Sigma^\ep(\Cscr_{\geq 0})$. Unless specified otherwise, the derived
    approximation $\L T=\L_1 T$ of a connective excisively polynomial $\Sigma$-monad will denote the $1$-derived
    approximation; after some discussion in the rest of this section and the next, this is the only case that will occur in the
    rest of this book.
\end{definition}

\begin{remark}
    The derived approximations $\L_n T$ can be viewed as either monads on
    $\Cscr_{\geq 0}$ or on $\Cscr$, an admissible move based on the monoidal
    inclusion $\Fun_\Sigma^\ep(\Cscr_{\geq 0},\Cscr_{\geq
    0})\subseteq\Fun_\Sigma^\ep(\Cscr,\Cscr)$ of Corollary~\ref{cor:connective}.
\end{remark}

\begin{warning}
    Despite the notation,
    if $T$ is a filtered excisively polynomial $\Sigma$-monad, its connective
    covers $\tau_{\geq 0}T$ and their derived approximations $\L_n\tau_{\geq 0}T$ crucially
    depend on the choice the auxiliary $t$-structure $\Cscr_{\geq
    0}\subseteq\Cscr$.
\end{warning}

All in all, there is a diagram of left adjoint monoidal functors as in
Figure~\ref{fig:apparatus}. This diagram will be the apparatus used to pass between
filtered $\Escr$-monads where $\Escr$ varies from among the monoidal
$\infty$-categories above.

\begin{figure}[h]
    \centering
    $$\xymatrix{
        \End_\Sigma^\ep(\Cscr_{\geq
        0})\ar[r]\ar[d]&\cdots\ar[r]&\End_\Sigma^\ap(\Cscr_{[0,n-1]})\ar[r]&\End_\Sigma^\ap(\Cscr_{[0,n-2]})\ar[r]&\cdots\ar[r]&\End_\Sigma^\ep(\Cscr^\heart)\\
        \End_\Sigma^\ep(\Cscr)\ar[d]\\
        \End^\ep(\Cscr)\\
    }$$
    \caption{The displayed arrows are monoidal
    left adjoint functors; the vertical arrows are fully faithful.}
    \label{fig:apparatus}
\end{figure}

\begin{theorem}\label{thm:apparatus}
    Fix a stable presentable $\infty$-category $\Cscr$ with an accessible
    $t$-structure $(\Cscr_{\geq 0},\Cscr_{\leq 0})$ which is right complete and compatible with
    filtered colimits and where $\Cscr_{\geq 0}$ is compact projectively generated.
    Let $T$ be a monad on $\Cscr$. If the sifted colimit preserving
    approximation $T_\Sigma$ is realized as a colimit $\colim\F_\star
    T_\Sigma\we T_\Sigma$ where
    $\F_\star T_\Sigma$ is a filtered excisively polynomial $\Sigma$-monad, then
    there is a diagram
    $$\xymatrix{
        \tau_{\geq 0}\F_\star T_\Sigma\ar[d]\ar[r]&\cdots\ar[r]& \L_n(\tau_{\geq 0}\F_\star
        T_\Sigma)\ar[r]& \L_{n-1}(\tau_{\geq 0}\F_\star T_\Sigma)\ar[r]&
        \cdots\ar[r]& \L_1(\tau_{\geq 0}\F_\star T_\Sigma)\\
        \F_\star T_\Sigma
    }$$
    of filtered monads in $\End_\Sigma^\ep(\Cscr)$ which yields
    a diagram
    $$\xymatrix{
        \tau_{\geq 0}T_\Sigma\ar[d]\ar[r]&\cdots\ar[r]& \L_n \tau_{\geq 0}T_\Sigma\ar[r]& \L_{n-1}
        \tau_{\geq 0}T_\Sigma\ar[r]&
        \cdots\ar[r]& \L_1 \tau_{\geq 0}T_\Sigma\\
        T_\Sigma
    }$$
    of $\Sigma$-monads on $\Cscr$ by taking colimits.
\end{theorem}

\begin{proof}
    The horizontal arrows are obtained by the units of the adjunction of
    Proposition~\ref{prop:derivedapproximation}. The vertical arrow is the
    counit of the adjunction Corollary~\ref{cor:connective}. The colimits exist
    as monads by Lemma~\ref{lem:colimit}.
\end{proof}

\begin{lemma}\label{lem:siftedcolimits}
    Suppose that $\Cscr$ is a $\kappa$-compactly generated presentable $\infty$-category and $T$ is
    a $\Sigma$-monad on $\Cscr$.
    \begin{enumerate}
        \item[{\em (a)}] The forgetful functor $\LMod_T(\Cscr)\rightarrow\Cscr$
            admits a left adjoint (the free $T$-module functor) and preserves
            sifted colimits.
        \item[{\em (b)}] The $\infty$-category $\LMod_T(\Cscr)$ is presentable and $\kappa$-compactly
            generated.
    \end{enumerate}
\end{lemma}

\begin{proof}
    See~\cite[Cor.~4.2.3.5]{ha} for part (a) and see~\cite[Prop.~4.1.10]{raksit} for part (b).
\end{proof}

Theorem~\ref{thm:apparatus} induces a large diagram of left
adjoint functors of module categories. The right adjoints all commute with the forgetful functors
to $\Cscr$.

\begin{corollary}\label{cor:apparatusmodule}
    Fix a stable presentable $\infty$-category $\Cscr$ with an accessible
    $t$-structure $(\Cscr_{\geq 0},\Cscr_{\leq 0})$ which is right complete and compatible with
    filtered colimits and where $\Cscr_{\geq 0}$ is compact projectively generated.
    Let $T$ be a monad on $\Cscr$. If the sifted colimit preserving
    approximation $T_\Sigma$ to $T$ is realized as a colimit $\colim\F_\star
    T_\Sigma\we T_\Sigma$, where
    $\F_\star T_\Sigma$ is a filtered excisively polynomial $\Sigma$-monad,
    then
    there is a commutative diagram of left adjoint functors
    $$\xymatrix{
        \LMod_{\tau_{\geq 0}T_\Sigma}(\Cscr^\heart)\ar[r]&\cdots\ar[r]&\LMod_{\L_n\tau_{\geq 0}T_\Sigma}(\Cscr^\heart)\ar[r]&\cdots\\
        \vdots\ar@{->>}[u]&\cdots&\vdots\ar@{->>}[u]\\
        \LMod_{\tau_{\geq 0}T_\Sigma}(\Cscr_{[0,m-1]})\ar@{->>}[u]\ar[r]&\cdots\ar[r]&\LMod_{\L_n\tau_{\geq 0}T_\Sigma}(\Cscr_{[0,m-1]})\ar@{->>}[u]\ar[r]&\cdots\\
        \vdots\ar@{->>}[u]&\cdots&\vdots\ar@{->>}[u]\\
        \LMod_{\tau_{\geq 0}T_\Sigma}(\Cscr_{\geq 0})\ar@{^{(}->}[d]\ar@{->>}[u]\ar[r]&\cdots\ar[r]&\LMod_{\L_n\tau_{\geq 0}T_\Sigma}(\Cscr_{\geq 0})\ar@{->>}[u]\ar[r]\ar@{^{(}->}[d]&\cdots\\
        \LMod_{\tau_{\geq 0}T_\Sigma}(\Cscr)\ar[d]\ar[r]&\cdots\ar[r]&\LMod_{\L_n\tau_{\geq 0}T_\Sigma}(\Cscr)\ar[r]&\cdots\\
        \LMod_{T_\Sigma}(\Cscr)
    }$$
    where the functors labeled $\hookrightarrow$ are fully faithful and those
    labeled $\twoheadrightarrow$ are localizations (the right adjoints are
    fully faithful).
    Moreover, if $m\leq n$, then the natural functor
    $$\LMod_{\tau_{\geq 0}T_\Sigma}(\Cscr_{[0,m-1]})\rightarrow\LMod_{\L_n \tau_{\geq 0}T_\Sigma}(\Cscr_{[0,m-1]}).$$
    is an equivalence. Finally, the diagram is natural in $\F_\star T_\Sigma$.
\end{corollary}

Modules, algebras, and operads give examples to which Theorem~\ref{thm:apparatus} applies.

\begin{example}
    Suppose that $k$ is a connective $\bE_\infty$-ring for
    and that $A$ is an $\bE_1$-algebra over $k$. Let $T$ be the
    free $A$-module monad on $\Mod_k$, which is an excisively polynomial $\Sigma$-monad.
    Then, $\tau_{\geq 0}T$ is the free $\tau_{\geq 0}A$-module
    monad and $\L_n \tau_{\geq 0}T$ is the free $\tau_{[0,n-1]}A$-module monad. On categories
    of algebras for these monads in $\Mod_k$, one obtains the zig-zag of left adjoint functors
    $$\xymatrix{\Mod_{\tau_{\geq 0}A}\ar[rrr]^{\tau_{[0,n-1]}A\otimes_{\tau_{\geq
    0}A}(-)}\ar[d]_{A\otimes_{\tau_{\geq 0}A}(-)}&&&\Mod_{\tau_{[0,n-1]}A}\\
    \Mod_A.}$$
\end{example}

\section{Derived $\Oscr$-algebras}\label{sec:derived_algebras}

We discuss the machinery of the previous section in the context of monads associated to operads. We
work with operads internal to $\Cscr$, i.e.,
algebra objects in the monoidal $\infty$-category of symmetric sequences in $\Cscr$.
For details, see~\cite[Sec.~4.1]{raksit}.

\begin{lemma}\label{lem:operads}
    Let $\Cscr$ be a presentably symmetric monoidal stable
    $\infty$-category and let $\Oscr$ be an operad in $\Cscr$.
    The free $\Oscr$-monad associated to $\Oscr$ naturally admits the structure of a filtered
    excisively polynomial $\Sigma$-monad.
\end{lemma}

\begin{proof}
    Follow the proof of~\cite[Const.~4.1.8]{raksit}. Specifically, Raksit
    constructs a functor from $\Op(\Cscr)$, the $\infty$-category of operads in
    $\Cscr$, to the $\infty$-category of lax symmetric monoidal functors
    $\bZ_{\geq 0}^\times\rightarrow\mathrm{SSeq}(\Cscr)$ with respect to the
    composition monoidal structure. Given a symmetric sequence $\Oscr$ the
    resulting filtered object $\psi(\Oscr)$ has $$\psi(\Oscr)_d(i)\we\begin{cases}
        \Oscr(i)&\text{if $i\leq d$,}\\
        \emptyset&\text{otherwise.}
    \end{cases}$$ There is a monoidal functor
    $\mathrm{SSeq}(\Cscr)\rightarrow\End(\mathrm{SSeq}(\Cscr))$ given by
    tensoring (via the composition product) with a symmetric sequence. The
    resulting endomorphisms all canonically preserve the full subcategory
    $\Cscr\subseteq\mathrm{SSeq}(\Cscr)$ obtained by left Kan extending along
    the inclusion $\ast\xrightarrow{\emptyset}\bF$. (Here, $\bF$ is the
    groupoid of finite sets and bijections.) The composition
    $\Op(\Cscr)\rightarrow\Alg(\End(\mathrm{SSeq}(\Cscr)))\rightarrow\Alg(\End(\Cscr))$
    takes an operad $\Oscr$ to the associated free $\Oscr$-algebra monad.
    Via the $\psi$ functor above, the operad $\Oscr$ upgrades to a lax monoidal
    functor $\bZ_{\geq 0}^\times\rightarrow\mathrm{SSeq}(\Cscr)$ and the
    composition of monoidal functors
    $\mathrm{SSeq}(\Cscr)\rightarrow\End(\mathrm{SSeq}(\Cscr))\rightarrow\End(\Cscr)$
    is used to obtain a filtered monad whose underlying filtered endofunctor is of the form
    $$\mathrm{Free}_{\Oscr,\leq\star}\we\bigoplus_{i\leq\star}(\Oscr(i)\otimes(-)^{\otimes
    i})_{\Sigma_i}.$$
    Each functor $(\Oscr(i)\otimes(-)^{\otimes i})_{\Sigma_i}$ is excisively polynomial and
    preserves sifted colimits, which completes the proof.
\end{proof}

\begin{definition}[$n$-derived $\Oscr$-algebras]
    Fix a stable presentable $\infty$-category $\Cscr$ with an accessible
    $t$-structure $(\Cscr_{\geq 0},\Cscr_{\leq 0})$ which is right complete and compatible with
    filtered colimits and where $\Cscr_{\geq 0}$ is compact projectively generated.
    An \longdefidx{$n$-derived $\Oscr$-algebra}{derived $\Oscr$-algebra} is a left module for the monad $\L_n\mathrm{Free}_{\Oscr}$.
    Write $\Alg_{\L_n\Oscr}(\Cscr)$ for $\LMod_{\L_n\mathrm{Free}_\Oscr}(\Cscr)$, the
    $\infty$-category of $n$-derived $\Oscr$-algebras. When $n=1$, these will be called simply
    derived $\Oscr$-algebras.
\end{definition}

\begin{example}\label{ex:derivedem}
    Let $k$ be a connective $\bE_\infty$-ring, such as the sphere spectrum
    $\bS$ or the ring of integers $\bZ$. In this case, the standard
    $t$-structure on $\Mod_k$ is accessible, right complete, and compatible with
    filtered colimits and the prestable $\infty$-category $\Mod_k^\cn$ of
    connective objects is compact projectively generated (by the compact projective
    $k$-modules). Fix $1\leq r\leq\infty$ and consider the
    $\bE_r$-operad internal to $\Mod_k$ (given by taking the $k$-chains of the
    spaces in the little $r$-cubes operads). Write $\bE_r$ for the free
    $\bE_r$-algebra monad. By Lemma~\ref{lem:operads}, $\bE_r$ admits the
    structure of a filtered excisively polynomial $\Sigma$-monad. It is already
    connective. In other words,
    $$\tau_{\geq 0}(\bE_r)_\Sigma\we(\bE_r)_\Sigma\we\bE_r.$$
    On the other hand, the tower of $n$-derived approximations
    $$\bE_r\rightarrow\cdots\rightarrow\L_n\bE_r\rightarrow\L_{n-1}\bE_r\rightarrow\cdots\rightarrow\L_1\bE_r$$
    provides an interpolation between the classical notion of $\bE_r$-algebra
    in $\D(k)$ and the notion of a (1-)derived $\bE_r$-algebra.
    For $r\geq 2$ and  $n\leq r-2$,
    $\L_n\bE_r\we\L_n\bE_\infty$ using that the $i$th anima $\bE_r(i)$ is
    equivalent to the homotopy type of the configuration space of $i$ points in $\bR^r$ and that the
    bottom non-trivial homology group of this space occurs in degree $r-1$;
    see~\cite[Cor.~2.1]{fadell-neuwirth}.

    The maps of monads above induce commutative squares of left adjoint
    functors
    $$\xymatrix{
        \Alg_{\bE_r}(\Mod_k^\cn)\ar[r]\ar[d]^{\tau_{\leq r-1}}&\Alg_{\L_n\bE_r}(\Mod_k^\cn)\ar[d]^{\tau_{\leq r-1}}\\
        \Alg_{\bE_r}((\Mod_k)_{[0,r-1]})\ar[r]^\we&\Alg_{\L_n\bE_r}((\Mod_k)_{[0,r-1]})\\
    }$$
    where the bottom arrow is an equivalence.

    The free $\L_n\bE_r$-algebra on a flat $k$-module
    $M$ is $\tau_{\leq n-1}\bE_r(M)$, i.e., the $(n-1)$-truncation of the
    free $\bE_r$-algebra on $M$. In particular,
    every $\L_n\bE_r$-algebra is canonically an $\L_n\bE_r$-algebra (and
    hence an $\bE_r$-algebra) over $\tau_{\leq n-1}k$.
\end{example}

\begin{example}
    If $k$ is $(n-1)$-truncated, then
    $\bE_1\we\L_m\bE_1$ for $m\geq n$. For $k$ a static commutative ring,
    $\bE_1\we\L_1\bE_1$. 
\end{example}

\begin{example}
    If $k$ is a static commutative ring, then the natural map
    $\L_1\bE_n\rightarrow\L_1\bE_\infty$ is an equivalence for $n\geq 2$.
    Indeed, in $\Mod_k^\heart$, there is no distinction between an
    $\bE_n$-algebra and a commutative algebra.
\end{example}

\section{Derived commutative rings}\label{sec:dalg}

We introduce derived commutative rings and prove a basic functoriality result
(Proposition~\ref{prop:functoriality}).

\begin{definition}[Condition $(*)$]
    Let $\Cscr$ be a symmetric monoidal stable $\infty$-category with a $t$-structure $(\Cscr_{\geq
    0},\Cscr_{\leq 0})$. We say that $\Cscr$ with its symmetric monoidal structure and
    $t$-structure satisfies condition $(*)$ if the tensor product on $\Cscr$
    commutes with small colimits in each variable and if the $t$-structure is
    accessible, right complete, compatible with filtered colimits, compatible with the symmetric
    monoidal structure on $\Cscr$ and if, moreover, $\Cscr_{\geq 0}$ is compact projectively generated.
    We will say that a stable $\infty$-category $\Cscr$ satisfies $(*)$ if it is equipped with a
    symmetric monoidal structure and a $t$-structure such that the triple
    $(\Cscr,\Cscr^\otimes,(\Cscr_{\geq 0},\Cscr_{\leq 0}))$ satisfies $(*)$.
\end{definition}

This is a weakening of Raksit's notion of derived algebraic context, for which see
Definition~\ref{def:context} below.

\begin{definition}[Derived commutative rings]\label{def:derivedcommutative}
    Suppose that $\Cscr$ is a stable $\infty$-category satisfying $(*)$.
    There is a
    derived commutative ring monad $\LSym_\Cscr=\L_1\bE_\infty$. Write $\DAlg_\Cscr$ for
    $\Alg_{\L_1\bE_\infty}(\Cscr)=\LMod_{\LSym_\Cscr}(\Cscr)$; the objects of $\DAlg_\Cscr$ are
    called derived commutative rings
    in $\Cscr$; there is a forgetful functor from $\DAlg_\Cscr$ to $\CAlg_\Cscr$, the
    $\infty$-category of $\bE_\infty$-algebras in $\Cscr$; the forgetful functor preserves all limits
    and sifted colimits. Under additional hypotheses, it preserves all colimits by Proposition~\ref{prop:colimits} below.
    In any case, $\DAlg_\Cscr$ is presentable by Lemma~\ref{lem:siftedcolimits}.
\end{definition}

\begin{example}
    If $k$ is a commutative ring, then $\Mod_k$ with its usual symmetric monoidal
    structure and $t$-structure satisfies $(*)$.
    We write $\LSym_k$ for $\LSym_{\Mod_k}=\L_1\bE_\infty$.
    A \longdefidx{derived commutative $k$-algebra}{derived commutative
    ring!derived commutative $k$-algebra} is an $\LSym_k$-algebra over $k$.
    The $\infty$-category of derived commutative $k$-algebras is written $\DAlg_k$.
    The full subcategory
    $\DAlg_k^\cn\subseteq\DAlg_k$ consisting of connective derived
    commutative $k$-algebras is equivalent to
    $\CAlg^{\mathrm{an}}_k=\s\CAlg^\heart_k[W^{-1}]$, the $\infty$-category of
    animated commutative $k$-algebras, which is also realized concretely as the
    $\infty$-categorical localization of the category of simplicial commutative
    $k$-algebras at the weak equivalences. Note that a static commutative $k$-algebra
    is canonically a derived commutative ring. It follows that any limit or colimit
    of static commutative $k$-algebras canonically admits the structure of a
    derived commutative ring. Particularly interesting cases include
    $\bZ^X=\lim_X\bZ$ when $X$ is an anima, which satisfies $\pi_i(\bZ^X)\iso\H^{-i}(X;\bZ)$, or $\R\Gamma(X,\Oscr)$ when $X$ is a
    site and $\Oscr$ is a sheaf of commutative rings on $X$ (or even a sheaf of
    animated commutative rings). These examples are striking because they are known not to
    admit cdga models in general thanks to the behavior of Steenrod operations.
\end{example}

\begin{example}
    If $k$ is a commutative $\bQ$-algebra, then
    $\Sym^{\leq\star}_{\Mod_k}\rightarrow\LSym^{\leq\star}_{\Mod_k}$ is a
    filtered equivalence of filtered monads, so the theories of derived commutative $k$-algebras and
    $\bE_\infty$-$k$-algebras agree; they both recover the homotopy theory of
    commutative differential graded $k$-algebras and quasi-isomorphisms.
\end{example}

\begin{example}\label{ex:low_free}
    Here are some computations of homotopy rings of
    free derived commutative $k$-algebras on
    generators of small degree:
    \begin{enumerate}
        \item[(0)] $\pi_*\LSym_k(k[0])\iso k[x]$ where $|x|=0$;
        \item[(1)] $\pi_*\LSym_k(k[1])\iso k[y]/(y^2)$ where $|y|=1$;
        \item[(2)] $\pi_*\LSym_k(k[2])\iso k\langle z\rangle$, the free
            divided power algebra on $z$ where $|z|=2$.
    \end{enumerate}
\end{example}

\begin{example}\label{ex:f2_free}
    Paul Goerss summarized in~\cite{goerss_f2} the computation of the homotopy ring of $\LSym_{\bF_2}(\bF_2[n])$ for
    all $n\geq 0$ following the work of several previous researchers.
    These are in terms of ``higher divided power'' operations $\delta_i\colon\pi_n
    R\rightarrow\pi_{n+i}R$ which exist for $n\geq 1$ and
    $2\leq i\leq n$ on the homotopy groups of any derived commutative $\bF_2$-algebra $R$.
    By composing, one can construct operations $\delta_I=\delta_{i_1}\cdots\delta_{i_k}$.
    Such a sequence $I$ is admissible if $i_s\geq 2i_{s+1}$ for all $s$.
    The excess of $I$ is $i_1-i_2-\cdots-i_k$. Goerss shows that for $n\geq 1$ there is an
    isomorphism $$\pi_*\LSym_{\bF_2}(\bF_2[n])\iso\Lambda_{\bF_2}(\delta_I(\iota_n)\colon\text{$I$
    admissible of excess at most $n$}),$$
    where $\iota_n\in\pi_n\LSym_{\bF_2}(\bF_2[n])$ is the canonical class.
    See also the proof of~\cite[Thm.~7.14]{brantner-mathew} which explains the work of
    Nakaoka~\cite{nakaoka} and Priddy~\cite{priddy}.
\end{example}

\begin{example}
    If $k$ is a commutative
    $\bF_p$-algebra, then one can compute $\pi_*\LSym_k(k[-1])$ as the group
    cohomology of $\Ga$ with coefficients in $k$. This follows from the
    computation of the space $\B\Ga(R)\we
    R[1]\we\Map_{\DAlg_k}(\LSym_k(k[-1]),R)$ for $R$ an animated commutative
    ring. It follows that  the stack $\B\Ga$ is an affine stack in the sense of
    To\"en and corresponds to 
    $\Spec\LSym_k(k[-1])$; see~\cite[Lem.~2.2.5]{toen-affines}. In particular,
    $\R\Gamma(\B\Ga,\Oscr)\we\LSym_k(k[-1])$. The calculation of the
    $\Ga$-cohomology of $k$ is achieved by Jantzen
    in~\cite[4.27]{jantzen} who shows that if $p$ is odd, then
    $\pi_*\LSym_k(k[-1])$ is polynomial on countably many degree $-2$
    generators and exterior on countably many degree $-1$ generators.
    When $p=2$, Jantzen obtains $\pi_*\LSym_k(k[-1])$ is polynomial on countably many degree $-1$
    generators.
\end{example}

The following lemma will be used heavily in Sections~\ref{sec:grdalg} and~\ref{sec:crystallization}.

\begin{lemma}\label{lem:toramplitude}
    Let $k$ be a commutative ring.
    If $P$ is a perfect complex of $k$-modules with Tor-amplitude in
    $[a,b]$ with $a\geq 0$, then $\LSym^r(P)$ is perfect with Tor-amplitude in
    $[\max(a+2r-2,0),rb]$. In particular, when $k=\bZ$, $\H_{rb}(P)$ is torsion-free.
\end{lemma}

\begin{proof}
    The functor $Q\mapsto Q^{\otimes r}$ is additively polynomial of degree $\leq r$ as a functor
    on projective $k$-modules. Using that $Q\mapsto (Q^{\otimes r})_{\Sigma_r}$ is a
    quotient of the functor $Q\mapsto Q^{\otimes r}$ implies that its cross effects are quotients
    of the cross effects of $Q\mapsto Q^{\otimes r}$, so that $Q\mapsto (Q^{\otimes r})_{\Sigma_r}$
    is additively polynomial of degree $\leq r$. The upper bound of the Tor-amplitude of
    $rb$ in the first part of the lemma now follows from~\cite[Lem.~3.3]{gillet-soule} or~\cite[Hilfssatz~4.23]{dold-puppe}.
    The fact about torsion-freeness applies to any perfect complex with the given Tor-amplitude.
    The lower bound for the Tor-amplitude follows from using the d\'ecalage equivalences
    $\LSym^r(P[1])\we\L\Lambda(P)[r]$ and $\L\Lambda^r(P[1])\we\L\Gamma^r(P)[r]$; see for
    instance~\cite[I.4.3.2.1]{illusie-cotangent-1}.
\end{proof}

We also need the upper bounds for derived alternating and antisymmetric powers.
If $k^\epsilon$ denotes the sign representation of $\Sigma_r$,
$\LAntiSym^r$ denotes the derived functor of $P\mapsto(k^\epsilon\otimes P^{\otimes
r})_{\Sigma_r}$. Note also that $\L\Lambda^r$ is the derived functor of
$P\mapsto(k^\epsilon\otimes P^{\otimes r})^{\Sigma_r}$. Here, $P$ is a finitely
presented project $k$-module and both orbits and fixed points
are computed in $\Mod_k^\heart$.

\begin{lemma}\label{lem:toramplitude_alternating}
    Let $k$ be a commutative ring.
    If $P$ is a perfect complex of $k$-modules with Tor-amplitude in $[a,b]$ with $a\geq 0$, then
    $\LAntiSym^r(P)$ and $\L\Lambda^r(P)$ are perfect with Tor-amplitudes bounded above by $rb$.
\end{lemma}

\begin{proof}
    The case of $\L\Lambda^r(P)$ follows from Lemma~\ref{lem:toramplitude} and d\'ecalage.
    For $\LAntiSym^r$, one argues as in the beginning of the proof of Lemma~\ref{lem:toramplitude}
    using its description as a quotient.
\end{proof}

\begin{example}
    If $R$ is
    a derived commutative $\bZ$-algebra, then $\pi_*R$ is a
    graded-commutative ring since this is true for $\bE_\infty$-rings.
    However, more can be said: $\pi_*R$ is strictly graded-commutative in positive
    degrees: if
    $x\in\pi_{2s+1}R$ for some $s\geq 1$, then $x^2=0$ in $\pi_{4s+2}R$. For $s=0$, this is easy to
    see using the fact that $\LSym^2(\bZ[1])\simeq(\Lambda^2\bZ)[2]\simeq 0$. For
    $s>0$, Akhil Mathew pointed out the following simple argument.
    It is enough to see that $\LSym^2(\bZ[2s+1])$ has no torsion in degree
    $4s+2$, but this follows from Lemma~\ref{lem:toramplitude}.

    Note however that if $x\in\pi_{2s+1}R$ and $s<0$, then it need not be the case
    that $x^2=0$. For example, consider the case of the Eilenberg--Mac Lane anima
    $K(\bZ,3)$ and take $R=\bZ^{K(\bZ,3)}=\lim_{K(\bZ,3)}\bZ$, which is a derived
    commutative ring. The tautological class $x\in\pi_{-3}R\iso\H^3(K(\bZ,3),\bZ)\iso\bZ$ has non-zero
    square, which is a generator of $\pi_{-6}R\iso\H^6(K(\bZ,3),\bZ)\iso\bZ/2$.
\end{example}

\begin{definition}
    Let $\Cscr$ be a presentably symmetric monoidal $\infty$-category. Say that a
    filtered monad $\F_\star T$ with a filtered monad map
    $\Sym^{\leq\star}\rightarrow\F_\star T$ is
    \longdefidx{binomial}{monad!binomial} if $\F_d T$ naturally decomposes into a
    direct sum $\F_dT=\bigoplus_{0\leq e\leq d} T^e$ for functors $\{T^e\}_{e\in\bN}$ compatible with the
    corresponding direct sum decomposition of $\Sym^{\leq d}$ and if the filtered
    multiplication maps induce equivalences $\bigoplus_{0\leq e\leq
    d}(T^e(M)\otimes_{\1_\Cscr}T^{d-e}(N))\we T^d(M\oplus N)$ for all $d\geq 0$ and
    all objects $M,N$ in $\Cscr$.
\end{definition}

The following proposition was observed by Bhatt and Mathew and proven under
slightly stronger hypotheses in~\cite[Prop.~4.2.27]{raksit}.

\begin{proposition}\label{prop:colimits}
    Suppose that $\Cscr$ is a stable $\infty$-category which satisfies $(*)$.
    Assume that the unit object $\1_\Cscr$ is static.
    Let $\Sym^{\leq\star,\heart}$ denote the filtered symmetric
    algebra monad on $\Cscr^\heart$ and suppose that $\Sym^{\leq\star,\heart}\rightarrow\F_\star T$
    is a map of additively polynomial $\Sigma$-monads on $\Cscr^\heart$ where
    $\F_\star T$ is binomial. If additionally $T(0)\we\1_\Cscr$
    and for each compact projective $P$ in $\Cscr^\heart$, the object $T^d(P)$
    is flat over $\1_\Cscr$ for all $d\geq 0$,
    then the forgetful functor
    $$\CAlg_\Cscr\leftarrow\LMod_T(\Cscr)$$ preserves all limits and colimits.
\end{proposition}

\begin{proof}
    The functor in question preserves limits since it is a right adjoint. It
    preserves sifted colimits by Lemma~\ref{lem:siftedcolimits}. The rest of the proof follows that
    of Raksit in~\cite[Prop.~4.2.27]{raksit}, which proves the case that $T$ is the
    $\LSym_\Cscr$-monad and $\Cscr_{\geq 0}$ defines a derived algebraic context.
    Specifically, the assumption that $T(0)\we\1_\Cscr$ (which means in particular that $\1_\Cscr$ with its
    $T$-module structure is the initial object of $\LMod_T(\Cscr)$) implies that the forgetful functor preserves
    the unit objects. Arguing as in the proof of~\cite[Prop.~4.2.27]{raksit}, one reduces to
    checking that the
    natural map $T(P)\otimes_{\1_\Cscr}T(P)\rightarrow T(P\oplus Q)$ 
    is an equivalence for $P$ and $Q$ finitely presented projective objects in
    $\Cscr_{\geq 0}$. Thus, it is enough to show that $T^d(P\oplus
    Q)\we\bigoplus_{e=0}^n T^e(P)\otimes_{\1_\Cscr}T^{d-e}(Q)$ for all $d\geq 0$, which
    follows from the binomial condition on $\F_\star T$ as a monad on
    $\Cscr^\heart$ together with the flatness hypothesis.
\end{proof}

\begin{corollary}
    Suppose that $\Cscr$ is a stable $\infty$-category which satisfies $(*)$.
    If the unit object $\1_\Cscr$ is static
    and if for each compact projective $P$ in $\Cscr^\heart$, the object $\pi_0\Sym^d(P)$
    is flat over $\1_\Cscr$ for all $d\geq 0$,
    then the forgetful functor
    $$\CAlg_\Cscr\leftarrow\DAlg_\Cscr$$ preserves all limits and colimits.
\end{corollary}

In~\cite{raksit}, the construction of derived commutative rings is carried out under slightly
stricter conditions and it will be convenient to refer to these conditions later.

\begin{definition}\label{def:context}
    A \defidx{derived algebraic context} is a presentably symmetric monoidal
    stable $\infty$-category $\Cscr$ equipped with a $t$-structure
    $(\Cscr_{\geq 0},\Cscr_{\leq 0})$ satisfying the following properties:
    \begin{enumerate}
        \item[(a)] the $t$-structure is right complete and compatible with
            filtered colimits;
        \item[(b)] the $t$-structure is compatible with the symmetric monoidal
            structure on $\Cscr$ in that the unit object of $\Cscr$ is
            connective and tensor products of connective objects are
            connective;
        \item[(c)] $\Cscr_{\geq 0}\we\Pscr_\Sigma(\Cscr^0)$ where
            $\Cscr^0\subseteq\Cscr^\heart$ is the full subcategory of compact
            projective objects;
        \item[(d)] $\Cscr^0\subseteq\Cscr^\heart$ is a symmetric monoidal
            subcategory which is closed under symmetric powers.
    \end{enumerate}
\end{definition}

\begin{remark}
    Any derived algebraic context $\Cscr$ comes with a canonical $\bZ$-linear
    structure, i.e., a symmetric monoidal left adjoint functor
    $\Mod_\bZ\rightarrow\Cscr$. Moreover, $\Cscr$ satisfies $(*)$, so one can speak
    of derived commutative rings in a derived algebraic context.
\end{remark}

If $\Cscr$ is a derived algebraic context,
write for a moment $\DAlg^{\mathrm{R}}_\Cscr$ for the $\infty$-category of derived commutative
rings in $\Cscr$ as constructed in~\cite[Sec.~4.2]{raksit}.

\begin{lemma}
    If $\Cscr$ is a derived algebraic context, then $\DAlg^{\mathrm{R}}_\Cscr\we\DAlg_\Cscr$.
\end{lemma}

\begin{proof}
    This follows by unwinding the definitions and using (d) in the definition of a derived
    algebraic context.
\end{proof}

We will therefore treat our construction as a generalization of Raksit's and omit the distinction
in the notation.

Finally, we come to main payoff of the work in Section~\ref{sec:polynomial_monad}.
The following result will allow us to construct comparison functors between $\LSym$-monads arising
from different stable $\infty$-categories satisfying $(*)$.

\begin{proposition}\label{prop:functoriality}
    Suppose that $\Cscr$ and $\Dscr$ are stable $\infty$-categories satisfying $(*)$ and that
    $F\colon\Cscr\rightarrow\Dscr$ is an exact functor such that
    \begin{enumerate}
        \item[{\em (a)}] $F$ preserves colimits,
        \item[{\em (b)}] the right adjoint $U$ preserves filtered colimits (or, equivalently, $F$ takes
            compact objects to compact objects),
        \item[{\em (c)}] $F$ is right $t$-exact,
        \item[{\em (d)}] $U(\LSym_\Dscr^{\leq\star}(F(P)))$ is coconnective for every compact projective
            generator $P$ of $\Cscr_{\geq 0}$, and
        \item[{\em (e)}] $F$ is symmetric monoidal.
    \end{enumerate}
    Then there is an induced left adjoint
    functor $\DAlg_\Cscr\xrightarrow{F}\DAlg_\Dscr$ fitting into a commutative diagram
    $$\xymatrix{
        \Cscr\ar[r]^F\ar[d]^{\LSym_\Cscr}&\Dscr\ar[d]^{\LSym_\Dscr}\\
        \DAlg_\Cscr\ar[r]^F&\DAlg_\Dscr.
    }$$
\end{proposition}

\begin{proof}
    To prove this, consider the filtered excisively polynomial monad $T^{\leq\star}=U\LSym_\Dscr^{\leq\star} F$ on $\Cscr$
    constructed using
    Lemma~\ref{lem:oplax}. There is a map of filtered monads $U\Sym_\Dscr^{\leq\star} F\rightarrow
    T^{\leq\star}$ by functoriality and another map of filtered excisively polynomial monads
    $\Sym_\Cscr^{\leq\star}\rightarrow U\Sym_\Dscr^{\leq\star} F$ by symmetric monoidality of $F$.
    The latter map factors through $\tau_{\geq 0}(U\Sym_\Dscr^{\leq\star} F)\rightarrow
    U\Sym_\Dscr^{\leq\star}F$ since $\Sym_\Cscr^{\leq\star}$ is connective and hence there is a
    commutative diagram
    $$\xymatrix{
        \Sym^{\leq\star}_\Cscr\ar[r]&\tau_{\geq
        0}(U\Sym_\Dscr^{\leq\star}F)\ar[d]\ar[r]&U\Sym_\Dscr^{\leq\star}F\ar[d]\\
        &\tau_{\geq 0}(U\LSym_\Dscr^{\leq\star}F)\ar[r]&U\LSym^{\leq\star}_\Dscr F
    }$$
    of filtered excisively polynomial monads. By assumption (d), $\tau_{\geq
    0}(U\LSym_\Dscr^{\leq\star}F)$ takes static values on compact projective generators of
    $\Cscr_{\geq 0}$, which implies that $\Sym^{\leq\star}_\Cscr\rightarrow\tau_{\geq
    0}(U\LSym_\Dscr^{\leq\star}F)$ factors through
    $\Sym^{\leq\star}_\Cscr\rightarrow\LSym^{\leq\star}_\Cscr$ by
    Proposition~\ref{prop:derivedapproximation}.

    In general, if $F\colon\Cscr\rightleftarrows\Dscr\colon U$ is an adjunction where $U$
    preserves sifted colimits and $T$ is a
    sifted colimit preserving
    monad on $\Dscr$, then there is an induced adjunction
    $\LMod_{UTF}(\Cscr)\rightleftarrows\LMod_T(\Dscr)$ which is compatible with forgetful functors in
    the obvious sense and where the right adjoint preserves sifted colimits. In the present situation, this adjunction can be precomposed with the
    adjunction
    $\LMod_{\LSym_\Cscr}(\Cscr)\rightleftarrows\LMod_{U\LSym_\Dscr F}(\Cscr)$
    arising from the map of monads constructed above. The composition of the left adjoints is the
    desired $F$ from the statement of the proposition.
\end{proof}

\begin{remark}
    A special case of Proposition~\ref{prop:functoriality} is the situation of a symmetric monoidal left adjoint functor $\Cscr\rightarrow\Dscr$
    of derived algebraic contexts which is $t$-exact and takes $\Cscr^0$ to $\Dscr^0$.
\end{remark}

It is necessary to consider situations where there is an $\LSym$ monad on
modules or comodules even when this monad does not arise from the process of deriving monads
considered above.

\begin{definition}[Modules]\label{def:derivedonmodules}
    If $\Cscr$ is a stable $\infty$-category satisfying $(*)$ and $R\in\DAlg_\Cscr$, then the
    forgetful functor $(\DAlg_\Cscr)_{R/}\rightarrow\DAlg_\Cscr$ preserves limits and sifted
    colimits, the latter because, in general, the limit of a constant $I$-indexed diagram on $x$ is
    equivalent to $x$ if $I$ is cosifted. It follows that $(\DAlg_\Cscr)_{R/}$ is monadic over $\Mod_R(\Cscr)$. We let
    $\LSym_R$ denote the associated monad. The left adjoint
    $\DAlg_\Cscr\rightarrow(\DAlg_\Cscr)_{R/}$ is given by the coproduct with $R$. In good cases,
    such as when Proposition~\ref{prop:colimits} applies, this coproduct is computed as the tensor product
    with $R$. We will only apply this construction in those cases.
\end{definition}

\begin{remark}
    Even if $R\in\DAlg_\Cscr$ is connective, so that $\Mod_R(\Cscr)$ satisfies $(*)$,
    the monad $\LSym$-monad constructed above on $\Mod_R(\Cscr)$ might
    differ from the $\LSym_{\Mod_R(\Cscr)}$ monad constructed using derived approximations.
    Indeed, suppose that $R$ is an $(n-1)$-truncated connective derived commutative $k$-algebra for
    some $n\geq 2$ and some commutative ring $k$. Then, $\pi_\ast\LSym_R(R)\iso(\pi_\ast R)[t]$, the
    flat polynomial algebra on a degree $0$ class over the graded-commutative ring $\pi_\ast R$.
    On the other hand, the best possible guess for how to obtain this $\LSym$ from
    Theorem~\ref{thm:apparatus} is to take $\L_n\bE_\infty$. However, the value of this monad on $R$
    itself is the free $\bE_\infty$-algebra on a degree $0$ generator in the symmetric monoidal
    abelian $n$-category $(\Mod_R)_{[0,n-1]}$, which is the $(n-1)$-truncation of the free
    $\bE_\infty$-algebra on $R$ and hence can have contributions from the homology of the
    symmetric groups for $n\geq 2$.
\end{remark}

The more important case for us will be the case of comodules over a derived bicommutative bialgebra.

\begin{definition}[Comodules]\label{def:comodules}
    If $\Cscr$ is a stable $\infty$-category satisfying $(*)$, a \longdefidx{derived commutative
    bialgebra}{derived commutative ring!derived commutative bialgebra} is an $\bE_1$-coalgebra $A$ in $\DAlg_\Cscr$; a
    \longdefidx{derived bicommutative bialgebra}{derived commutative
    ring!derived bicommutative bialgebra} is an $\bE_\infty$-coalgebra $A$
    in $\DAlg_\Cscr$. Let
    $\cDAlg_\Cscr$ denote the $\infty$-category of derived commutative
    bialgebras. If $A$ is a derived commutative bialgebra, let
    $\cMod_A(\Cscr)$ denote the $\infty$-category of $A$-comodules in $\Cscr$
    and let $\cMod_A(\DAlg_\Cscr)$ denote the $\infty$-category of
    $A$-comodules in $\DAlg_\Cscr$. Since
    $\LSym_\Cscr\colon\Cscr\rightarrow\DAlg_\Cscr$ is naturally oplax
    symmetric monoidal, one obtains a monadic adjunction
    $$\LSym_\Cscr\colon\cMod_A(\Cscr)\rightleftarrows\cMod_A(\DAlg_\Cscr)\colon\text{forget}.$$
    The forgetful functors fit into a commutative diagram
    $$\xymatrix{
        \cMod_A(\Cscr)\ar[r]^{\LSym_\Cscr}\ar[d]&\cMod_A(\DAlg_\Cscr)\ar[d]\\
        \Cscr\ar[r]^{\LSym_\Cscr}&\DAlg_\Cscr.
    }$$
    See~\cite[Prop.~4.2.32]{raksit} for details.
\end{definition}

\begin{example}[Group actions]
    The picture the reader should have in mind is the following. If $G$ is a
    finite group, then the Hopf algebra $\Oscr(G)=\bZ[G]^\vee$ of functions on $G$ is a
    commutative bialgebra. The objects of $\cMod_{\Oscr(G)}(\Mod_\bZ)$ are
    $\bZ$-modules with a $G$-action and similarly the objects of
    $\cMod_{\Oscr(G)}(\DAlg_{\bZ})$ are derived commutative rings with a
    $G$-action. The tensor product in each case is via that of the underlying
    action; this is where the comultiplication on $\Oscr(G)$ is used. Intuitively,
    it is clear that if $M$ is a $\bZ$-module spectrum with a $G$-action, then $\LSym_{\bZ}(M)$
    is a derived commutative ring with a $G$-action. The diagram above makes this
    precise.
\end{example}

\section{Graded derived commutative rings}\label{sec:grdalg}

Let $k$ be a commutative ring and let $\Mod_k=\Mod_k(\Sp)$ be the derived $\infty$-category of $k$.
The stable $\infty$-category $\GrMod_k$ admits a rich
structure, including infinitely many different $t$-structures to which the methods of
the previous sections apply. Much of the discussion below applies to more general symmetric
monoidal stable $\infty$-categories. Specifically, we could replace $\Mod_k$
by any derived algebraic context. However, the operation of shearing recalled in
Proposition~\ref{prop:shear} is not symmetric monoidal over $\Mod_\bS$ when $\bS$ denotes the sphere
spectrum.

\begin{definition}[Graded derived $\infty$-categories]
    The \defidx{graded derived $\infty$-category} of a commutative ring $k$ is the functor category
    $$\GrMod_k=\Fun(\bZ^\delta,\Mod_k)\we\prod_{n\in\bZ}\Mod_k,$$ where $\bZ^\delta$ denotes the abelian group $\bZ$ viewed as
    a symmetric monoidal category with no non-identity morphisms. The objects $X^\ast$ of $\GrMod_k$ are called
    \defidx{graded $k$-modules}, with $X^n$ the \longdefidx{weight $n$
    piece}{graded $k$-module!weight $n$ piece} of $X^*$. Addition in $\bZ$ endows
    $\GrMod_k$ with the structure of a symmetric monoidal $\infty$-category under Day
    convolution~\cite{glasman,ha}:
    given graded objects $X^\ast$ and $Y^\ast$, the tensor product $(X\otimes_kY)^\ast$ has
    weight $n$ piece equivalent to $(X\otimes_kY)^n\we\bigoplus_{i+j=n}X^i\otimes_kY^j$.
    Evaluation at $n\in\bZ^\delta$ defines a functor $\GrMod_k\rightarrow\Mod_k$ which preserves
    all limits and colimits. Its left and right adjoints are equivalent and will be denoted by
    $M\mapsto M(n)$, the graded $k$-module which is $M$ in weight $n$ and $0$ elsewhere.
    The unit object is $k(0)$, denoted simply by $k$.
\end{definition}

\begin{definition}[The neutral $t$-structure]
    The $\infty$-category $\GrMod_k$ is equivalent to the derived $\infty$-category of the abelian category
    $\Gr(\Mod_k^\heart)=\Fun(\bZ^\delta,\Mod_k^\heart)$ of graded static $k$-modules. There is in particular a
    $t$-structure on $\GrMod_k$ with connective objects the full subcategory $(\GrMod_k)_{\geq 0}^\N$ of graded objects
    $X^\ast$ where $X^n\in(\Mod_k)_{\geq 0}$ for all $n\in\bZ$. This is the pointwise $t$-structure
    associated to the standard $t$-structure on $\Mod_k$ and the presentation of $\GrMod_k$ as a
    functor category. We will call it the \longdefidx{neutral $t$-structure}{$t$-structure!neutral} or
    \longdefidx{$[0]$-sheared neutral $t$-structure}{$t$-structure!$[0]$-sheared} on $\GrMod_k$. It is
    compatible with the symmetric monoidal structure on $\GrMod_k$ and $(\GrMod_k)^\N_{\geq 0}$ is
    compact projectively generated. The heart
    $\GrMod_k^{\N\heart}$ is equivalent to the abelian category $\Gr(\Mod_k^\heart)$ of graded static $k$-modules. As this
    $t$-structure is compatible with the symmetric monoidal structure on $\GrMod_k$, there is an
    induced symmetric monoidal structure on $\GrMod_k^{\N\heart}$, which agrees with the Day
    convolution symmetric monoidal structure $\Gr(\Mod_k^\heart)^{\otimes_\D}$ on $\Gr(\Mod_k^\heart)$. This is the symmetric monoidal
    structure in which the commutative algebra objects are the algebraic geometer's graded
    commutative rings, so there is no sign when commuting odd weight classes.
\end{definition}

The other basic $t$-structure is the Beilinson, or negative, $t$-structure, although this terminology is usually
reserved for the lift of this $t$-structure to filtered $k$-modules.

\begin{definition}[The Beilinson $t$-structure]
    Let $(\GrMod_k)_{\geq 0}^\B\subseteq\GrMod_k$ be the full subcategory of graded $k$-modules
    $X^\ast$ where $X^n\in(\Mod_k)_{\geq -n}$ for all $n\in\bZ$. This defines the connective part of a
    $t$-structure, the \longdefidx{Beilinson}{$t$-structure!Beilinson} or
    \longdefidx{negative}{$t$-structure!Beilinson} $t$-structure, whose coconnective part
    $(\GrMod_k)_{\leq 0}^\B$ consists of the graded
    objects $Y^\ast$ where $Y^n\in(\Mod_k)_{\leq -n}$ for all $n\in\bZ$. The heart
    $\GrMod_k^{\B\heart}$ is again equivalent to $\Gr(\Mod_k^\heart)$. And, again, this $t$-structure is
    compatible with the symmetric monoidal structure on $\GrMod_k$, so there is an induced symmetric
    monoidal structure $\Gr(\Mod_k^\heart)^{\otimes_\K}$ on the heart. This time it is the Koszul symmetric monoidal structure, i.e.,
    the one where the braiding induces a sign for products of odd classes. Note that the underlying
    monoidal structures of $\Gr(\Mod_k^\heart)^{\otimes_\D}$ and $\Gr(\Mod_k^\heart)^{\otimes_\K}$ are equivalent.
\end{definition}

Besides these basic $t$-structures, shearing induces a plethora of other $t$-structures.

\begin{proposition}[Shearing up]\label{prop:shear}
    The autoequivalence $[2\ast]\colon\GrMod_k\rightarrow\GrMod_k$, which takes a graded object
    $X^\ast$
    to the graded object $X^\ast[2\ast]$ with weight $n$ piece $X^n[2n]$, admits a canonical symmetric monoidal structure.
\end{proposition}

\begin{proof}
    See~\cite[Prop.~3.3.4]{raksit}.
\end{proof}

\begin{definition}[$a$-fold shearing]
    It follows that the inverse $[-2\ast]$ of $[2\ast]$ admits a canonical symmetric monoidal
    structure and hence that the equivalence $[2a\ast]$, which assigns to $X^\ast$ the graded object
    $X^\ast[2a\ast]$, with weight $n$ piece $X^n[2an]$, admits a canonical symmetric monoidal
    structure for all $a\in\bZ$. If $a\geq 1$, this is called \defidx{shearing up} $a$-times or
    $a$-fold shearing up; if $a\leq -1$, this is called \defidx{shearing down} $a$-times or $a$-fold
    shearing down.
\end{definition}

\begin{definition}[Sheared neutral $t$-structures]
    The $a$-sheared neutral $t$-structure on $\GrMod_k$ has connective part $(\GrMod_k)^{\N[a]}_{\geq
    0}$ the full subcategory of graded $k$-modules $X^\ast$ where $X^n\in(\Mod_k)_{\geq 2an}$ for all
    $n\in\bZ$. In other words, it is the image of $(\GrMod_k)^\N_{\geq 0}$ under the $a$-fold shear
    functor $[2a\ast]$. As above, this $t$-structure is compatible with the symmetric monoidal
    structure; the heart $\GrMod_k^{\N[a]\heart}$ is equivalent to the symmetric monoidal abelian category
    $\Gr(\Mod_k^\heart)^{\otimes_\D}$. Shearing yields symmetric monoidal equivalences $(\GrMod_k)^{\N}_{\geq
    0}\we(\GrMod_k)^{\N[a]}_{\geq 0}$ and
    $\GrMod_k^{\N\heart}\we\GrMod_k^{\N[a]\heart}$, although these are distinct
    as subcategories of $\GrMod_k$. For example,
    $\GrMod_k^{\N[a]\heart}\subseteq\GrMod_k$ is the full subcategory of graded
    $k$-modules $X^\ast$ where $X^n\in(\Mod_k)_{[2an,2an]}$ for all $n$,
    or equivalently, $X^n[-2an]\in\Mod_k^\heart$ for all $n$.
\end{definition}

\begin{definition}[Sheared Beilinson $t$-structures]
    The $a$-sheared Beilinson $t$-structure on $\GrMod_k$ has connective part $(\GrMod_k)^{\B[a]}_{\geq
    0}$ the full subcategory of graded $k$-modules $X^\ast$ where $X^n\in(\Mod_k)_{\geq
    -n+2an}$ for all $n\in\bZ$. It is the image of $(\GrMod_k)^\B_{\geq 0}$ under the $a$-fold shear
    functor $[2a\ast]$. The special case of $a=1$ is called the \defidx{positive $t$-structure}
    in~\cite{raksit}. In general, $[2a\ast]\colon(\GrMod_k)^{\B}_{\geq 0}\we(\GrMod_k)^{\B[a]}_{\geq 0}$
    and $[2a\ast]\colon\GrMod_k^{\B\heart}\we\GrMod_k^{\B[a]\heart}$. These $t$-structures are also
    compatible with the symmetric monoidal structure on $\GrMod_k$ and the equivalences above are
    symmetric monoidal. In particular, all hearts are equivalent to $\Gr(\Mod_k^\heart)^{\otimes_\K}$.
    Again, these are distinct subcategories of $\GrMod_k$; the heart
    $\GrMod_k^{\B[a]\heart}$ is the full subcategory of objects $X^\ast$ where
    $X^n\in(\Mod_k)_{[-n+2an,-n+2an]}$ for all $n$.
\end{definition}

The following lemma is straightforward and left to the reader.

\begin{lemma}
    Each of the $t$-structures $(\GrMod_k)^{\N[a]}_{\geq 0}$ and
    $(\GrMod_k)^{\B[a]}_{\geq 0}$ is left and right complete, accessible, compatible with filtered
    colimits, and compatible with the symmetric monoidal structure on
    $\GrMod_k$. Moreover, the connective parts of these $t$-structures are compact
    projectively generated by static objects, of the form $P(n)[2an]$ (in the neutral case) or
    $P(n)[-n+2an]$ (in the Beilinson case) for $P$ a finitely presented project $k$-module and $n\in\bZ$.
\end{lemma}

It follows from the lemma that the methods of the previous sections apply to the $[a]$-sheared
neutral and Beilinson $t$-structures.

\begin{definition}[Graded commutative rings]
    Let $\Gr\CAlg_k=\CAlg(\GrMod_k)$ be the $\infty$-category of graded $\bE_\infty$-rings.
\end{definition}

\begin{definition}[Graded derived commutative rings]
    The \longdefidx{$[a]$-sheared graded derived commutative ring}{graded
    derived commutative ring!$[a]$-sheared} monad $\LSym^{\N[a]}$ is the derived monad
    $\L_1^{\N[a]}\bE_\infty$ associated to the
    $\bE_\infty$ monad with respect to the $[a]$-sheared neutral $t$-structure
    as in Definition~\ref{def:derivedapproximation}. The
    $\infty$-categories of algebras for these monads are written $\Gr\DAlg^{\N[a]}_k$.
\end{definition}

\begin{definition}[Derived graded-commutative rings]
    The \longdefidx{$[a]$-sheared derived graded-commutative ring}{graded
    derived commutative ring!$[a]$-sheared graded-commutative} monad $\LSym^{\B[a]}$ is the derived monad
    $\L_1^{\B[a]}\bE_\infty$ associated to the
    $\bE_\infty$ monad with respect to the $[a]$-sheared Beilinson $t$-structure. The
    $\infty$-categories of algebras for these monads are written
    $\Gr\DAlg^{\B[a]}_k$.
\end{definition}

\begin{remark}
    The monad $\LSym^{\N[a]}$ is an example of a derived commutative ring monad associated to a derived
    algebraic context in the sense of Raksit; in this case the $[a]$-sheared neutral $t$-structure
    defines such a context. However, $\LSym^{\B[a]}$ does not fit into this framework unless $2$ is
    invertible in $k$. Indeed, if $2$ is not a unit, symmetric powers of compact projective generators $P(n)$ for $n$ odd are not
    projective in the heart of the $a$-sheared Beilinson $t$-structure.
    Nevertheless, the $[a]$-sheared derived
    graded-commutative ring monad is constructed in the same way as $\LSym^{\N[a]}$ using Theorem~\ref{thm:apparatus}.
\end{remark}

\begin{remark}[Comparison to Raksit]
    The notation in~\cite{raksit} for what we call $\Gr\DAlg_k^{\N[0]}$ is $\mathrm{GrDAlg}(\Mod_k)$
    and is introduced in~\cite[Const.~4.3.4]{raksit}.
\end{remark}

\begin{definition}[Strict derived graded-commutative rings]
    The monad $\LSym^{\B[a]}$ is derived off of the free commutative ring monad on the heart
    $\GrMod_k^{\B[a]\heart}\we\Gr(\Mod_k^\heart)^{\otimes_\K}$. This is the free graded-commutative ring
    monad; there is also a free strict graded-commutative ring monad obtained by enforcing that all
    odd classes square to zero and this monad can in turn be derived off of the heart to obtain
    $\LSym^{\B[a]\s}$, whose $\infty$-category of algebras $\Gr\DAlg^{\B[a]s}_k$ are the
    \longdefidx{$[a]$-sheared strict derived graded-commutative rings}{graded
    derived commutative!$[a]$-sheared strict graded-commutative}.
\end{definition}

\begin{notation}[Special cases]
    A $[0]$-sheared graded derived commutative ring is called simply either a
    \defidx{graded derived commutative ring} or an \longdefidx{infinitesimal graded derived
    commutative ring}{graded derived commutative ring!infinitesimal}. A
    $[-1]$-sheared graded derived commutative ring is also
    called a \longdefidx{crystalline graded derived commutative ring}{graded
    derived commutative ring!crystalline}. A $[0]$-sheared
    derived graded-commutative ring is also called a \longdefidx{derived
    graded-commutative ring}{graded derived commtuative ring!
    graded-commutative} and a strict $[0]$-sheared
    derived graded-commutative ring is also a \longdefidx{derived strict
    graded-commutative ring}{graded derived commutative ring!strict
    graded-commutative}.
\end{notation}

\begin{remark}
    By shearing all of the $\infty$-categories $\Gr\DAlg^{\N[a]}_k$ are equivalent as are the
    graded-commutative variants. Specifically, in the context of Lemma~\ref{lem:oplax}, as $a$-fold shearing
    $$[2a*]\colon\GrMod_k\rightleftarrows\GrMod_k\colon[-2a*]$$ fits into an adjoint equivalence,
    the induced functor $\End(\GrMod_k)\rightarrow\End(\GrMod_k)$ given by $$T\mapsto
    T((-)[-2a*])[2a*]$$ is in fact opmonoidal instead of merely oplax monoidal, so it sends monads
    to monads. Applying this functor to $\LSym^{\N[0]}$ produces $\LSym^{\N[a]}$.
    This does not mean, however, that a given graded derived commutative
    ring admits the structure of an $[a]$-sheared derived graded commutative ring for all $a$. For
    example, the free $[0]$-sheared graded derived commutative ring over $\bZ$ on $\bZ(1)$ is the graded polynomial
    ring $\bZ[x]$ where $\bZ$ has weight $1$ and homological degree $0$. This does not admit the structure
    of a $[-1]$-sheared derived commutative ring because then $x$ would admit
    divided powers.
\end{remark}

\begin{remark}[Comparison of the monads]
    Let $M\in\Mod_k$ be a $k$-module and let $M(n)$ be the
    corresponding graded $k$-module supported in weight $n$. Figure~\ref{fig:monad_table} displays a
    table of the value of the various graded derived commutative ring monads on
    $M(n)$.
    \begin{figure}[h]
        \begin{align*}
            \LSym^{\N[a]}(M(n))&\we\bigoplus_{r\geq 0}\LSym^r(M[-2an])[2anr](nr),\\
            \LSym^{\B[a]}(M(n))&\we\begin{cases}
                \bigoplus_{r\geq 0}\LSym^r(M[n-2an])[-nr+2anr](nr)&\text{if $n$ is even,}\\
                \bigoplus_{r\geq 0}\LAntiSym^r(M[n-2an])[-nr+2anr](nr)&\text{if $n$ is odd,}
            \end{cases}\\
            \LSym^{\B[a]\s}(M(n))&\we\begin{cases}
                \bigoplus_{r\geq 0}\LSym^r(M[n-2an])[-nr+2anr](nr)&\text{if $n$ is even,}\\
                \bigoplus_{r\geq 0}\L\Lambda^r(M[n-2an])[-nr+2anr](nr)&\text{if $n$ is odd,}
            \end{cases}
        \end{align*}
        \label{fig:monad_table}
        \caption{The values of various graded $\LSym$ monads at $M(n)$. }
    \end{figure}
\end{remark}

\begin{remark}[Insertion and summation]\label{rem:insertion}
    The insertion functor $\Mod_k\rightarrow\GrMod_k$ given by $M\mapsto M(0)$ commutes
    with the $\LSym$ functor on $\Mod_k$ acting on the left and any of the
    $\LSym$-monads discussed above acting on the right, so that if $R$ is a derived commutative
    $k$-algebra, then $R(0)$ is an $a$-sheared derived graded commutative
    ring for all $a$ and also an $a$-sheared derived graded-commutative ring
    for all $a$ (and is even strict). On the other hand, the direct sum functor
    $\bigoplus\colon\GrMod_k\rightarrow\Mod_k$ commutes in the above sense only with
    $\LSym^{\N[0]}$. Thus, if $R^\ast$ is a $[0]$-sheared graded derived
    commutative ring, then $\bigoplus_{n\in\bZ} R^n$ is a derived commutative
    ring.

    One perspective on this fact is that graded $k$-modules are equivalent to
    comodules over the bicommutative Hopf algebra $k[t^{\pm 1}]$, which is naturally a derived commutative
    ring but even more strongly an $\bE_\infty$-coalgebra in derived
    commutative rings. One can define the $\infty$-category of $k[t^{\pm
    1}]$-comodules in $\DAlg_k$ to obtain a notation of graded derived
    commutative ring (see~\cite[Const.~4.2.32]{raksit}). It turns out that this
    recovers $\Gr\DAlg^{\N[0]}_k$. This
    is what singles out the $[0]$-sheared, or infinitesimal, graded derived commutative rings as
    special in the assortment of different flavors described here.
\end{remark}

\begin{variant}[Graded derived commutative $A$-algebras]\label{var:over_a} 
    Suppose now that $A$ is a derived commutative $k$-algebra. Then, we let
    $\DAlg_A=(\DAlg_k)_{A/}$. By Remark~\ref{rem:insertion}, we may view $A(0)$ as a graded
    derived commutative ring or a (strict) derived graded-commutative ring. We define
    $\Gr\DAlg_A^{\N[a]}$, $\Gr\DAlg_A^{\B[a]}$, and $\Gr\DAlg_A^{\B[a]\s}$ as the slice categories
    $$(\Gr\DAlg_k^{\N[a]})_{A(0)/},\quad(\Gr\DAlg_k^{\B[a]})_{A(0)/},\quad\text{and}\quad
    (\Gr\DAlg_k^{\B[a]\s})_{A(0)/}.$$
    Each of these is monadic over $\Gr\CAlg_A$, which can be defined as either the slice category
    under $A$ in $\Gr\CAlg_k$ or as the algebras for the $\bE_\infty$-monad in
    $\Gr\Mod_A=\Fun(\bZ^\delta,\Mod_A)$.
\end{variant}

Each of the monads above preserves the full subcategory $\Gr^+\Mod_k$ of objects concentrated
in nonnegative weights. And there are corresponding full subcategories $\Gr^+\DAlg^{\bN[0]}_k$,
etc., of graded derived commutative rings concentrated in nonnegative weights.

\begin{theorem}[Graded crystallization]\label{thm:gradedcrystallization}
    When restricted to $\Gr^+\Mod_k$, there are maps of monads
    $$\bE_\infty\rightarrow\LSym_k^{\N[0]}\rightarrow\LSym_k^{\B[0]}\rightarrow\LSym_k^{\B[0]\s}\rightarrow\LSym_k^{\N[-1]},$$
    which induce forgetful functors
    \begin{equation}\label{eq:fildalg}
            \Gr^+\CAlg_k\leftarrow\Gr^+\DAlg^{\N[0]}_k\leftarrow\Gr^+\DAlg^{\B[0]}_k\leftarrow\Gr^+\DAlg^{\B[0]s}_k\leftarrow\Gr^+\DAlg_k^{\N[-1]}.
    \end{equation}
    These functors are conservative and preserve limits and sifted colimits. The functors
    \begin{align*}
        \Gr^+\CAlg_k&\leftarrow\Gr^+\DAlg^{\N[0]}_k,\\
        \Gr^+\DAlg^{\N[0]}_k&\leftarrow\Gr^+\DAlg^{\B[0]\s}_k, \text{and}\\
        \Gr^+\DAlg^{\B[0]\s}_k&\leftarrow\Gr^+\DAlg_k^{\N[-1]}
    \end{align*}
    preserve all limits and colimits.
\end{theorem}

\begin{proof}
    By construction, there are canonical maps from $\bE_\infty$ to each of the
    other monads, and there is a map
    $\LSym_k^{\B[0]}\rightarrow\LSym_k^{\B[0]\s}$.
    The identity functor satisfies the hypotheses of Proposition~\ref{prop:functoriality}
    for the right $t$-exact identity functors
    $\Gr^+\Mod_k^{\N[0]}=\Gr^+\Mod_k^{\B[0]}=\Gr^+\Mod_k^{\N[-1]}$ by Lemmas~\ref{lem:toramplitude}
    and~\ref{lem:toramplitude_alternating}. Thus, there
    are maps of monads 
    $$\bE_\infty\rightarrow\LSym_k^{\N[0]}\rightarrow\LSym_k^{\B[0]}\rightarrow\LSym_k^{\N[-1]}.$$
    That $\LSym_k^{\B[0]}\rightarrow\LSym_k^{\N[-1]}$ factors through $\LSym_k^{\B[0]\s}$ uses the fact that the connective cover of
    $\LSym_k^{\N[-1]}$ with respect to the Beilinson $t$-structure takes torsion-free values on compact projective
    generators of the heart of the Beilinson $t$-structure on $\Gr^+\Mod_k$ by
    Lemma~\ref{lem:toramplitude}. Given the maps of monads,
    conservativity follows from the fact that these are $\infty$-categories of modules for monads.
    The claim about preservation of limits and sifted colimits follows from
    Lemma~\ref{lem:siftedcolimits}.

    The functors $\Gr^+\CAlg_k\leftarrow\Gr^+\DAlg^{\N[0]}_k$ and $\Gr^+\CAlg_k\leftarrow\Gr^+\DAlg_k^{\N[-1]}$
    arise from examples of derived algebraic contexts, so~\cite[Prop.~4.2.27]{raksit} applies to
    show that they preserve all colimits.
    The free strict graded-commutative ring monad is binomial on the heart
    of the Beilinson $t$-structure and the exterior powers of finitely
    presented projective objects in the heart are projective, so
    Proposition~\ref{prop:colimits} applies to show that the forgetful functor
    $\Gr^+\CAlg_k\leftarrow\Gr^+\DAlg^{\B[0]\s}_k$ preserves colimits. By conservativity and what
    we already observed at the beginning of this paragraph, 
    $\Gr^+\DAlg_k^{\N[0]}\leftarrow\Gr^+\DAlg^{\B[0]\s}_k$ preserves all colimits.
\end{proof}

\begin{corollary}
    If $R^\ast$ is a crystalline graded derived commutative ring concentrated
    in nonnegative weights, then
    $\bigoplus_{n\in\bN}R^n$ is a derived commutative ring.
\end{corollary}

\begin{remark}
    By shearing the map of monads $\LSym_k^{\N[0]}\rightarrow\LSym_k^{\N[-1]}$ from
    Theorem~\ref{thm:gradedcrystallization}, we can produce a sequence of maps of monads
    $$\cdots\rightarrow\LSym_k^{\N[1]}\rightarrow\LSym_k^{\N[0]}\rightarrow\LSym_k^{\N[-1]}\rightarrow\cdots\rightarrow\LSym_k^{\N[-a]}\rightarrow\cdots$$
    on $\Gr^+\Mod_k$
    and hence induced forgetful functors at the level of $\infty$-categories of
    modules for these monads. In particular, if $a\geq 0$ and if $R^\ast$ is a $[-a]$-sheared
    graded derived commutative ring concentrated in nonnegative weights, then
    $\bigoplus_{n\in\bN}R^n$ is a derived commutative ring.
\end{remark}

\begin{variant}\label{var:nonpositive}
    Let $\Gr^-\Mod_k$ be the $\infty$-category of graded objects in $\Mod_k$ concentrated in
    nonpositive weights and let $\Gr^-\DAlg^{\N[0]}_k$ denote the $\infty$-category of
    infinitesimal graded derived commutative rings concentrated in nonpositive weights,
    and similarly for the other notions of graded derived commutative ring.
    The proof of Theorem~\ref{thm:gradedcrystallization} works here too, but in reverse, to give
    maps of monads
    $$\bE_\infty\rightarrow\LSym_k^{\N[-1]}\rightarrow\LSym_k^{\B[0]}\rightarrow\LSym_k^{\B[0]\s}\rightarrow\LSym_k^{\N[0]}$$
    and hence limit and sifted colimit preserving forgetful functors
    $$\Gr^-\CAlg_k\leftarrow\Gr^-\DAlg^{\N[-1]}_k\leftarrow\Gr^-\DAlg^{\B[0]}_k\leftarrow\Gr^-\DAlg^{\B[0]\s}_k\leftarrow\Gr^-\DAlg^{\N[0]}_k.$$
    As above, one uses that the identity functor $\Gr^-\Mod_k\rightarrow\Gr^-\Mod_k$ is right
    $t$-exact for $$\Gr^-\Mod_k^{\N[-1]}\rightarrow\Gr^-\Mod_k^{\B[0]}\rightarrow\Gr^-\Mod_k^{\N[0]}.$$
    To see that $\LSym_k^{\B[0]}\rightarrow\LSym_k^{\N[0]}$ factors through the strict
    graded-commutative monad, it is enough to see that when restricted to compact projective
    objects $P[-n](n)$
    of $\Gr^-\Mod_k^{\B[0],\heart}$ the homology of the connective cover of
    $\LSym_k^{\N[0]}(P[-n](n))$ with resect to the Beilinson $t$-structure is torsion-free.
    But, this follows from Lemma~\ref{lem:toramplitude}. As in
    Theorem~\ref{thm:gradedcrystallization}, the forgetful functors not involving
    $\Gr^-\DAlg_k^{\B[0]}$ preserve all colimits.
\end{variant}

\begin{example}
    Suppose that $X$ is a projective variety over a field $k$ with a
    line bundle $L$. The graded object $L^{\otimes\ast}$ for $\ast\geq 0$ defines a graded derived
    commutative ring; indeed, it is $\LSym_{\Oscr_X}^{\N[0]}(L(1))$.
    It follows that $\R\Gamma(X,L^{\otimes\ast})$ is a graded derived commutative ring and that
    $\bigoplus_{n\in\bN}\R\Gamma(X,L^{\otimes
    n})$ admits the structure of a derived commutative ring, refining the previously understood
    $\bE_\infty$-structure. Applying $\pi_0$ recovers the usual graded
    commutative ring used to
    define the map to the projective space associated to $L$.
\end{example}

\begin{example}
    The cohomology of a Dirac scheme in the sense of Hesselholt--Pstr\k{a}gowski~\cite{hesselholt-pstragowski-i} is naturally a
    derived graded-commutative ring.
\end{example}

\begin{example}
    The functor $\Gr^+\CAlg_k\leftarrow\Gr^+\DAlg^{\B[0]}_k$ does not preserve
    general coproducts. To see this,
    consider $R=\LSym_k^{\B[0]}(k(1))$, which is static in the Beilinson $t$-structure with homotopy
    ring $k[x]/(2x^2)$ where $x$ has weight $1$. Then,
    $R\otimes_kR\rightarrow\LSym_k^{\B[0]}(k(1)\oplus k(1))$ is not an equivalence. Indeed, the
    latter is static in the Beilinson $t$-structure, but the former is not because of Tor
    interference. Specifically, the former is not static in the Beilinson $t$-structure in weight $4$.
\end{example}

\section{Filtered derived commutative rings}\label{sec:fdalg}

Let $k$ be a commutative ring.

\begin{definition}[Filtered derived $\infty$-categories]
    The filtered derived $\infty$-category of $k$ is $\FMod_k=\Fun(\bZ^\op,\Mod_k)$. The objects are decreasing filtrations
    $$\F^\star M\colon\cdots\rightarrow\F^{s+1}M\rightarrow\F^s
    M\rightarrow\F^{s-1}M\rightarrow\cdots$$ in $\Mod_k$.

    A filtration $\F^\star M$ is \longdefidx{complete}{filtration!complete} if $\lim_s\F^s M\we 0$.
    The object $\F^{-\infty}M=\colim_{s}\F^s M$ is called the underlying object and is also written
    as $|\F^\star M|$. The inclusion of the full subcategory
    $\widehat{\FMod}_k\subseteq\FMod_k$ of complete
    filtrations is closed under limits and admits a left adjoint called \defidx{completion}. The
    kernel of the completion functor $\FMod_k\rightarrow\widehat{\FMod}_k$ is the full subcategory of
    constant filtrations, which is equivalent to $\Mod_k$.
\end{definition}

\begin{example}
    Given $M\in\Mod_k$, the filtration $\ins^sM$ is the value at $M$ of the left Kan extension
    $\Mod_k\rightarrow\FMod_k$ along
    $\{s\}\rightarrow\bZ^\op$, which has value $0$ for $i>s$ and $M$ for $i\leq s$ with transition
    maps the identity on $M$ in weights at most $s$. It is
    complete with underlying object $M$.
\end{example}

\begin{construction}[Symmetric monoidal structure]
    The symmetric monoidal structure on $\bZ^\op$ given by addition gives rise to a symmetric
    monoidal structure on $\FMod_k$ by Day convolution.
    The unit object of $\FMod_k$ is $\ins^0k$ which we will simply write as $k$. The tensor product
    on $\FMod_k$ is just written as $\otimes_k$ and there is a pointwise description of the tensor
    product of two filtered $k$-modules:

    $$(\F^\star M\otimes_k\F^\star N)^s=\colim_{i+j\geq s}\F^iM\otimes_k\F^jN.$$
    The symmetric monoidal structure on $\FMod_k$ induces one on $\widehat{\FMod}_k$
    obtained by taking the tensor product first in $\F\Cscr$ and then completing; this will be
    written as $\F^\star M\tensorhat_{k}\F^\star N$.
\end{construction}

\begin{construction}[Associated graded]
    Given a filtration $\F^\star M\in\FMod_k$, there is associated graded object: $\gr^\ast
    M\in\GrMod_k$ whose weight $s$ piece is $\gr^sM\we\cofib(\F^{s+1}M\rightarrow\F^sM)$.
    Taking associated graded objects induces a functor $\gr^\ast\colon\FMod_k\rightarrow\GrMod_k$
    which is naturally symmetric
    monoidal and factors through the completion functor $\FMod_k\rightarrow\widehat{\FMod}_k$ to give a 
    conservative symmetric monoidal functor $\gr^\ast\colon\widehat{\FMod}_k\rightarrow\GrMod_k$.
\end{construction}

\begin{definition}[Nonnegative filtrations]
    Say that a filtration $\F^\star M\in\FMod_k$ is nonnegative if $\gr^iM\we 0$ for $i<0$.
    Equivalently, the maps $\F^0M\rightarrow\F^{-1}M\rightarrow\F^{-2}M\rightarrow\cdots$ are
    equivalences. Or, equivalently, $\F^\star M$ is left Kan extended from its restriction along
    $\bN^\op\rightarrow\bZ^\op$.
\end{definition}

There are three perspectives on the relationship between filtered and graded objects which will be
used below. These ideas are developed in~\cite{ariotta,lurie-rotation,raksit}, where we refer the
reader for proofs.

\begin{construction}[Rees construction]\label{const:rees}
    Restriction along the functor $\bZ^\delta\rightarrow\bZ^\op$ induces a forgetful functor
    $\GrMod_k\leftarrow\FMod_k$ which takes $\F^\star M$ and views it as a graded object $\F^\ast M$ with weight
    $s$ piece given by $\F^s M$ for all $s\in\bZ$. We will call this the underlying graded object
    construction.
    This forgetful functor is lax symmetric monoidal
    as it is the right adjoint of a symmetric monoidal functor $\ins$ given by left Kan extension along
    the symmetric monoidal functor
    $\bZ^\delta\rightarrow\bZ^\op$. Here, $\ins(M^\ast)$ is the filtered object with
    weight $s$ piece given by $\bigoplus_{i\geq s}M(i)$. This is a filtration on
    $\bigoplus_{i\in\bZ}M(i)$, which is typically not complete. The unit object $k=\ins^0k$ of $\FMod_k$ has
    underlying graded object the graded commutative algebra object $k[t]$ of $\Mod_k$ where $t$ has weight $-1$ and
    homological degree $0$, so it is actually in the heart $\GrMod_k^{\heart}$. The adjunction
    $\ins\colon\Gr\Cscr\rightleftarrows\F\Cscr\colon\F^\ast$
    is monadic and induces a symmetric monoidal equivalence $$\Mod_{k[t]}(\GrMod_k)\we\FMod_k.$$
    From this perspective, it is natural to write the image $\F^\ast M$ of $\F^\star M$ as the graded
    $k[t]$-module $\bigoplus_{i\in\bZ}\F^iM\cdot t^{-i}$, as in the classical Rees construction.
    The associated graded object is computed as $$\gr^\ast M\we\F^\ast M\otimes_{k[t]}k,$$
    where $k[t]\rightarrow k$ necessarily sends $t$ to zero.

    The neutral $t$-structure on $\GrMod_k$ induces a $t$-structure $(\FMod_k)_{\geq 0}^\N$ on
    $\FMod_k$, also called
    \longdefidx{neutral}{filtered object!neutral $t$-structure},
    with connective objects $\FMod_k$ those filtrations $\F^\star M$ such that
    $\F^sM\in(\Mod_k)_{\geq 0}$ for each $n$ and coconnective objects those
    such that $\F^sM\in(\Mod_k)_{\leq 0}$ for each $s$. This is also the pointwise $t$-structure
    from the description of $\FMod_k$ as $\Fun(\bZ^\op,\Mod_k)$ and it is compatible with the
    symmetric monoidal structure. The heart $\FMod_k^{\N\heart}$ is
    $\F(\Mod_k^\heart)$, the Grothendieck abelian category of decreasing
    filtrations in $\Mod_k^\heart$.\footnote{In our work, filtrations on objects in an abelian
    category are not required to have injective transition maps; when they do, they are called
    strict filtrations. This is why $\F(\Mod_k^\heart)$ is in fact an abelian category.}
\end{construction}

\begin{construction}[Coherent cochain complexes]\label{const:coherent}
    Ariotta proves in~\cite{ariotta}
    that $\widehat{\FMod}_k$ is equivalent to the pointed functor category
    $\Fun_*(\Xi^\op,\Mod_k)$, where $\Xi$ is the pointed category with objects
    $\bZ\sqcup\{*\}$, the set of integers plus a basepoint, where there is a unique
    non-trivial morphism $\partial\colon n\rightarrow n-1$ such that
    $\partial\circ\partial=*$ when the composition makes sense. (In fact, this holds more generally
    for any stable $\infty$-category with countable limits.)

    A \longdefidx{coherent chain complex}{coherent (co)chain complex}
    consists of a graded object $M^\ast$ together with maps $$\cdots \xrightarrow{\partial}
    M^{-1}[-1]\xrightarrow{\partial} M^0\xrightarrow{\partial}
    M^1[1]\xrightarrow{\partial} M^2[2]\xrightarrow{\partial}\cdots$$ as well as nullhomotopies
    $\partial^2\we 0$ and
    coherent higher nullhomotopies as well. Ariotta's equivalence sends a complete filtration
    $\F^\star M$ to
    the coherent chain complex
    $$\cdots\rightarrow\gr^{-1}M[-1]\rightarrow\gr^0M\rightarrow\gr^1M[1]\rightarrow\cdots.$$
    Taking homology yields the $\E^1$-page of the usual spectral sequence associated to the
    filtration $\F^\star M$. This connection between spectral sequences and Ariotta's theorem is
    explored further in~\cite{antieau-decalage}.

    The induced pointwise $t$-structure on $\widehat{\FMod}_k$ arising from its description as a
    (pointed) functor category is
    called the Beilinson $t$-structure. The $\infty$-category of connective
    objects $(\widehat{\FMod}_k)_{\geq 0}^\B$ consists of the complete
    filtrations $\F^\star M$ where $\gr^s M\in(\Mod_k){\geq -s}$ for all $s\in\bZ$.
    The coconnective objects are those complete filtrations $\F^\star M$ where $\F^sM\in(\Mod_k)_{\leq -s}$ for
    all $s\in\bZ$. The heart of the Beilinson $t$-structure is $\FMod_k^{\B\heart}\we\Ch^\bullet(\Mod_k^\heart)$, the
    abelian category of cochain complexes in $\Cscr$.
\end{construction}

\begin{variant}
    It is useful to also consider the Beilinson $t$-structure on all of $\FMod_k$.
    It is defined in the same way: $\F^\star M\in\FMod_k$ is connective if $\gr^sM\in\Cscr_{\geq -s}$ for all
    $s\in\bZ$, while $(\FMod_k)_{\leq 0}^\B=(\widehat{\FMod}_k)_{\leq 0}^\B$. The constant filtrations are infinitely connective in this $t$-structure.
\end{variant}

\begin{example}
    It is nice that $\Xi^\op$ is a $1$-category. This implies for example that
    $\widehat{\FMod}_k\we\D(\Ch^\bullet(k))$, the derived $\infty$-category of the Grothendieck
    abelian category of cochain complexes of static $k$-modules.
\end{example}

\begin{construction}[$\bD_-$-modules and $\bD_-^\vee$-comodules]\label{const:dminus}
    The functor $\gr^\ast\colon\widehat{\FMod}_k\rightarrow\GrMod_k$ is conservative and preserves all limits
    and colimits. Let $\bD_-^\vee$ be the graded bicommutative bialgebra $k\otimes_{k[t]}k$, which
    has underlying graded object $k\oplus k[1](-1)$ and let $\bD_-$ be its graded $k$-linear dual,
    which has underlying graded object $k\oplus k[-1](1)$.
    Using Koszul duality arguments, one proves that
    $\Mod_{\bD_-}(\GrMod_k)\we\cMod_{\bD_-^\vee}(\GrMod_k)\we\widehat{\FMod}_k$;
    see~\cite[Thm.~3.2.14]{raksit}.

    Note that $\bD_-$ is in the heart
    of the Beilinson $t$-structure, so that
    $\widehat{\FMod}_k^{\B\heart}\we\Mod_{\bD_-}(\GrMod_k^{\B\heart})\we\Mod_{\bD_-}(\Gr(\Mod_k^\heart))\we\Ch^\bullet(k)$,
    which gives another identification of the heart of the Beilinson $t$-structure on filtered
    complexes.
\end{construction}

\begin{construction}[Infinitesimal filtered derived commutative rings]
    Let $k$ be a commutative ring. The $\infty$-category $\FMod_k$ together
    with the neutral $t$-structure is a derived algebraic context:
    as $P$ ranges over the compact projective $k$-modules and
    $n$ ranges over $\bZ$, the objects
    $\ins^n(P)$ give the compact projective objects of
    $\F(\Mod_k^\heart)$. Let $\LSym_k^\inf$ denote the associated derived commutative ring monad on
    $\FMod_k$. One checks that $$\LSym_k^\inf(\ins^n(M))\we\bigoplus_{s\geq
    0}\ins^{ns}\LSym^s_k(M).$$
    Let $\FDAlg_k^\inf$ be the $\infty$-category of $\LSym_k^\inf$-algebras.
    The
    \longdefidx{associated graded functor}{associated graded}
    $$\gr^\ast\colon\FMod_k\rightarrow\GrMod_k$$ is a morphism of derived
    algebraic contexts where we view $\GrMod_k$ as a derived algebraic context via the
    $[0]$-sheared neutral $t$-structure. Thus, there is an induced colimit preserving
    associated graded functor
    $$\gr^\ast\colon\FDAlg_k^\inf\rightarrow\Gr\DAlg_k^{\N[0]}$$
    which fits into a commutative diagram
    $$\xymatrix{
        \FDAlg_k^\inf\ar[r]^{\gr^\ast}\ar[d]&\Gr\DAlg_k^{\N[0]}\ar[d]\\
        \FMod_k\ar[r]^{\gr^\ast}&\GrMod_k,
    }$$
    where the vertical arrows are the forgetful functors.
    The $\LSym$ functor on $\FMod_k$ preserves constant filtrations and thus
    induces an $\LSym$ functor on $\widehat{\FMod}_k$; see~\cite[Prop.~4.1.10]{raksit} for the
    argument that one can transport the monad over the localization in this
    situation. We let $\widehat{\FDAlg}_k^\inf$ be the $\infty$-category of complete infinitesimal
    filtered derived commutative rings. Under the equivalence
    $\gr^\ast\colon\widehat{\FMod}_k\we\Mod_{\bD_-}(\GrMod_k)\we\cMod_{\bD^\vee_-}(\GrMod_k)$,
    we have $\widehat{\FDAlg}^\inf_k\we\cMod_{\bD^\vee_-}(\Gr\DAlg^\inf_k)$. 
\end{construction}

\begin{remark}[Comparison to Raksit]
    The notation in~\cite{raksit} for what we call $\FDAlg_k^\inf$ is $\mathrm{FilDAlg}(\Mod_k)$
    and is introduced in~\cite[Const.~4.3.4]{raksit}.
\end{remark}

\begin{remark}
    We view the notion of infinitesimal filtered derived commutative ring as the basic, possibly
    naive, notion of filtered derived commutative ring. The free objects on the basic filtered
    projective $k$-modules $\ins^ik$ are the static filtered commutative rings $k[x]$ where $x$ has
    weight $i$.
\end{remark}

Recall that $\bD_-^\vee$ is the Bar complex of the graded Rees algebra and that
$\widehat{\FMod}_k\we\cMod_{\bD_-^\vee}(\GrMod_k)$. The $k$-linear dual of $\bD_-^\vee$ is denoted
by $\bD_-$ and has underlying graded $k$-module $k\oplus k[-1](1)$. Note that there is a unique
$\bE_1$-$k$-algebra structure on $k\oplus k[-1](1)$, for instance by~\cite[Thm.~1.3]{dwyer-greenlees-iyengar-exterior}.

\begin{lemma}\label{lem:dminus}
    There are unique infinitesimal and crystalline derived graded bicommutative bialgebra
    structures on $\bD^\vee_-$. 
\end{lemma}

\begin{proof}
    As an infinitesimal graded derived commutative $k$-algebra, $\bD^\vee_-$ is equivalent to
    $\LSym_{\Gr\Mod_k}^{[0]}(k[1](-1))$, i.e., it is free on a generator $\d$ in weight $-1$ and
    homological degree $1$. As with any free algebra, it follows that $\bD^\vee_-$ admits an
    $\bE_\infty$-coalgebra structure in $\Gr\CAlg^\inf_k$. In fact, this structure is unique thanks
    to the grading. Indeed, $\bD^\vee_-$ is in the heart of the Beilinson $t$-structure and, when
    restricted to nonpositive graded objects, $\LSym^\inf$ preserves Beilinson-connective objects.
    Thus, the infinitesimal $\bE_\infty$-coalgebra structures on $\bD^\vee_-$ are determined by
    such structures on $\bD^\vee_-$ when viewed as a commutative algebra in the heart of the
    Beilinson $t$-structure. Thus, we see the claimed uniqueness. A weight argument shows that the
    generator $\d$ must be mapped to $\d\otimes 1+1\otimes\d$.

    Using Variant~\ref{var:nonpositive}, there is an induced crystalline graded derived
    bicommutative bialgebra structure on $\bD_-^\vee$. 
    For uniqueness of the crystalline derived graded bicommutative bialgebra, we appeal
    to~\cite[Prop.~5.1.7]{raksit}. This result says that the sheared up version of $\bD_-^\vee$, which is
    denoted by $\bD_+^\vee$ admits a unique structure of an infinitesimal derived graded
    bicommutative bialgebra. But, shearing implies that the crystalline derived graded
    bicommutative algebra structures on $\bD_-^\vee$ correspond to infinitesimal filtered derived
    graded bicommutative bialgebra structures on $\bD_+^\vee$. This completes the proof.
\end{proof}

\begin{construction}[Crystalline filtered derived commutative rings]
    Using that $\bD^\vee_-$ is a crystalline graded derived bicommutative bialgebra, we let
    $$\widehat{\FDAlg}^\crys_k=\cMod_{\bD^\vee_-}(\Gr\DAlg^\crys_k).$$
    If $\widehat{\LSym}^\crys_k$
    denotes the crystalline filtered derived commutative ring monad on $\widehat{\FMod}_k$, then
    there is an induced monad on $\FMod_k$ by Lemma~\ref{lem:oplax}. This monad does not preserve
    sifted colimits, but as $\FMod_k$ is compactly generated, it may be approximated by a filtered
    colimit preserving functor by left Kan extending from the value on compact objects using a
    variant of the proof of Proposition~\ref{prop:siftedapproximation}. As the restriction of
    $\widehat{\LSym}^\crys_k$ to $\FMod_k^\omega$ preserves finite geometric realizations, the left
    Kan extension $\LSym^\crys_k$ does too and hence it preserves all sifted colimits.
    We let $\FDAlg_k^\crys$ denote the $\infty$-category of $\LSym^\crys_k$-modules in $\FMod_k$.
    There is a completion adjunction
    $$\FDAlg_k^\crys\rightleftarrows\widehat{\FDAlg}_k^\crys$$
    with fully faithful right adjoint.
\end{construction}

\begin{remark}[Comparison to Raksit and Moulinos--Robalo--To\"en]
    Raksit and Moulinos--Robalo--To\"en work ``sheared up''. Our $\infty$-category $\widehat{\FDAlg}^\crys_k$ is equivalent to
    Raksit's $\infty$-category of $h_+$-differential graded derived commutative $k$-algebras
    from~\cite[Def.~5.1.10]{raksit}. These are also related to the ``simplicial commutative
    mixed graded algebras'' of~\cite[Rem.~5.3.2]{mrt}.
\end{remark}

\begin{variant}[Filtered derived commutative $A$-algebras]\label{var:filtered_over_a} 
    Suppose that $A$ is a derived commutative $k$-algebra. Then, $\ins^0A$ admits both an
    infinitesimal and a crystalline filtered derived commutative ring structure and we let
    $$\FDAlg_A^\inf\quad\text{and}\quad\widehat{\FDAlg}_A^\inf\quad\text{and}\quad\FDAlg_A^\crys\quad\text{and}\quad\quad\widehat{\FDAlg}_A^\crys$$
    denote the slice categories under $\ins^0A$.
\end{variant}

\section{Hochschild homology and derived de Rham cohomology}\label{sec:raksit}

In~\cite{raksit}, Raksit proves universal properties of HKR-filtered Hochschild homology and
Hodge-filtered Hodge-complete derived de Rham cohomology. We review those ideas in our language.

Let $S^1$ denote the anima underlying the circle, let $\bT=\bZ[S^1]$ be its coalgebra of
$k$-linear chains, and let $\bT^\vee=\bZ^{S^1}$ be the $\bZ$-linear dual of $\bT$, which is naturally a derived
commutative $\bZ$-algebra.

\begin{lemma}\label{lem:cochain_existence}
    The derived commutative ring $\bT^\vee\in\DAlg_\bZ$ admits the structure of a
    cocommutative coalgebra in $\DAlg_\bZ$.
\end{lemma}

\begin{proof}
    Let $\Sscr_{\mathrm{ft}}\subseteq\Sscr$ be the full subcategory of anima of finite type.
    By K\"unneth, taking cochains yields a symmetric monoidal functor
    $\Sscr^\op_{\mathrm{ft}}\rightarrow\DAlg_\bZ$ where we equip $\Sscr^\op_{\mathrm{ft}}$ with the coCartesian symmetric
    monoidal structure. As $S^1$ admits the structure of an
    $\bE_\infty$-algebra in $\Sscr_{\mathrm{ft}}$, taking cochains produces an
    $\bE_\infty$-coalgebra on $\bT^\vee$ in $\DAlg_{\bZ}$.
\end{proof}

\begin{remark}
    In fact, the anima of $\bE_\infty$-algebra structures on $S^1$ in $\Sscr$ is contractible.
    By using binomial derived $\lambda$-rings $\DAlg_\bZ^{\mathrm{bin}}$ as
    in~\cite{antieau-spherical,horel,kubrak-shuklin-zakharov}, we see that $\bT^\vee$
    admits a unique structure of a cocommutative coalgebra in $\DAlg_\bZ^{\mathrm{bin}}$.
\end{remark}

The crucial structure needed to define the filtered circle is the infinitesimal filtered derived 
bicommutative bialgebra structure discovered by~\cite{raksit}.

\begin{proposition}[Functions on the filtered circle]\label{prop:existence}
    Let $\bT^\vee_\fil$ denote the Whitehead tower $\tau_{\geq\star}\bT^\vee$ of $\bT^\vee$. Then,
    $\bT^\vee_\fil$ admits the structure of a cocommutative coalgebra in infinitesimal filtered
    derived commutative rings lifting the object $\bT^\vee\in\cCAlg(\DAlg_\bZ)$ under the colimit
    functor $\colim\colon\cCAlg(\FDAlg_\bZ^\inf)\rightarrow\cCAlg(\DAlg_\bZ)$.
\end{proposition}

\begin{proof}[Sketch of proof]
    See~\cite[Thm.~6.1.5]{raksit} for the original proof. Here is another argument, inspired by the
    filtered Cartier duality perspective of~\cite{moulinos-cartier} on~\cite{mrt}, which we
    discovered with Sanath Devalapurkar. Consider the group ring $\bZ[\bZ]=\bZ[u^{\pm 1}]$ and equip
    it with the $(u-1)$-adic filtration $(u-1)^\star\bZ[u^{\pm 1}]$.
    The comultiplication $\bZ[u^{\pm 1}]\rightarrow\bZ[u_1^{\pm 1}]\otimes\bZ[u_2^{\pm 1}]$ determined
    by $u\mapsto u_1\otimes u_2$ sends $u-1$ to $(u_1-1)(u_2-1)+(u_1-1)+(u_2-1)$ and hence is
    compatible with the $(u-1)$-adic filtration on the left-hand side and the $(u_0-1,u_1-1)$-adic
    filtration on the right-hand side. Thus, $(u-1)^\star\bZ[u^{\pm 1}]$ is naturally a
    bicommutative bialgebra in static filtered $k$-modules. As it is moreover flat in each weight,
    it defines an object of $\cCAlg(\F\DAlg_\bZ^\inf)$.
    The bar construction with respect to the augmentation
    $\bZ[u^{\pm 1}]\rightarrow\bZ$ sending $u$ to $1$ therefore inherits the same structure
    and it is immediately seen to be $\bT_\fil=\tau_{\geq\star}\bZ[S^1]$. Thus, the filtered
    $\bZ$-linear dual $\bT_\fil^\vee$ also admits such a structure. Specifically, the $\bZ$-linear dual
    of $\bZ\otimes_{(u-1)^\star\bZ[u^{\pm 1}]}\bZ$ is computed as
    $$\Tot\left(\Hom_{\FMod_{\bZ}}(\bZ\otimes((u-1)^\star\bZ[u^{\pm
    1}])^{\otimes\bullet+1}\otimes\bZ,\bZ)\right).$$
    which is thus a sifted limit of a cosimplicial object $R^\bullet$ in $\cCAlg(\DAlg_\bZ)$. It suffices to see that the
    pointwise comultiplication induces a comultiplication on the limit. For this, it is enough to
    show that the natural map $\Tot(R^\bullet)^{\otimes n}\rightarrow\Tot((R^\bullet)^{\otimes n})$
    is an equivalence for all $n\geq 1$. Ignoring filtrations, the left-hand side computes the $n$-fold tensor product of
    $\bT^\vee$ while the right-hand side computes the dual of the bar construction of
    $\bZ[\bZ^n]$. These agree; checking that the filtrations match up is left to the reader.
\end{proof}

\begin{remark}
    Our proof sketch of Proposition~\ref{prop:existence} shows that $\bT$ also admits the
    structure of a derived bicommutative bialgebra and that there is a filtered refinement
    $\bT_\fil\in\cCAlg(\FDAlg_k^\inf)$ lifting $\bT\in\cCAlg(\DAlg_k)$. We remark that the
    cocommutative coalgebra structure on $\bT$ is not unique over $\bZ$. Indeed, if we replace
    $(u-1)^\star\bZ[u^{\pm 1}]$ with $x^\star\bZ[x]$, which corresponds to the group scheme
    $\bG_a$, then its bar complex is equivalent to $\bT_\fil$ as an infinitesimal filtered derived commutative ring, but
    not with the cocommutative comultiplication. They do become equivalent over $\bQ$ via a
    logarithm.
\end{remark}

\begin{definition}[The filtered circle]
    For us, the filtered circle is $\bT_\fil^\vee$, the infinitesimal filtered derived
    bicommutative bialgebra of Proposition~\ref{prop:existence}. Since this is a dualizable
    filtered $k$-module, there is an equivalence
    $\cMod_{\bT_\fil^\vee}(\FMod_\bZ)\we\Mod_{\bT_\fil}(\FMod_\bZ)$ and it is convenient to pass
    back and forth between these two perspectives. However, for the purposes of Hochschild
    homology, it is really the comodule perspective which is more central.
    If $k$ is a derived commutative ring,
    then by base change $k\otimes_\bZ\bT_\fil^\vee$ is an infinitesimal filtered derived
    bicommutative bialgebra over $k$, i.e., an object of $\cCAlg(\FDAlg_k^\inf)$.
    We will typically denote it by $\bT_\fil^\vee$ leaving the base implicit.
\end{definition}

\begin{corollary}\label{cor:gr_circle}
    The graded circle $\bT_\gr=\gr^*\bT_\fil$ and its dual $\bT_\gr^\vee=\gr^*\bT_\fil^\vee$ admit
    infinitesimal graded derived bicommutative bialgebra structures. Moreover,
    $\gr^*\bT_\fil^\vee[-2*]\we\bD^\vee_-$ as infinitesimal and crystalline graded derived bicommutative bialgebras.
\end{corollary}

\begin{proof}
    Only the second statement needs to be proved, but it follows from the uniqueness part of
    Lemma~\ref{lem:dminus}.
\end{proof}

\begin{definition}[Hochschild homology]
    Let $k$ be a derived commutative ring.
    The forgetful functor
    $$\DAlg_k\leftarrow\cMod_{\bT^\vee}(\DAlg_k)$$ admits a left adjoint
    $$\HH(-/k)\colon\DAlg_k\rightarrow\cMod_{\bT^\vee}(\DAlg_k)$$ called
    \defidx{Hochschild homology} relative to $k$.
\end{definition}

The next result is due to McClure--Schw\"anzl--Vogt~\cite{mcclure-schwanzl-vogt} in the $\bE_\infty$-case.
It says in particular that our definition of $\HH(R/k)$ agrees with the classical definition.

\begin{lemma}[Tensor product formula]
    Let $k$ be a derived commutative ring.
    There is an equivalence $\cMod_{\bT^\vee}(\DAlg_k)\we\Fun(S^1,\DAlg_k)$ and hence $\HH(R/k)\we
    S^1\otimes_k R$, the copowering of $R$ by $S^1$ in $\DAlg_k$. As a derived commutative ring,
    $\HH(R/k)\we R\otimes_{R\otimes_kR}R$.
\end{lemma}

\begin{proof}[Sketch of proof]
    The first claim follows from a comonadicity argument, and this also proves that $\HH(R/k)\we
    S^1\otimes_k R$. We can use the compatibility of copowering with colimits to present $S^1$ as
    the pushout of the span $*\leftarrow S^0\rightarrow *$ so that $S^1\otimes_kR$ is the pushout
    of $*\otimes_k R\leftarrow S^0\otimes_kR\rightarrow *\otimes_kR$, which gives the tensor
    product formula.
\end{proof}

\begin{remark}
    While $k$-linear duality induces an equivalence
    $\cMod_{\bT^\vee}(\Mod_k)\we\Mod_{\bT}(\Mod_k)$, the symmetric monoidal structure on
    $\Mod_{\bT}(\Mod_k)\we\Fun(S^1,\Mod_k)$ is the pointwise one corresponding to parametrized
    $k$-modules over $S^1$. Thus, $\CAlg(\Fun(S^1,\Mod_k))\we\Fun(S^1,\CAlg_k)\we\cMod_{\bT^\vee}(\CAlg_k)$. This is
    not equivalent to $\Mod_\bT(\CAlg_k)$.
\end{remark}

\begin{definition}[Filtered Hochschild homology]
    Let $k$ be a derived commutative ring. The functor
    $$\DAlg_k\leftarrow\cMod_{\bT_\fil^\vee}(\FDAlg_k^\inf)$$
    which forgets the $\bT_\fil^\vee$-coaction and takes $\F^0$ 
    admits a left adjoint
    $$\HH_\fil(-/k)\colon\DAlg_k\rightarrow\cMod_{\bT^\vee_\fil}(\FDAlg_k^\inf)$$
    called \longdefidx{filtered Hochschild homology}{Hochschild
    homology!filtered}.
\end{definition}

\begin{remark}
    By functoriality, the colimit functor
    $\cMod_{\bT_\fil^\vee}(\FDAlg_k^\inf)\rightarrow\cMod_{\bT}(\DAlg_k)$
    takes $\HH_\fil(R/k)$ to $\HH(R/k)$. In particular, $\HH_\fil(R/k)$ is some
    filtered version of Hochschild homology.
\end{remark}

\begin{remark}[Comparison to Raksit]
    Raksit uses the shorthand $\mathrm{Fil}_{S^1}\DAlg_k$ for what we have denoted by
    $\cMod_{\bT_\fil^\vee}(\FDAlg_k)$.
\end{remark}

The main move in the proof that $\HH_\fil(R/k)$ {\em is} HKR-filtered
Hochschild homology is to take the \longdefidx{sheared down associated
graded}{associated graded!sheared down}. Let $M$ be a
$\bT_\fil^\vee$-comodule in
$\FMod_k$. It has an associated graded $\gr^\ast M$ which is a
$\bT_\gr^\vee$-comodule. Shearing down,
$\gr^\star M[-2\ast]$ has a $\bT_\gr[-2\ast]\we\bD_-^\vee$-coaction. In other
words, it is a coherent chain complex. By applying this construction
to $\HH_\fil(R/k)$ and identifying the universal property of the result, Raksit
verifies that the filtration on Hochschild homology is indeed the HKR
filtration.

\begin{theorem}[Raksit]\label{thm:raksit}
    Let $k$ be a derived commutative ring.
    There is an adjunction
    $$\F^\star_\H\widehat{\dR}_{-/k}\colon\DAlg_k\rightleftarrows\cMod_{\bD_-^\vee}(\GrDAlg_k^\crys)\we\widehat{\FDAlg}_k^\crys\colon\gr^0,$$
    where $\F^\star_\H\widehat{\dR}_{R/k}$ is the Hodge-complete
    derived de Rham complex
    $$(R\rightarrow\L_{R/k}\rightarrow\Lambda^2\L_{R/k}\rightarrow\cdots).$$
\end{theorem}

\begin{proof}
    To prove the theorem, we restrict to complete and graded objects where $\gr^i\we 0$ for $i<0$ and
    decorate the corresponding $\infty$-categories with a $+$ as in
    Theorem~\ref{thm:gradedcrystallization}. Let $F$ denote the left adjoint to
    $\DAlg_k\leftarrow\F^+\DAlg^\crys_k\colon\gr^0$ and consider the commutative diagram
    $$\xymatrix{
        \Mod_k\ar[r]^<<<<<{(-)\otimes\bD_-}\ar[d]_\LSym&\cMod_{\bD_-^\vee}(\Gr^+\Mod_k)\ar[r]^<<<<<{\gr^\ast}\ar[d]&\Gr^+\Mod_k\ar[d]\\
        \DAlg_k\ar[r]_<<<{F}&\cMod_{\bD_-^\vee}(\Gr^+\DAlg_k^\crys)\ar[r]_<<<<<{\gr^\ast}&\Gr^+\DAlg_k^\crys.
    }$$
    We extend this to the right by restricting to the $\infty$-category $\Gr^{\{0,1\}}\Mod_k$ of objects of weights $0$ and $1$. This is a
    symmetric monoidal localization and the resulting theories of infinitesimal and crystalline
    derived commutative algebras in $\Gr^{\{0,1\}}\Mod_k$ agree. The result is the following
    commutative diagram:
    $$\xymatrix{
        \Mod_k\ar[r]^<<<<<{(-)\otimes\bD_-}\ar[d]_\LSym&\cMod_{\bD_-^\vee}(\Gr^+\Mod_k)\ar[r]^<<<<<{\gr^\ast}\ar[d]&\Gr^+\Mod_k\ar[d]\ar[r]^{\ev_{0,1}}&\Gr^{\{0,1\}}\Mod_k\ar[d]\\
        \DAlg_k\ar[r]_<<<{F}&\cMod_{\bD_-^\vee}(\Gr^+\DAlg_k^\crys)\ar[r]_<<<<<{\gr^\ast}&\Gr^+\DAlg_k^\crys\ar[r]_{\ev_{0,1}}&\Gr^{\{0,1\}}\DAlg_k.
    }$$
    The derived commutative algebras in $\Gr^{\{0,1\}}\Mod_k$ are precisely the
    pairs consisting of a derived commutative $k$-algebra $A$ (in weight $0$) and an $A$-module
    $M$ (in weight $1$). The total right adjoint $\Gr^{\{0,1\}}\DAlg_k\rightarrow\DAlg_k$
    sends $(A,M)\we A(0)\oplus M(1)$ to the square-zero extension $A\oplus M[1]$.
    It follows by definition of the cotangent complex that the total left
    adjoint $$\DAlg_k\rightarrow\Gr^{\{0,1\}}\DAlg_k$$
    sends $R$ to $R(0)\oplus\L_{R/k}[-1](1)$. (See for
    example~\cite[Sec.~4.4]{raksit}, but note the shift because Raksit is
    working `sheared up'.)

    This shows that there are natural identifications $\gr^0 F(R)\we R$ and
    $\gr^1 F(R)\we \L_{R/k}[-1]$. Since $\gr^\ast F(R)$ is naturally a
    $[-1]$-sheared derived commutative $R\we \gr^0 F(R)$-algebra, it follows that
    there is a canonical map
    $$\LSym^{\N[-1]}_{R}(\L_{R/k}[-1](1))\rightarrow\gr^\ast F(R).$$
    Now, it is easy to see, and left to the reader to check, that this map is
    an equivalence when $R\we\LSym_k(M)$ is itself free on an object
    $M\in\Mod_k$. It follows that it is an equivalence for all $R\in\DAlg_k$.
    This shows that in general $$\gr^\ast
    F(R)\we\Lambda^\ast\L_{R/k}[-\ast](\ast).$$ To complete the proof, it is
    enough to see that the coherent chain complex
    $$R\rightarrow\L_{R/k}\rightarrow\Lambda^2\L_{R/k}\rightarrow\cdots$$ is
    the Hodge-complete derived de Rham complex. For this, one reduces to the
    case where $k$ is static and $R$ is smooth over $k$. In this case, $F(R)$
    is static in the Beilinson $t$-structure on complete filtered complexes
    and it has a universal property with respect to strict
    commutative dgas; i.e., it is the classical de Rham complex
    $\Omega^\bullet_{R/k}$. In particular, $R\rightarrow\L_{R/k}$ is the
    universal derivation, and the theorem follows.
\end{proof}

\begin{corollary}
    There is a natural equivalence
    $$\gr^\ast\HH_\fil(R/k)[-2\star]\we\F^\star_\H\dRhat_{R/k}\we\left(R\rightarrow\L_{R/k}\rightarrow\Lambda^2\L_{R/k}\rightarrow\cdots\right)$$
    in $\cMod_{\bD_-^\vee}(\GrDAlg^\crys_k)$.
\end{corollary}

\begin{proof}
    The commutative diagram
    \small
    $$\xymatrix{
        \Mod_k\ar[r]^>>>{\ins^0(-)\otimes\bT_\fil}\ar[d]_{\LSym}&\cMod_{\bT_\fil^\vee}(\FMod_k)\ar[r]^{\gr^\ast}\ar[d]_{\LSym^\inf}&\cMod_{\bT_\gr^\vee}(\GrMod_k)\ar[r]_\we^{\text{shear
        down}}\ar[d]_{\LSym^\inf}&\cMod_{\bD_-^\vee}(\GrMod_k)\ar[r]^<<<<<<{|-|}_<<<<<<\we\ar[d]_{\LSym^\crys}&\widehat{\FMod}_k\ar[d]_{\LSym^\crys}\\
        \DAlg_k\ar[r]_>>>{\HH_\fil(-/k)}&\cMod_{\bT_\fil^\vee}(\FDAlg^\inf_k)\ar[r]_{\gr^\ast}&\cMod_{\bT_\gr^\vee}(\GrDAlg^\inf_k)\ar[r]^\we_{\text{shear
        down}}&\cMod_{\bD_-^\vee}(\GrDAlg^\crys_k)\ar[r]^<<<<<\we_<<<<<{|-|}&\widehat{\FDAlg}_k^\crys\\
    }$$\normalsize
    captures the relevant functoriality.
    The total compositions of the top $\Mod_k\rightarrow\widehat{\FMod}_k$ and of the bottom
    $\DAlg_k\rightarrow\widehat{\DAlg}^\crys_k$ have
    right adjoints given by $\gr^0$. Thus, by composition of left adjoints, the claim follows from
    Theorem~\ref{thm:raksit}.
\end{proof}

Let $k$ be a commutative ring. If $R$ is a smooth commutative $k$-algebra, let
$\F^\star_\HKR\HH(R/k)=\tau_{\geq\star}\HH(R/k)$, which is a filtered $\bE_\infty$-$k$-algebra.
For a general animated commutative $k$-algebra
$R$, let $\F^\star_\HKR\HH(R/k)$ denote the filtration obtained by left Kan extension from the smooth
case. We call this the derived HKR filtration for the moment.

\begin{corollary}\label{cor:hkr_comparison}
    If $k$ is a commutative ring and $R$ is an animated commutative $k$-algebra, then
    $\HH_\fil(R/k)\we\F^\star_\HKR\HH(R/k)$ as filtered $\bE_\infty$-$k$-algebras.
\end{corollary}

\begin{proof}
    Given filtered Hochschild homology $\HH_\fil(R/k)$, the sheared down associated
    graded $\gr^\star\HH_\fil(R/k)[-2\star]$ is a coherent chain complex
    $$\gr^0\HH_\fil(R/k)\rightarrow\gr^1\HH_\fil(R/k)[1]\rightarrow\gr^2\HH_\fil(R/k)[2]\rightarrow\cdots.$$
    Moreover, $\HH_\fil(R/k)$ will be a filtered derived commutative ring with
    $\bT_\fil$-action so that
    $\gr^\star\HH_\fil(R/k)$ will be a $[0]$-sheared graded derived commutative
    ring with a $\bT_\gr$-action and the shear down $\gr^\star\HH_\fil(R/k)[-2\star]$ will be a
    $[-1]$-sheared graded derived commutative ring with a $\bD_-^\vee$-coaction using the
    $[-1]$-sheared derived bicommutative bialgebra structure of
    Lemma~\ref{lem:dminus}. The theorem says that the sheared down associated
    graded of $\HH_\fil(R/k)$ is $\F^\star_\H\dRhat_{R/k}$.
    (The shear functors also play an important role in the
    paper~\cite{bzn} of Ben-Zvi and Nadler on loop spaces and
    connections in characteristic zero.)
    In particular, this
    implies that $$\gr^n\HH_\fil(R/k)\we\Lambda^n\L_{R/k}[n].$$
    When $k$ is static and $R$ is smooth over $k$, it follows that the filtration
    is precisely the Postnikov filtration on $\HH(R/k)$; since $\HH_\fil(R/k)$
    commutes with sifted colimits as a functor to filtered complexes, it follows
    that $\HH_\fil(R/k)$ is the derived HKR filtration.
\end{proof}

\begin{corollary}\label{cor:completeness}
    If $k\rightarrow R$ is a map of animated commutative rings, then $\F^s\HH_\fil(R/k)$ is
    $s$-connective for all $s\in\bZ$ and hence the filtration $\HH_\fil(R/k)$ is complete.
\end{corollary}

\begin{proof}
    First, we assume that $k$ is static. It is straightforward then to check that
    $\F^s\HH_\fil(R/k)$ is $s$-connective for all polynomial
    $k$-algebras. Since $\HH_\fil(R/k)$ commutes with sifted
    colimits, it follows that $\F^s\HH_\fil(R/k)$ is $s$-connective in general for animated
    commutative $k$-algebras. This implies completeness of the filtration when $k$ is static.
    Now, in general, if $k$ and $R$ are both animated, then we can compute $\HH_\fil(R/k)$ as
    $\HH_\fil(R/\bZ)\otimes_{\HH_\fil(k/\bZ)}k$. Using that $\F^s\HH_\fil(R/\bZ)$ and
    $\F^s\HH_\fil(k/\bZ)$ are both $s$-connective, together with the bar construction of the
    relative tensor product and the formula for the Day convolution symmetric monoidal structure,
    it follows that $\F^s\HH_\fil(R/k)$ is $s$-connective as well. This implies completeness.
\end{proof}

Using Raksit's work~\cite{raksit}, or the approach of Moulinos--Robalo--To\"en~\cite{mrt}, one can prove the following theorem.
It was first proven for animated commutative rings in~\cite{antieau-derham}, following
work in the $p$-adic case in~\cite{bms2}. In characteristic $0$ it is due
to~\cite{toen-vezzosi-simpliciales}. The functors $(-)_{\h\bT_\fil}$, $(-)^{\h\bT_\fil}$, and
$(-)^{\t\bT_\fil}$ are filtered analogues of the homotopy $S^1$-orbits, homotopy $S^1$-fixed
points, and $S^1$-Tate construction.

\begin{theorem}\label{thm:derhamhp}
    Let $k\rightarrow R$ be a morphism of derived commutative rings and let
    $\HH_\fil(R/k)$ be $\bT_\fil$-equivariant Hochschild homology of $R$ over
    $k$.
    There are natural filtrations
    $$\HH_\fil(R/k)_{\h\bT_\fil},\quad\HH_\fil(R/k)^{\h\bT_\fil},\quad\text{and}\quad\HH_\fil(R/k)^{\t\bT_\fil}$$
    with the following properties:
    \begin{enumerate}
        \item[{\em (i)}] the filtrations on $\HH_\fil(R/k)^{\h\bT_\fil}$ and
            $\HH_\fil(R/k)^{\t\bT_\fil}$ are multiplicative and
            $\HH_\fil(R/k)^{\h\bT_\fil}\rightarrow\HH_\fil(R/k)^{\t\bT_\fil}$ is a map of
            filtered $\bE_\infty$-algebras;
        \item[{\em (ii)}] the fiber sequence
            $$\HH_\fil(R/k)_{\h\bT_\fil}(1)[1]\rightarrow\HH_\fil(R/k)^{\h\bT_\fil}\rightarrow\HH_\fil(R/k)^{\t\bT_\fil}$$
            naturally admits the structure of a fiber sequence of
            $\HH_\fil(R/k)^{\h\bT_\fil}$-modules in $\FMod_k$.
        \item[{\em (iii)}] there are natural equivalences
            \begin{align*}
                \gr^s\HH_\fil(R/k)_{\h\bT_\fil}&\we\F^{0\leq\star\leq s}_\H\dRhat[2s],\\
                \gr^s\HH_\fil(R/k)^{\h\bT_\fil}&\we\F^s_\H\dRhat_{R/k}[2s],\\
                \gr^s\HH_\fil(R/k)^{\t\bT_\fil}&\we\dRhat_{R/k}[2s].
            \end{align*}
    \end{enumerate}
\end{theorem}

\begin{proof}
    The multiplicativity and
    module structure follows the multiplicativity of the Tate construction,
    which in this case is~\cite[Prop.~2.4.10]{raksit}. The identification of
    the graded pieces is~\cite[Sec.~6.3]{raksit}.
\end{proof}

\begin{remark}
    Note that the theorem does not say anything about the exhaustiveness or
    completeness of the filtration in general and indeed, in general, it is neither
    complete nor exhaustive. By~\cite{antieau-derham}, one does have completeness and exhaustiveness when $k$
    is static and $\L_{R/k}$ has Tor-amplitude in $[0,1]$.
    In characteristic $0$, a thorough discussion of exhaustiveness is provided by the work of Bals~\cite{bals}.
\end{remark}

\begin{variant}\label{variant:einfinityderham}
    There is an $\bE_\infty$-version of the results of this section and the next, at
    least when working over a base $\bE_\infty$-ring where the filtered circle
    $\bT_\fil$ exists as a cocommutative bialgebra in filtered modules. One
    such base is $\bZ$ or any $\bE_\infty$-algebra over $\bZ$. If $R\rightarrow S$ is a map of
    $\bE_\infty$-algebras over $\bZ$, then the $\bE_\infty$-algebra $\HH(S/R)$ is equipped with a filtration
    $\HH_\fil^{\bE_\infty}(S/R)$ where
    $\gr^\ast\HH_\fil(S/R)\we\Sym_S^\ast(\L_{S/R}^{\bE_\infty}[1])$ is the free graded
    $\bE_\infty$-algebra over $S$ on the $\bE_\infty$-cotangent complex
    $\L_{S/R}^{\bE_\infty}$ of $S$ over $R$. There is an associated $\bE_\infty$-version
    of the de Rham
    complex, which looks like
    $$S\rightarrow\L_{S/R}^{\bE_\infty}\rightarrow\Sym_S^2(\L_{S/R}^{\bE_\infty}[1])[-2]\rightarrow\cdots$$
    as a coherent cochain complex and admits a natural universal property as
    such. Because ordinary commutative rings have complicated
    $\bE_\infty$-cotangent complexes, and because the group homology of the
    symmetric groups appears in the symmetric powers appearing in the complex, it
    looks difficult to meaningfully compute the $\bE_\infty$-de Rham cohomology, even in simple
    cases.

    Note however that on the Hochschild side, if $R$ and $S$ are derived
    commutative rings, then
    the derived version of $\HH(S/R)$ agrees with the $\bE_\infty$-version, even
    as $\bE_\infty$-rings with $S^1$-action. It is only the HKR filtrations which
    differ in this case. There is a natural map
    $\HH_\fil^{\bE_\infty}(S/R)\rightarrow\HH_\fil(S/R)$ which on associated graded
    pieces induces the composite map
    $\Sym^\ast_S(\L_{S/R}^{\bE_\infty}[1])\rightarrow\Sym^\ast_S(\L_{S/R}[1])\rightarrow\LSym^\ast_S(\L_{S/R}[1])$.
    We study the $\bE_\infty$-de Rham and $\bE_\infty$-infinitesimal cohomology theories in
    F\&C.III.
\end{variant}

\section{Infinitesimal cohomology}\label{sec:inf}

Here is our definition of (derived) infinitesimal cohomology.

\begin{definition}[Infinitesimal cohomology]\label{def:inf}
    Fix a commutative ring $k$.
    Let $$\F^\star_\H\Inf_{-/k}\colon\DAlg_k\rightarrow\F\DAlg^\inf_k$$
    denote the left adjoint to the functor
    $$\DAlg_k\leftarrow\F\DAlg^\inf_k\colon\gr^0.$$ This is the
    Hodge-filtered
    \defidx{derived infinitesimal cohomology} functor. Similarly, let
    $\F^\star_\H\Infhat_{-/k}\colon\DAlg_k\rightarrow\widehat{\F\DAlg}^{\inf}_k$ be the completion, which is the left adjoint to 
    $\DAlg_k\leftarrow\widehat{\F\DAlg}^{\inf}_k\colon\gr^0$.
    We let $\Inf_{R/k}=\F^0_\H\Inf_{R/k}$ and $\Infhat_{R/k}=\F^0_\H\Infhat_{R/k}$.
\end{definition}

\begin{lemma}\label{lem:inf_nonnegative}
    For any $R\in\DAlg_k$, $\F^\star_\H\Inf_{R/k}$ is nonnegative.
\end{lemma}

\begin{proof}
    The right adjoint $\gr^0$ factors as a composition
    $\DAlg_k\xleftarrow{\gr^0}\F^+\DAlg_k\xleftarrow{(-)^+}\FDAlg_k$, where the first functor in
    the composition is the right adjoint (colocalization) to the fully faithful inclusion
    $\F^+\DAlg_k\hookrightarrow\FDAlg_k$ (obtained via left Kan extension along
    $\bN^\op\rightarrow\bZ^\op$.
\end{proof}

\begin{lemma}
    For any map $k\rightarrow R$ of derived commutative rings,
    there are natural equivalences
    $$\gr^i\Inf_{R/k}\we\LSym_R^i(\L_{R/k}[-1])$$
    for all $i\geq 0$.
\end{lemma}

\begin{proof}
    This follows as in the proof of Raksit's Theorem~\ref{thm:raksit}.
\end{proof}

\begin{corollary}
    The functor
    $\F^\star_\H\Inf_{-/k}$ is fully faithful.
\end{corollary}

\begin{proof}
    Indeed, the unit of the adjunction is $R\rightarrow\gr^0\Inf_{R/k}\we R$, which is an
    equivalence.
\end{proof}

\begin{remark}
    As a coherent cochain complex, the picture of $\F^\star_\H\Infhat_{R/k}$ is
    $$R\rightarrow\L_{R/k}\rightarrow\LSym_R^2(\L_{R/k}[-1])[2]\rightarrow\cdots.$$
\end{remark}

The following result holds for $\bQ$-algebras for derived de Rham cohomology by the Poincar\'e lemma
and gives an easy way to see why one must complete; see~\cite[Cor.~2.5]{bhatt-padic}.
The proposition gives a new proof even in characteristic $0$.

\begin{proposition}[Poincar\'e lemma]\label{prop:poincare}
    For any $k\rightarrow R$, the natural map $k\rightarrow\F^0_\H\Inf_{R/k}=\Inf_{R/k}$ is an
    equivalence.
\end{proposition}

\begin{proof}
    By Lemma~\ref{lem:inf_nonnegative}, the maps
    $$\F^0_\H\Inf_{R/k}\rightarrow\F^{-1}_\H\Inf_{R/k}\rightarrow\F^{-2}_\H\Inf_{R/k}\rightarrow\cdots$$
    are all equivalences.
    In general, taking $\F^{-\infty}=\colim_{n\mapsto-\infty}\F^n$ gives a map $\FDAlg_k\rightarrow\DAlg_k$ which is left
    adjoint to a functor $\FDAlg_k\leftarrow\DAlg_k$ which takes a derived commutative
    ring $R$ in $\Mod_k$ to the constant filtration $\cdots=R=R=R=\cdots$. The composition
    $$\DAlg_k\xrightarrow{\F^\star_\H\Inf(-/k)}\FDAlg_k\xrightarrow{\F^0}\DAlg_k$$
    is therefore left adjoint to the functor which takes $R$ to $\gr^0$ of the constant filtration
    on $R$, which is the $0$ ring. Thus, the total left adjoint takes any $R$ to the initial derived
    commutative ring, i.e., $k$ itself.
\end{proof}

\begin{remark}
    Achim Krause points out that alternatively one can compute directly the infinitesimal
    cohomology of free derived commutative $k$-algebras. If $R=\LSym_k(M)$, then arguing as in
    Theorem~\ref{thm:raksit} we find that $\F^\star_\H\Inf_{R/k}$ is equivalent to the free infinitesimal
    filtered derived commutative ring on the filtration $$\underline{M}=\cdots 0\rightarrow M[-1]\rightarrow
    0\rightarrow\cdots,$$ where the only non-vanishing term $M[-1]$ occurs in filtration weight
    $1$. By the usual formula for the Day convolution tensor product, it follows that
    $\F^0(\underline{M}^{\otimes i})\we 0$ for $i>0$ and hence the same is true of
    $\F^0((\underline{M}^{\otimes i})_{\h\Sigma_i})$ for $i>0$.
\end{remark}

We thus view the completed version $\Infhat_{R/k}$ as the fundamental form of infinitesimal
cohomology and then it becomes an interesting question of when the Hodge filtration on $\Inf_{R/k}$
is complete. These questions are considered in detail in F\&C.II.

\section{Infinitesimal cohomology in characteristic $p$}\label{sec:charp}

The results of this section were pointed out to me by Akhil Mathew and independently discovered by
Jiaqi Fu.

\begin{lemma}\label{lem:infcoh_qrsp}
    Let $k$ be a perfect commutative $\bF_p$-algebra, and let $R$ be a quasiregular semiperfect
    $k$-algebra with inverse limit perfection $R^\flat$. Let $I=\ker(R^\flat\rightarrow R)$.
    Then, there is a natural equivalence $I^\star R^\flat\rightarrow\F^\star_\H\Infhat_{R/k}$.
    In particular, $R^\flat\we\Infhat_{R/k}$.
\end{lemma}

\begin{proof}
    As $\L_{R^\flat/k}\we 0$, it follows that $\Infhat_{R^\flat/k}\xrightarrow{\gr^0} R^\flat$ is
    an equivalence. Thus, functoriality induces a natural map
    $R^\flat\we\Infhat_{R^\flat/k}\rightarrow\Infhat_{R/k}$ and since $R^\flat$ is perfect,
    $\F^\star_\H\Infhat_{R/k}\we\F^\star_\H\Infhat_{R/R^\flat}$. Since $\L_{R/R^\flat}\we I/I^2[1]$
    is a shift of a flat module, it follows that
    $\LSym^n_R(\L_{R/R^\flat}[-1])\we\LSym^n_R(I/I^2)$ is flat for all $n\geq 0$. In particular, by
    induction, $\F^\star\Infhat_{R/R^\flat}$ is a strictly filtered commutative $R^\flat$-algebra. We obtain an
    induced map $I^\star R^\flat\rightarrow\F^\star\Infhat_{R/R^\flat}$. On the other hand, there
    are natural surjective maps $\LSym^n(I/I^2)\rightarrow I^n/I^{n+1}$ for all $n\geq 0$ arising
    from the classical Rees construction.
    The composition $I^n/I^{n+1}\rightarrow\LSym^n(I/I^2)\rightarrow
    I^n/I^{n+1}$ is the identity. Thus, $I^\star R^\flat\iso\F^\star_\H\Infhat_{R/R^\flat}$.
\end{proof}

\begin{remark}
    Let $S$ be a commutative ring with an ideal $J$.
    In general, $\pi_0\LSym^n(J/J^2)\rightarrow J^n/J^{n+1}$ is not an injection. (When this
    happens is a question studied in~\cite{huneke}.) Unwinding the
    argument above, it follows that, without the flatness condition on $J/J^2$, one does not
    typically expect to have a map $J^\star S\rightarrow\Infhat_{(S/J)/S}$.
\end{remark}

\begin{theorem}[Perfection comparison]\label{thm:perfection_comparison}
    Suppose that $k$ is a perfect commutative $\bF_p$-algebra and that $R$ is a quasisyntomic
    $k$-algebra. Then, there is a natural equivalence $R^\flat\we\Infhat_{R/k}$.
\end{theorem}

\begin{proof}
    The assignment $R\mapsto R^\flat$ satisfies quasisyntomic descent. This is also true of
    $R\mapsto\LSym_R^n(\L_{R/k}[-1])$ for all $n\geq 0$ by a variant of the argument
    in~\cite{bms2} for quasisyntomic descent for the cotangent complex and its exterior powers. (A
    proof will be given in joint forthcoming work with Ryomei Iwasa and Achim Krause.) Thus, it is true of $R\mapsto\Infhat_{R/k}$.
    Now, the theorem follows from Lemma~\ref{lem:infcoh_qrsp} via unfolding in the sense
    of~\cite[Sec.~4]{bms2}.
\end{proof}

\begin{remark}
    Jiaqi Fu has pointed out to us that one can remove the quasisyntomicity condition of
    Theorem~\ref{thm:perfection_comparison} at the cost of imposing some finiteness conditions. We
    will return to this in F\&C.II after developing further tools for analyzing
    $\Infhat$.
\end{remark}

\begin{corollary}\label{cor:charp}
    If $R$ is a $p$-completely flat quasisyntomic $\bZ_p$-algebra, then there is a natural $p$-adic
    equivalence $\bA_\inf(R)\we\Infhat_{R/\bZ_p}$, where $\bA_\inf(R)=W((R/p)^\flat)$.
\end{corollary}

\begin{proof}
    It suffices to check modulo $p$, where it follows from Theorem~\ref{thm:perfection_comparison}.
\end{proof}

\begin{remark}
    In the setting of Corollary~\ref{cor:charp}, if $R$ is the $p$-completion of a
    finite type $\bZ_{(p)}$-algebra, then $\L_{R/\bZ_p}$ and all of its exterior powers are
    $p$-complete. It follows that $\Infhat_{R/\bZ_p}$ is $p$-complete and the $p$-completion in the
    statement of the corollary is unnecessary.
\end{remark}

\begin{remark}
    The fiber $I$ of $R^\flat\rightarrow R$ also satisfies quasisyntomic descent on
    $\QSyn_{\bF_p}$. Unfolding gives a formula for $\F^1_\H\Infhat_{R/k}$ and similarly
    $\F^n_\H\Infhat_{R/k}$ is obtained by unfolding $I^n$ from the qrsp case.
\end{remark}

\begin{remark}[Grothendieck's infinitesimal site]
    Finally, we justify the name infinitesimal cohomology. If $R$ is a smooth commutative
    $k$-algebra where $k$ is a perfect field, then $\Infhat_{R/k}$ is equivalent to the
    cohomology of the infinitesimal site relative to $k$ introduced by Grothendieck
    in~\cite{grothendieck-crystals}. We also recover the results of
    Ogus~\cite{ogus_infinitesimal}. For example, if $X$ is a smooth proper $k$-scheme, then
    $\Infhat_{X/k}=\R\Gamma(X,\Infhat_{\Oscr_X})\we\R\Gamma(X,\Oscr_X^\perf)\we\R\Gamma(X,\Oscr)^\perf\we\R\Gamma(X,k)$.
\end{remark}

\section{The crystallization functor}\label{sec:crystallization}

The object $\bD_-^\vee\in\Gr^-\Mod_k$ is the free $[0]$-sheared graded derived commutative ring on
$\1_\gr[1](-1)$; it inherits its bialgebra structure from the canonical diagonal
$\1_\gr[1](-1)\rightarrow\1_\gr[1](-1)\oplus\1_\gr[1](-1)$. Thus, it defines an
$\bE_\infty$-coalgebra in $\Gr^-\DAlg^{\N[0]}_k$. Using Variant~\ref{var:nonpositive}, this forgets to
$\bE_\infty$-coalgebra structures on $\bD_-^\vee$ in
$$\Gr^-\CAlg_k,\quad\Gr^-\DAlg^{\N[-1]}_k,\quad\Gr^-\DAlg^{\B[0]}_k,\quad\text{and}\quad\Gr^-\DAlg^{\B[0]\s}_k$$
as well.

\begin{definition}
    Let $k$ be a commutative ring.
    \begin{enumerate}
        \item[(a)]   A complete filtered $\bE_\infty$-$k$-algebra is an object of
            $\cMod_{\bD_-^\vee}(\Gr\CAlg_k)$, which is denoted below as $\widehat{\F\CAlg}_k$.
        \item[(b)]   A complete infinitesimal filtered derived commutative $k$-algebra is an object
            of $\cMod_{\bD_-^\vee}(\Gr\DAlg^{\N[0]}_k)$, denoted below by $\widehat{\F\DAlg}^\inf_k$.
        \item[(c)]   A \longdefidx{coherent cdga}{filtered derived commutative ring!coherent cdga}
            over $k$ is an object of $\cMod_{\bD_-^\vee}(\Gr\DAlg^{\B[0]}_k)$, denoted below by
            $\widehat{\F\DAlg}^{\mathrm{cdga}}_k$.
        \item[(d)]   A \longdefidx{strict coherent cdga}{filtered derived commutative ring!strict
            coherent cdga} over $k$ is an object of $\cMod_{\bD_-^\vee}(\Gr\DAlg^{\B[0]\s}_k)$,
            denoted below by
            $\widehat{\F\DAlg}^{\mathrm{scdga}}_k$.
        \item[(e)]   A complete crystalline filtered derived commutative $k$-algebra is an object
            of $\cMod_{\bD_-^\vee}(\Gr\DAlg^{\N[-1]}_k)$, denoted below by $\widehat{\F\DAlg}^\crys_k$.
    \end{enumerate}
    A filtered commutative $k$-algebra in of any of the above flavors is
    \longdefidx{nonnegative}{filtered commutative ring!nonnegative} if the underlying graded
    object is in $\Gr^+\Mod_k$. It is \longdefidx{nonpositive}{filtered commutative ring!nonnegative} if the underlying graded
    object is in $\Gr^-\Mod_k$. There are then variants on the notation such as
    $\widehat{\F^+\DAlg}^\crys_k$, etc.
\end{definition}

The main result of this section is a lift of Theorem~\ref{thm:gradedcrystallization} to
the filtered context.

\begin{theorem}[Crystallization]\label{thm:filteredcrystallization}
    Let $k$ be a commutative ring.
    When restricted to nonnegative complete objects, there are forgetful functors
    \begin{equation}\label{eq:grdalg}
        \widehat{\F^+\CAlg}_k\leftarrow\widehat{\F^+\DAlg}^\inf_k\leftarrow\widehat{\F^+\DAlg}^{\mathrm{cdga}}_k\leftarrow\widehat{\F^+\DAlg}^{\mathrm{scdga}}_k\leftarrow\widehat{\F^+\DAlg}^\crys_k
    \end{equation}
    which are conservative and preserve limits and sifted colimits and which commute with the
    forgetful functors to $\widehat{\F^+\Cscr}$.
\end{theorem}

\begin{proof}
    We prove the existence of $\widehat{\F^+\DAlg}^\inf_k\leftarrow\widehat{\F^+\DAlg}^\crys_k$
    and leave the other cases to the reader. Once constructed, the rest
    of the properties are routine to check and similar to the case of
    Theorem~\ref{thm:gradedcrystallization}.
    Our proof is inspired by the proof of that theorem and the use of
    Proposition~\ref{prop:functoriality}. However, this proposition does not quite apply in the
    situation we have here, so we make a more complicated argument.

    We have by construction a forgetful functor $\Gr\CAlg^\crys_k\leftarrow\Gr\DAlg^\crys_k$ and
    taking $\bD_-^\vee$-comodules produces a forgetful functor
    $\widehat{\FCAlg}_k\leftarrow\widehat{\FDAlg}_k^\crys$. We can restrict this to nonnegatively
    graded objects to obtain $\widehat{\F^+\CAlg}_k\leftarrow\widehat{\F^+\DAlg}_k^\crys$. We will
    show that this functor factors canonically through $\widehat{\F^+\DAlg}_k^\inf$.

    To this end, note that $$\Gr^+\Mod_k\we\lim_{N\mapsto\infty}\Gr^{[0,N-1]}\Mod_k,$$ where
    $\Gr^{[0,N-1]}\Mod_k\subseteq\Gr^+\Mod_k$ is the full subcategory of objects concentrated in
    weights $0,\ldots,N-1$. We call $\Gr^{[0,N-1]}\Mod_k$ the $\infty$-category of $N$-stubs in
    $\Mod_k$. It is a symmetric monoidal localization of $\Gr^+\Mod_k$. The tensor product of two
    $N$-stubs is obtained by computing the tensor product in $\Gr^+\Mod_k$ and then forgetting
    everything in weights $N$ and more. We also have, for each $N\geq 0$, an $\infty$-category of crystalline
    derived commutative $k$-algebras in $N$-stubs obtained as a pullback
    $$\xymatrix{
        \Gr^{[0,N-1]}\DAlg_k^\crys\ar[r]\ar[d]&\Gr^+\DAlg_k^\crys\ar[d]\\
        \Gr^{[0,N-1]}\Mod_k\ar[r]&\Gr^+\Mod_k.
    }$$
    There is an $\infty$-category $\Gr^{[0,N-1]}\CAlg_k$ defined in the same way, which agrees with
    the $\infty$-category of $\bE_\infty$-algebras in $N$-stubs. And, we have
    $\Gr^{[0,N-1]}\DAlg_k^\inf$ as well. The limiting property on modules extends to categories of
    algebras as well:
    \begin{align*}
        \Gr^+\CAlg_k &\we \lim_N\Gr^{[0,N-1]}\CAlg_k,\\
        \Gr^+\DAlg_k^\inf &\we \lim_N\Gr^{[0,N-1]}\DAlg_k^\inf,\\
        \Gr^+\DAlg_k^\crys &\we \lim_N\Gr^{[0,N-1]}\DAlg_k^\crys.
    \end{align*}
    Taking $\bD_-^\vee$-comodules we obtain
    \begin{align*}
        \cMod_{\bD_-^\vee}(\Gr^+\CAlg_k) &\we \lim_N\cMod_{\bD_-^\vee}(\Gr^{[0,N-1]}\CAlg_k),\\
        \cMod_{\bD_-^\vee}(\Gr^+\DAlg_k^\inf) &\we
        \lim_N\cMod_{\bD_-^\vee}(\Gr^{[0,N-1]}\DAlg_k^\inf),\\
        \cMod_{\bD_-^\vee}(\Gr^+\DAlg_k^\crys) &\we \lim_N\cMod_{\bD_-^\vee}(\Gr^{[0,N-1]}\DAlg_k^\crys).
    \end{align*}
    These are precisely the $\infty$-categories $\widehat{\F^+\CAlg}_k$,
    $\widehat{\F^+\DAlg}_k^\inf$, and $\widehat{\F^+\DAlg}_k^\crys$.

    Now, it suffices to factor the forgetful functor
    $\cMod_{\bD_-^\vee}(\Gr^{[0,N-1]}\CAlg_k)\leftarrow\cMod_{\bD_-^\vee}(\Gr^{[0,N-1]}\DAlg_k^\crys)$
    through $\cMod_{\bD_-^\vee}(\Gr^{[0,N-1]}\DAlg_k^\inf)$ compatibly in $N$.
    Fix $N\geq 0$ and $0\leq i\leq N$. Let $P$ be a projective $k$-module and let $\ins^iP$
    denote the object of $\cMod_{\bD_-^\vee}(\Gr^{[0,N-1]}\Mod_k)$ corresponding to the filtration
    $$\cdots\rightarrow 0\rightarrow P=P=P=\cdots,$$ where the first non-zero term occurs in filtration weight
    $P$. The free object in $\cMod_{\bD_-^\vee}(\Gr^{[0,N-1]}\DAlg_k^\crys)$ on $\ins^iP$ is
    $$\bigoplus_{ir\leq N}\ins^{ir}(\LSym^r(P[2i])[-2ir]).$$ By Lemma~\ref{lem:toramplitude}, this
    object is coconnective. It follows that the natural map $$\bigoplus_{ir\leq
    N}\ins^{ir}(\Sym^r(P))\rightarrow\bigoplus_{ir\leq N}\ins^{ir}(\LSym^r(P[2i])[-2ir])$$
    factors through $$\bigoplus_{ir\leq
    N}\ins^{ir}(\Sym^r(P))\rightarrow\bigoplus_{ir\leq N}\ins^{ir}(\LSym^r(P)).$$
    As a result, thanks to Theorem~\ref{thm:apparatus}, the map of monads on $\Sym_k\rightarrow\LSym^\crys_k$
    $\cMod_{\bD_-^\vee}(\Gr^{[0,N-1]}\Mod_k)$ factors through $\LSym_k^\inf$, the $1$-derived
    approximation of $\Sym$. This proves the main case of the theorem.
\end{proof}

\begin{remark}
    In the proof of Theorem~\ref{thm:filteredcrystallization}, we use the notion of $N$-stubs simply as a
    way of working around the fact that the $\infty$-category of connective objects with respect to the neutral $t$-structure on $\widehat{\F^+\DAlg}_k$ is not
    compactly generated, so the results of Section~\ref{sec:polynomial_monad} do not apply directly.
\end{remark}

\begin{corollary}\label{cor:crystallization}
    Let $k$ be a derived commutative ring. There are adjunctions
    $$\divideontimes\colon\F^+\DAlg^\inf_k\rightleftarrows\F^+\DAlg^\crys_k\colon\text{forget}$$
    and
    $$\divideontimes\colon\widehat{\F^+\DAlg}^\inf_k\rightleftarrows\widehat{\F^+\DAlg}^\crys_k\colon\text{forget}.$$
\end{corollary}

\begin{definition}
    The left adjoints appearing in Corollary~\ref{cor:crystallization} are called the filtered crystallization functors.
\end{definition}

\begin{corollary}
    Let $k$ be a derived commutative ring.
    If $\F^\star R\in\F^+\DAlg^\crys_k$ or $\widehat{\F^+\DAlg}^\crys_k$, then the
    underlying filtered $\bE_\infty$-ring $\F^\star R$ is naturally an infinitesimal filtered derived commutative
    ring and the underlying $\bE_\infty$-ring of $\F^0R$ is naturally a derived commutative ring.
\end{corollary}

The following result gives a decompletion of one of the main points of Raksit's work.

\begin{proposition}\label{prop:decompletedraksit}
    Suppose that $k$ is a commutative ring and
    consider the adjunction
    $$\F^\star_\H\dR_{R/k}\colon\DAlg_k\rightleftarrows\F\DAlg^\crys_k\colon\gr^0.$$
    The left adjoint identifies with the Hodge-filtered derived de Rham
    cohomology when forgetting to $\F\CAlg_k$.
\end{proposition}

\begin{proof}
    By completing, Theorem~\ref{thm:raksit} implies that the completion of the left
    adjoint gives Hodge-filtered Hodge-complete derived de Rham cohomology.
    Since Hodge-filtered derived de Rham cohomology is left Kan extended from
    polynomial algebras, it is enough to see that the left adjoint takes the
    correct value on $k[x]$ (where $x$ has degree $0$) and specifically that the result is
    already complete. This follows from the construction of the decompleted $\LSym$ functor since it
    is left Kan extended from its values on compact objects of $\FMod_k$ as the left adjoint to
    $\Mod_k\leftarrow\FMod_k\colon\gr^0$ takes $k$ to a compact object.
\end{proof}

\begin{corollary}
    For any map $k\rightarrow R$ of derived commutative rings, $\dR_{R/k}$ and $\dRhat_{R/k}$ are
    derived commutative $k$-algebras.
\end{corollary}

\begin{theorem}[Infinitesimal--de Rham comparison]
    If $R\in\DAlg_k$, then
    $$\divideontimes(\F^\star_\H\Inf_{R/k})\we\F^\star_\H\dR_{R/k}$$ and
    $$\divideontimes(\F^\star_\H\Infhat_{R/k})\we\F^\star_\H\dRhat_{R/k}.$$
\end{theorem}

\begin{proof}
    On $\F^+\DAlg^\crys_k$, the $\gr^0$ functor factors through the
    forgetful functor to $\F^+\DAlg^\inf_k$ and similarly in the
    complete case. The result now follows from composition of left adjoints.
\end{proof}

By adjunction, there is a canonical map
$$\F^\star_\H\Inf_{R/k}\rightarrow\F^\star_\H\dR_{R/k}$$ of infinitesimal
filtered derived commutative $k$-algebras which on associated
graded pieces induces canonical maps
$$\LSym_R^i(\L_{R/k}[-1])\rightarrow\Lambda^i\L_{R/k}[-i].$$
For $\bQ$-algebras $R$, these maps are equivalences, but in general there is a
significant difference between the theories.

\begin{examples}
    \begin{enumerate}
        \item[(1)] It is easy to see that $(\bZ,(p))$ is the universal pair with a map
            $\bF_p\rightarrow\gr^0(\bZ,(p))$. It follows that
            $\F^\star_\H\Inf_{\bF_p/\bZ}\we p^\star\bZ$ and
            $\divideontimes(p^\star\bZ)\we\F^\star_\H\dR_{\bF_p/\bZ}\we
            \bZ\langle x\rangle/(x-p)$ with the pd filtration.
        \item[(2)] More generally, if $I\subseteq R$ is an ideal locally generated by a regular
            sequence and $S=R/I$, then,
            $$\F^\star_\H\Inf_{S/R}\we \LSym^\star_RI$$ and
            the completion is $\widehat{\Inf}_{S/R}\we R_I^\wedge$.
            On the other hand, the Hodge-complete derived de Rham cohomology
            $\F^\star_\H\dRhat_{S/R}$ is the pd envelope $\widehat{D}_R(I)$ completed with
            respect to the pd filtration.
        \item[(3)] If $R$ is a perfectoid $\bZ_p$-algebra, isomorphic to $\bA_\inf(R)/(\xi)$, then there is a
            $p$-adic equivalence
            $$\F^\star_\H\Infhat_{R/\bZ_p}\we\xi^\star\bA_\mathrm{inf}(R)$$
            and the natural map
            $\F^\star\Infhat_{R/\bZ_p}\rightarrow\F^\star\dRhat_{R/\bZ_p}$ is
            identified after $p$-completion with
            $\xi^\star\bA_{\mathrm{inf}}(R)\rightarrow\F^\star_{\mathrm{pd}}\widehat{\bA}_{\mathrm{crys}}(R)$,
            the completion of $\bA_\crys(R)$ with respect to the pd filtration.
    \end{enumerate}
\end{examples}

\begin{construction}
    Given the classical construction of de Rham cohomology as a cdga, it is possibly more natural to use
    strict coherent cdgas as a basis for de Rham cohomology, rather than crystalline filtered derived
    commutative rings. However, the resulting theories are the equivalent.
    Let $\L\Omega^\bullet_{-/k}$ denote the left adjoint to
    $\gr^0\colon\widehat{\F\DAlg}^{\mathrm{scdga}}_k\rightarrow\DAlg_k$
    and let
    $\div\colon\widehat{\F^+\DAlg}^{\mathrm{scdga}}_k\rightarrow\widehat{\F^+\DAlg}^\crys_k$ be the free crystalline
    filtered derived commutative functor.
\end{construction}

\begin{theorem}
    Fix $R\in\DAlg_k$.
    \begin{enumerate}
        \item[{\em (a)}] There is a natural equivalence
            $\div\L\Omega^\bullet_{R/k}\we\F^\star_\H\dRhat_{R/k}$.
        \item[{\em (b)}] The unit map
            $\L\Omega^\bullet_{R/k}\rightarrow\F^\star_\H\dRhat_{R/k}$ of
            underlying coherent strict coherent cdgas is an equivalence.
    \end{enumerate}
\end{theorem}

\begin{proof}
    This follows because the map of monads from the free strict coherent cdga monad to the free crystalline filtered derived
    commutative ring monad is an equivalence on complete filtered objects concentrated in weights $0$ and
    $1$, so (a) and (b) are true on free derived commutative rings and hence on all by taking
    resolutions by free objects.
\end{proof}

\small
\bibliographystyle{amsplain}
\bibliography{infcoh}
\addcontentsline{toc}{section}{References}

\medskip
\noindent
\textsc{Department of Mathematics, Northwestern University}\\
{\ttfamily antieau@northwestern.edu}

\end{document}

%% file: preamble.tex
\usepackage[pdfstartview=FitH,
            pdfauthor={},
            pdftitle={},
            colorlinks,
            linkcolor=reference,
            citecolor=citation,
            urlcolor=e-mail,
            backref]{hyperref}

\usepackage{amsmath}
\usepackage{amscd}
\usepackage{amsbsy}
\usepackage{amssymb}
\usepackage{verbatim}
\usepackage{eufrak}
\usepackage{eucal}
\usepackage{microtype}
\usepackage{hyperref}
\usepackage{mathrsfs}
\usepackage{amsthm}
\usepackage{stmaryrd}
\usepackage{wasysym}
\usepackage{xy}
\input xy
\xyoption{all}
\usepackage{tikz}
\usetikzlibrary{matrix,arrows}
\usepackage{tikz-cd}
\usepackage{dsfont}
\usepackage[T1]{fontenc}
\usepackage{spectralsequences}

\usepackage{fancyhdr}
\usepackage{makeidx}
\makeindex
\newcommand{\defidx}[1]{#1\index{#1|textbf}}
\newcommand{\longdefidx}[2]{#1\index{#2|textbf}}

\newcommand{\stackspace}{2.5}
\newcommand{\stack}[2][1cm]{\;\tikz[baseline, yshift=.65ex]%
    {\foreach \k [evaluate=\k as \r using (.5*#2+.5-\k)*\stackspace] in {1,...,#2}{%
    \ifodd\k{\draw[->](0,\r pt)--(#1,\r pt);}%
    \else{\draw[<-](0,\r pt)--(#1,\r pt);}\fi
    }}\;}

\usepackage[bbgreekl]{mathbbol}
\usepackage{amsfonts}
\DeclareSymbolFontAlphabet{\mathbb}{AMSb} 
\DeclareSymbolFontAlphabet{\mathbbl}{bbold}

\newcommand{\Inf}{{\mathbbl{\Pi}}}
\newcommand{\Infhat}{\widehat{\Inf}}

\usepackage{color}
\definecolor{todo}{rgb}{1,0,0}
\definecolor{conditional}{rgb}{0,1,0}
\definecolor{e-mail}{rgb}{0,.40,.80}
\definecolor{reference}{rgb}{.20,.60,.22}
\definecolor{mrnumber}{rgb}{.80,.40,0}
\definecolor{citation}{rgb}{0,.40,.80}

\let\oldmarginpar\marginpar
\renewcommand\marginpar[1]{\-\oldmarginpar[\raggedleft\footnotesize #1]%
{\raggedright\footnotesize #1}}

\newcommand{\Ascr}{\mathcal{A}}

\newcommand{\Cscr}{\mathcal{C}}
\newcommand{\Dscr}{\mathcal{D}}
\newcommand{\Escr}{\mathcal{E}}
\newcommand{\Fscr}{\mathcal{F}}

\newcommand{\Mscr}{\mathcal{M}}

\newcommand{\Oscr}{\mathcal{O}}
\newcommand{\Pscr}{\mathcal{P}}

\newcommand{\Sscr}{\mathcal{S}}

\newcommand{\Xscr}{\mathcal{X}}

\newcommand{\B}{\mathrm{B}}

\renewcommand{\d}{\mathrm{d}}
\newcommand{\D}{\mathrm{D}}
\newcommand{\E}{\mathrm{E}}
\newcommand{\F}{\mathrm{F}}

\renewcommand{\H}{\mathrm{H}}
\newcommand{\h}{\mathrm{h}}

\newcommand{\K}{\mathrm{K}}
\renewcommand{\L}{\mathrm{L}}

\newcommand{\N}{\mathrm{N}}
\renewcommand{\P}{\mathrm{P}}
\newcommand{\R}{\mathrm{R}}
\renewcommand{\r}{\mathrm{r}}

\newcommand{\s}{\mathrm{s}}

\renewcommand{\t}{\mathrm{t}}

\renewcommand{\1}{\mathbf{1}}
\newcommand{\bA}{\mathbf{A}}

\newcommand{\bD}{\mathbf{D}}
\newcommand{\bE}{\mathbf{E}}
\newcommand{\bF}{\mathbf{F}}
\newcommand{\bG}{\mathbf{G}}

\newcommand{\bN}{\mathbf{N}}

\newcommand{\bQ}{\mathbf{Q}}
\newcommand{\bR}{\mathbf{R}}
\newcommand{\bS}{\mathbf{S}}
\newcommand{\bT}{\mathbf{T}}

\newcommand{\bZ}{\mathbf{Z}}

\newcommand{\Cat}{\Cscr\mathrm{at}}

\newcommand{\Op}{\Oscr\mathrm{p}}

\newcommand{\Sp}{\Sscr\mathrm{p}}

\newcommand{\op}{\mathrm{op}}
\newcommand{\cofib}{\mathrm{cofib}}

\newcommand{\cn}{\mathrm{cn}}

\newcommand{\ap}{\mathrm{ap}}
\newcommand{\ep}{\mathrm{ep}}

\newcommand{\Mod}{\mathrm{Mod}}
\newcommand{\LMod}{\mathrm{LMod}}
\newcommand{\FMod}{\mathrm{FMod}}
\newcommand{\GrMod}{\mathrm{GrMod}}

\newcommand{\Ind}{\mathrm{Ind}}

\newcommand{\Alg}{\mathrm{Alg}}
\newcommand{\CAlg}{\mathrm{CAlg}}
\newcommand{\cCAlg}{\mathrm{cCAlg}}
\newcommand{\cDAlg}{\mathrm{cDAlg}}

\newcommand{\cMod}{\mathrm{cMod}}

\newcommand{\DAlg}{\mathrm{DAlg}}

\newcommand{\FDAlg}{\mathrm{FDAlg}}
\newcommand{\FCAlg}{\mathrm{FCAlg}}
\newcommand{\GrDAlg}{\mathrm{GrDAlg}}

\newcommand{\Gr}{\mathrm{Gr}}
\newcommand{\gr}{\mathrm{gr}}
\newcommand{\ev}{\mathrm{ev}}
\newcommand{\ins}{\mathrm{ins}}
\newcommand{\fil}{\mathrm{fil}}

\newcommand{\LSym}{\mathrm{LSym}}
\newcommand{\Sym}{\mathrm{Sym}}

\newcommand{\LAntiSym}{\mathrm{LAntiSym}}

\newcommand{\mot}{\mathrm{mot}}

\newcommand{\QSyn}{\mathrm{QSyn}}

\newcommand{\tensorhat}{\widehat{\otimes}}

\newcommand{\HKR}{\mathrm{HKR}}
\newcommand{\perf}{\mathrm{perf}}

\newcommand{\heart}{\heartsuit}

\renewcommand{\geq}{\geqslant}
\renewcommand{\leq}{\leqslant}
\newcommand{\Ho}{\mathrm{Ho}}

\newcommand{\HC}{\mathrm{HC}}

\newcommand{\HH}{\mathrm{HH}}

\renewcommand{\inf}{\mathrm{inf}}
\newcommand{\crys}{\mathrm{crys}}

\newcommand{\dR}{\mathrm{dR}}
\newcommand{\dRhat}{\widehat{\dR}}

\newcommand{\Ch}{\mathrm{Ch}}

\newcommand{\Map}{\mathrm{Map}}

\newcommand{\Hom}{\mathrm{Hom}}
\newcommand{\Fun}{\mathrm{Fun}}
\newcommand{\End}{\mathrm{End}}

\newcommand{\Ga}{\bG_{a}}

\DeclareMathOperator*{\colim}{colim}

\DeclareMathOperator*{\Tot}{Tot}

\DeclareMathOperator{\Spec}{Spec}

\newcommand{\we}{\simeq}
\newcommand{\iso}{\cong}

\theoremstyle{plain}
\newtheorem{theorem}{Theorem}[section]
\newtheorem*{theorem*}{Theorem}
\newtheorem{lemma}[theorem]{Lemma}

\newtheorem{proposition}[theorem]{Proposition}

\newtheorem{corollary}[theorem]{Corollary}
\newtheorem*{corollary*}{Corollary}

\theoremstyle{plain}

\theoremstyle{definition}

\newtheoremstyle{named}{}{}{\itshape}{}{\bfseries}{.}{.5em}{#1 \thmnote{#3}}
\theoremstyle{named}

\theoremstyle{definition}

\newtheorem{definition}[theorem]{Definition}
\newtheorem{warning}[theorem]{Warning}
\newtheorem{variant}[theorem]{Variant}
\newtheorem{notation}[theorem]{Notation}

\newtheorem{example}[theorem]{Example}
\newtheorem*{example*}{Example}
\newtheorem{examples}[theorem]{Examples}

\newtheorem*{question*}{Question}
\newtheorem{construction}[theorem]{Construction}

\newtheorem{remark}[theorem]{Remark}

%% file: infcoh.bbl
\providecommand{\bysame}{\leavevmode\hbox to3em{\hrulefill}\thinspace}
\providecommand{\MR}{\relax\ifhmode\unskip\space\fi MR }
\providecommand{\MRhref}[2]{%
  \href{http://www.ams.org/mathscinet-getitem?mr=#1}{#2}
}
\providecommand{\href}[2]{#2}
\begin{thebibliography}{10}

\bibitem{antieau-derham}
Benjamin Antieau, \emph{Periodic cyclic homology and derived de {R}ham
  cohomology}, Ann. K-Theory \textbf{4} (2019), no.~3, 505--519. \MR{4043467}

\bibitem{antieau_whatis}
\bysame, \emph{What is de {R}ham cohomology?}, Oberwolfach Reports \textbf{35}
  (2022), 2027--2034.

\bibitem{antieau-spherical}
\bysame, \emph{Spherical {W}itt vectors and integral models for spaces}, arXiv
  preprint arXiv:2308.07288 (2023).

\bibitem{antieau-decalage}
\bysame, \emph{Spectral sequences, d\'ecalage, and the {B}eilinson
  t-structure}, arXiv preprint arXiv:2411.09115 (2024).

\bibitem{antieau-riggenbach}
Benjamin Antieau and Noah Riggenbach, \emph{Cyclotomic synthetic spectra},
  arXiv preprint arXiv:2411.19929 (2024).

\bibitem{ariotta}
Stefano Ariotta, \emph{Coherent cochain complexes and {B}eilinson t-structures,
  with an appendix by {A}chim {K}rause}, arXiv preprint arXiv:2109.01017
  (2021).

\bibitem{bals}
Konrad Bals, \emph{Periodic cyclic homology over {$\mathbf{Q}$}}, Ann. K-Theory
  \textbf{9} (2024), no.~1, 119--142. \MR{4756830}

\bibitem{bgmn}
Clark Barwick, Saul Glasman, Akhil Mathew, and Thomas Nikolaus, \emph{K-theory
  and polynomial functors}, arXiv preprint arXiv:2102.00936 (2021).

\bibitem{bzn}
David Ben-Zvi and David Nadler, \emph{Loop spaces and connections}, J. Topol.
  \textbf{5} (2012), no.~2, 377--430. \MR{2928082}

\bibitem{bhatt-completions}
Bhargav Bhatt, \emph{Completions and derived de {R}ham cohomology}, arXiv
  preprint arXiv:1207.6193 (2012).

\bibitem{bhatt-padic}
\bysame, \emph{p-adic derived de {R}ham cohomology}, arXiv preprint
  arXiv:1204.6560 (2012).

\bibitem{bms2}
Bhargav Bhatt, Matthew Morrow, and Peter Scholze, \emph{Topological
  {H}ochschild homology and integral {$p$}-adic {H}odge theory}, Publ. Math.
  Inst. Hautes \'{E}tudes Sci. \textbf{129} (2019), 199--310. \MR{3949030}

\bibitem{brantner-campos-nuiten}
Lukas Brantner, Ricardo Campos, and Joost Nuiten, \emph{{PD} operads and
  explicit partition {L}ie algebras}, arXiv preprint arXiv:2104.03870 (2021).

\bibitem{brantner-mathew}
Lukas Brantner and Akhil Mathew, \emph{Deformation theory and partition {L}ie
  algebras}, arXiv preprint arXiv:1904.07352 (2019).

\bibitem{dold-puppe}
Albrecht Dold and Dieter Puppe, \emph{Homologie nicht-additiver {F}unktoren.
  {A}nwendungen}, Ann. Inst. Fourier (Grenoble) \textbf{11} (1961), 201--312.
  \MR{150183}

\bibitem{dwyer-greenlees-iyengar-exterior}
W.~G. Dwyer, J.~P.~C. Greenlees, and S.~B. Iyengar, \emph{D{G} algebras with
  exterior homology}, Bull. Lond. Math. Soc. \textbf{45} (2013), no.~6,
  1235--1245. \MR{3138491}

\bibitem{eml2}
Samuel Eilenberg and Saunders Mac~Lane, \emph{On the groups {$H(\Pi,n)$}. {II}.
  {M}ethods of computation}, Ann. of Math. (2) \textbf{60} (1954), 49--139.
  \MR{65162}

\bibitem{fadell-neuwirth}
Edward Fadell and Lee Neuwirth, \emph{Configuration spaces}, Math. Scand.
  \textbf{10} (1962), 111--118. \MR{141126}

\bibitem{feigin-tsygan-additive}
B.~L. Fe\u{\i}gin and B.~L. Tsygan, \emph{Additive {$K$}-theory}, {$K$}-theory,
  arithmetic and geometry ({M}oscow, 1984--1986), Lecture Notes in Math., vol.
  1289, Springer, Berlin, 1987, pp.~67--209. \MR{923136}

\bibitem{fu-inf}
Jiaqi Fu, \emph{A duality between {L}ie algebroids and infinitesimal
  foliations}, arXiv preprint arXiv:2410.04950 (2024).

\bibitem{gepner-haugseng-kock}
David Gepner, Rune Haugseng, and Joachim Kock, \emph{{$\infty$}-operads as
  analytic monads}, Int. Math. Res. Not. IMRN (2022), no.~16, 12516--12624.
  \MR{4466007}

\bibitem{gillet-soule}
H.~Gillet and C.~Soul\'{e}, \emph{Intersection theory using {A}dams
  operations}, Invent. Math. \textbf{90} (1987), no.~2, 243--277. \MR{910201}

\bibitem{glasman}
Saul Glasman, \emph{Day convolution for {$\infty$}-categories}, Math. Res.
  Lett. \textbf{23} (2016), no.~5, 1369--1385. \MR{3601070}

\bibitem{goerss_f2}
Paul~G. Goerss, \emph{On the {A}ndr\'e-{Q}uillen cohomology of commutative
  {$\mathbf{F}_2$}-algebras}, Ast\'erisque (1990), no.~186, 169. \MR{1089001}

\bibitem{goodwillie-calc2}
Thomas~G. Goodwillie, \emph{Calculus. {II}. {A}nalytic functors}, $K$-Theory
  \textbf{5} (1991/92), no.~4, 295--332. \MR{1162445}

\bibitem{grothendieck-crystals}
A.~Grothendieck, \emph{Crystals and the de {R}ham cohomology of schemes}, Dix
  expos\'{e}s sur la cohomologie des sch\'{e}mas, Adv. Stud. Pure Math.,
  vol.~3, North-Holland, Amsterdam, 1968, Notes by I. Coates and O. Jussila,
  pp.~306--358. \MR{269663}

\bibitem{hedenlund-moulinos}
Alice Hedenlund and Tasos Moulinos, \emph{The synthetic {H}ilbert additive
  group scheme}, Selecta Mathematica \textbf{31} (2025), no.~5, 1--64.

\bibitem{hesselholt-pstragowski-i}
Lars Hesselholt and Piotr Pstr\k{a}gowski, \emph{Dirac geometry {I}:
  commutative algebra}, Peking Math Journal \textbf{8} (2025), 405--480.

\bibitem{holeman-derived}
Adam Holeman, \emph{Derived {$\delta$}-rings and relative prismatic
  cohomology}, arXiv preprint arXiv:2303.17447 (2023).

\bibitem{horel}
Geoffroy Horel, \emph{Binomial rings and homotopy theory}, J. Reine Angew.
  Math. \textbf{813} (2024), 283--305. \MR{4780987}

\bibitem{huneke}
Craig Huneke, \emph{On the symmetric and {R}ees algebra of an ideal generated
  by a {$d$}-sequence}, J. Algebra \textbf{62} (1980), no.~2, 268--275.
  \MR{563225}

\bibitem{illusie-cotangent-1}
Luc Illusie, \emph{Complexe cotangent et d\'{e}formations. {I}}, Lecture Notes
  in Mathematics, Vol. 239, Springer-Verlag, Berlin-New York, 1971.
  \MR{0491680}

\bibitem{jantzen}
Jens~Carsten Jantzen, \emph{Representations of algebraic groups}, Pure and
  Applied Mathematics, vol. 131, Academic Press, Inc., Boston, MA, 1987.
  \MR{899071}

\bibitem{johnson-mccarthy}
Brenda Johnson and Randy McCarthy, \emph{Taylor towers for functors of additive
  categories}, J. Pure Appl. Algebra \textbf{137} (1999), no.~3, 253--284.
  \MR{1685140}

\bibitem{kubrak-shuklin-zakharov}
Dmitry Kubrak, Georgii Shuklin, and Alexander Zakharov, \emph{Derived binomial
  rings {I}: integral {B}etti cohomology of log schemes}, arXiv preprint
  arXiv:2308.01110 (2023).

\bibitem{loday}
Jean-Louis Loday, \emph{Cyclic homology}, second ed., Grundlehren der
  Mathematischen Wissenschaften, vol. 301, Springer-Verlag, Berlin, 1998,
  Appendix E by Mar\'{\i}a O. Ronco, Chapter 13 by the author in collaboration
  with Teimuraz Pirashvili. \MR{1600246}

\bibitem{lurie-rotation}
Jacob Lurie, \emph{Rotation invariance in algebraic {$K$}-theory}, available at
  \url{https://www.math.ias.edu/~lurie/papers/Waldhaus.pdf}.

\bibitem{htt}
\bysame, \emph{Higher topos theory}, Annals of Mathematics Studies, vol. 170,
  Princeton University Press, Princeton, NJ, 2009. \MR{2522659}

\bibitem{ha}
\bysame, \emph{Higher algebra}, available at
  \url{https://www.math.ias.edu/~lurie/papers/HA.pdf}, version dated 18
  September 2017.

\bibitem{sag}
\bysame, \emph{Spectral algebraic geometry}, available at
  \url{https://www.math.ias.edu/~lurie/papers/SAG-rootfile.pdf}, version dated
  3 February 2018.

\bibitem{magidson-divided}
Kirill Magidson, \emph{Divided powers and derived de {R}ham cohomology}, arXiv
  preprint arXiv:2405.05153 (2024).

\bibitem{mao-crystalline}
Zhouhang Mao, \emph{Revisiting derived crystalline cohomology}, Bull. Soc.
  Math. France \textbf{152} (2024), no.~4, 659--784. \MR{4851406}

\bibitem{mcclure-schwanzl-vogt}
J.~McClure, R.~Schw\"anzl, and R.~Vogt, \emph{{$THH(R)\cong R\otimes S^1$} for
  {$E_\infty$} ring spectra}, J. Pure Appl. Algebra \textbf{121} (1997), no.~2,
  137--159. \MR{1473888}

\bibitem{moulinos-cartier}
Tasos Moulinos, \emph{Filtered formal groups, {C}artier duality, and derived
  algebraic geometry}, \'Epijournal G\'eom. Alg\'ebrique \textbf{8} (2024),
  Art. 2, 41. \MR{4717400}

\bibitem{mrt}
Tasos Moulinos, Marco Robalo, and Bertrand To\"en, \emph{A universal
  {H}ochschild-{K}ostant-{R}osenberg theorem}, Geom. Topol. \textbf{26} (2022),
  no.~2, 777--874. \MR{4444269}

\bibitem{mundinger}
Joshua Mundinger, \emph{On the differentials of the
  {H}ochschild--{K}ostant--{R}osenberg spectral sequence}, arXiv preprint
  arXiv:2410.01894 (2024).

\bibitem{nakaoka}
Minoru Nakaoka, \emph{Cohomology mod {$p$} of symmetric products of spheres.
  {II}}, J. Inst. Polytech. Osaka City Univ. Ser. A \textbf{10} (1959), 67--89.
  \MR{121789}

\bibitem{ogus_infinitesimal}
Arthur Ogus, \emph{Cohomology of the infinitesimal site}, Ann. Sci. \'{E}cole
  Norm. Sup. (4) \textbf{8} (1975), no.~3, 295--318. \MR{422280}

\bibitem{priddy}
Stewart Priddy, \emph{Mod {$p$} right derived functor algebras of the symmetric
  algebra functor}, J. Pure Appl. Algebra \textbf{3} (1973), 337--356.
  \MR{342592}

\bibitem{raksit}
Arpon Raksit, \emph{Hochschild homology and the derived de {R}ham complex
  revisited}, arXiv preprint arXiv:2007.02576 (2020).

\bibitem{toen-affines}
Bertrand To\"{e}n, \emph{Champs affines}, Selecta Math. (N.S.) \textbf{12}
  (2006), no.~1, 39--135. \MR{2244263}

\bibitem{toen-vezzosi-simpliciales}
Bertrand To{\"e}n and Gabriele Vezzosi, \emph{Alg\`ebres simpliciales
  {$S^1$}-\'{e}quivariantes, th\'{e}orie de de {R}ham et th\'{e}or\`emes {HKR}
  multiplicatifs}, Compos. Math. \textbf{147} (2011), no.~6, 1979--2000.
  \MR{2862069}

\bibitem{tv-book}
\bysame, \emph{Derived foliations}, arXiv preprint arXiv:2305.08212 (2023).

\bibitem{tv-inf}
\bysame, \emph{Infinitesimal derived foliations}, arXiv preprint
  arXiv:2305.13010 (2023).

\end{thebibliography}
